\documentclass[11pt]{article}
\usepackage[ansinew]{inputenc}
\usepackage{graphicx}
\usepackage{color}

\usepackage{enumerate,latexsym}
\usepackage{latexsym}
\usepackage{amsmath,amssymb}% CUIDADO: Si uso estos packages, $\widetilde{}$ no funciona a veces
\usepackage{graphicx}
\newfont{\bb}{msbm10 at 11pt}
\newfont{\bbsmall}{msbm8 at 8pt}
\def\rth{\mathbb{R}^3}
\def\R{\mathbb{R}}

\def\B{\mathbb{B}}
\def\N{\mathbb{N}}
\def\Z{\mathbb{Z}}
\def\C{\mathbb{C}}

\def\D{\mathbb{D}}

\def\esf{\mathbb{S}}

\newcommand{\re}{\mbox{\bb R}}

\newcommand{\ben}{\begin{enumerate}}
\newcommand{\bit}{\begin{itemize}}
\newcommand{\een}{\end{enumerate}}
\newcommand{\eit}{\end{itemize}}

\newcommand{\Int}{\mbox{Int}}

\newcommand{\Div}{\mbox{\rm div}}

\def\a{{\alpha}}
\def\lc{{\cal L}}
\def\t{{\theta}}

\def\g{{\gamma}}
\def\G{{\Gamma}}
\def\l{{\lambda}}

\def\de{{\delta}}
\def\be{{\beta}}
\def\ve{{\varepsilon}}

\def\centerbmp#1#2#3{\vskip#2\relax\centerline{\hbox to#1{\special
    {bmp:#3 x=#1, y=#2}\hfil}}}

\newtheorem{theorem}{Theorem}[section]
\newtheorem{lemma}[theorem]{Lemma}
\newtheorem{proposition}[theorem]{Proposition}
\newtheorem{thm}[theorem]{Theorem}
\newtheorem{remark}[theorem]{Remark}
\newtheorem{corollary}[theorem]{Corollary}
\newtheorem{definition}[theorem]{Definition}
\newtheorem{conjecture}[theorem]{Conjecture}
\newtheorem{assertion}[theorem]{Assertion}

\newtheorem{claim}[theorem]{Claim}

\newenvironment{proof}{\smallskip\noindent{\it Proof.}\hskip \labelsep}
{\hfill\penalty10000\raisebox{-.09em}{$\Box$}\par\medskip}

\begin{document}

\begin{title}
{Properly embedded minimal planar domains with infinite topology are
Riemann minimal examples}
\end{title}
\vskip .5in

\begin{author}
{William H. Meeks III\thanks{This material is based upon
   work for the NSF under Awards No. DMS -
   0405836 and DMS - 0703213. Any opinions, findings, and conclusions or recommendations
   expressed in this publication are those of the authors and do not
   necessarily reflect the views of the NSF.}
   \and Joaqu\'\i n P\' erez
\thanks{Research partially supported by a MEC/FEDER
grant no. MTM2007-61775 and a Junta de Andaluc\'\i a grant no. P06-FQM-01642.}, }
\end{author}
\maketitle

\nocite{me21,gny1,dhkw1}

\begin{abstract}
These notes outline recent developments in classical minimal
surface theory that are essential in classifying the properly
embedded minimal planar domains $M\subset \rth$ with infinite
topology (equivalently, with an infinite number of ends). This final
classification result by Meeks, P\'{e}rez, and Ros \cite{mpr6}
states that such an $M$ must be congruent to a homothetic scaling of
one of the classical examples found by Riemann \cite{ri2} in 1860.
These examples ${\cal R}_s, 0<s<\infty $, are defined in terms of the
Weierstrass ${\cal P}$-functions ${\cal P}_t$ on the rectangular
elliptic curve $\frac{\C}{\langle 1, t\sqrt{-1}\rangle }$, are
singly-periodic and intersect each horizontal plane in $\rth$ in a
circle  or a line parallel to the $x$-axis. Earlier work by Collin
\cite{col1}, L\'{o}pez and Ros \cite{lor1} and Meeks and Rosenberg
\cite{mr8} demonstrate that the plane, the catenoid and the helicoid
are the only properly embedded minimal surfaces of genus zero with
finite topology (equivalently, with a finite number of ends). Since the
surfaces ${\cal R}_s$ converge to a catenoid as $s\to 0$ and to  a
helicoid as $s\to \infty$, then the moduli space ${\cal M}$ of all
properly embedded, non-planar, minimal planar domains in $\rth$ is
homeomorphic to the closed unit interval $[0,1]$.
\par
\vspace{.1cm}

%\vspace{.1cm} \noindent{\it Key words and phrases:} Minimal surface,
% minimal lamination, locally simply-connected, finite total curvature,
% conformal structure, harmonic
%function, recurrence, transience,
% parabolic Riemann surface, harmonic measure,
%universal superharmonic function, Jacobi function, stability, index
%of stability, Shiffman function, Korteweg-de Vries equation, KdV
%hierarchy, algebro-geometric potential, curvature estimates, maximum
%principle at infinity, limit tangent cone at infinity, limit tangent
%plane at infinity, parking garage, local picture on the scale of topology,
%minimal planar domain, isotopy class.
\end{abstract}
% \pagebreak

\tableofcontents

\section{Introduction.}

In the last decade spectacular progress has been made in various
aspects of classical minimal surface theory. Some of the successes
obtained are the solutions of open problems which have been pursued
since the birth of this subject in the 19-th century, while others
have opened vast new horizons for future research. Among the first
such successes, we would like to highlight the achievement of a
deep understanding of topological aspects of proper minimal
embeddings in three-space including their complete topological
classification~\cite{ckmr1,fr2,fme3,fme2,fme1,fme5,fme4}. Equally
important in this progress has been a comprehensive analysis of the
behavior of limits of sequences of embedded minimal surfaces without
{\it a priori} area or curvature
bounds~\cite{cm25,cm21,cm22,cm24,cm23}, with outstanding
applications such as the classification of all simply-connected,
properly embedded minimal surfaces~\cite{mr8}. Also, many deep
results have been obtained on the subtle relationship between
completeness and properness for complete immersed minimal
surfaces~\cite{afm1,fmm1,mmn,marmor1,marmor2,mn1,na1}, and how
embeddedness introduces a strong dichotomy in this
relationship~\cite{cm35,mpr9}. While all of these results are
extremely interesting, they will not be treated in these notes (at
least, not in depth) but they do give an idea of the enormous
activity within this field; instead, we will explain the recent
solution to the following long standing problem in classical minimal
surface theory:
\begin{quote}
{\it Classify all
possible properly embedded minimal surfaces of genus zero in $\R^3$.}
\end{quote}
Research by various authors help to understand this problem, the
more relevant work being by Colding and Minicozzi~\cite{cm23,cm35},
Collin~\cite{col1}, L\'opez and Ros~\cite{lor1}, Meeks, P\'erez and
Ros~\cite{mpr6} and Meeks and Rosenberg~\cite{mr6,mr8}. In fact,
this problem has been one of main goals of the two authors of these
notes (in collaboration with A. Ros) for over the past 15 years, and
the long path towards its solution has been marked by the discovery
of powerful techniques which have proved useful in other
applications. Putting together all of these efforts, we now state
the final solution to the above problem, whose proof appears in
\cite{mpr6}.

\begin{theorem}%[Classification Theorem for Minimal Planar Domains]
\label{classthm} Up to scaling and rigid motion, any connected,
properly embedded, minimal planar domain in $\rth$ is a plane, a
helicoid, a catenoid or one of the Riemann minimal
examples\footnote{See Section \ref{subsecexamples} for further
discussion of these surfaces.}. In particular, for every such
surface there exists a foliation of $\R^3$ by parallel planes, where
each plane intersects the surface transversely in a circle or a
straight line.
\end{theorem}

In these notes we will try to pass on to the reader a glimpse of the
beauty of the arguments and different theories that intervene in the
proof of Theorem~\ref{classthm}. Among these auxiliary theories, we
highlight the theory of integrable systems, whose applications to
minimal and constant mean curvature surface theory have gone far
beyond the existence results in the late eighties
(Abresch~\cite{ab1}, Bobenko~\cite{bob4}, Pinkall and
Sterling~\cite{ps1}, based on the sinh-Gordon equation) to recent
uniqueness theorems like the one that gives the title of these notes
(based on the KdV equation), and the even more recent tentative
solution to the Lawson conjecture by Kilian and
Schmidt~\cite{KiSch1}.

Before proceeding, we make a few general comments about the
organization of this article, which relate to the proof of
Theorem~\ref{classthm}. In Section~\ref{Background} we
briefly  introduce the main definitions and background material, including
short discussions of the classical examples that appear in the statement of
the above theorem. Since a complete, immersed minimal surface $M$
without boundary in $\R^3$ cannot be compact,
$M$ must have {\it ends}\footnote{An {\it end} of a non-compact connected
topological manifold $M$
is an equivalence class in the set ${\mathcal A}=\{ \a
\colon [0,\infty )\to M\ | \ \a \mbox{
 is a proper arc}\} $, under the equivalence relation:
 $\a _1\sim \a _2$ if for every compact set $C \subset M$,
 $\a _1,\a _2$ lie eventually in the same component of $M-C$.
If $\a \in {\mathcal A}$ is a representative proper arc of and end
of $M$ and $\Omega \subset M$ is a proper subdomain with compact
boundary such that $\alpha \subset \Omega $, then we say that
$\Omega $ {\it represents} the end $e$.}. As we are interested in
planar domains, the allowed topology for our surfaces in that of the
two-sphere minus a compact totally disconnected set ${\cal E}\neq
\mbox{\O }$ corresponding to the space of ends of the surface. A
crucial result proved by Collin~\cite{col1} in 1997 (see
Conjecture~\ref{conjNitsche} below) implies that when the
cardinality $\# ({\cal E})$ of ${\cal E}$ satisfies $2\leq \# ({\cal
E})<\infty $, then the total Gaussian curvature of $M$ is finite.
Complete embedded minimal surfaces with finite total curvature
comprise the best understood minimal surfaces in $\R^3$; the main
reason for this is the fact discovered by Osserman that the
underlying complex structure for every such minimal surface $M$ is
that of a compact Riemann $\overline{M}$ surface minus a finite
number of points, and the classical analytic Weierstrass
representation data on $M$ extends  across the punctures to
meromorphic data on $\overline{M}$. Using the maximum principle
together with their result that every complete, embedded minimal
surface with genus zero and finite total curvature can be minimally
deformed, L\'opez and Ros characterized the plane and the catenoid
as being the unique embedded examples of finite total curvature and
genus zero, see Theorem~\ref{thmLR}. We also explain Collin's and
L\'opez-Ros' theorems in Section~\ref{ftc}.

Section~\ref{sec1conn} covers the one-ended case of
Theorem~\ref{classthm}, which was solved by Meeks and Rosenberg
(Theorem~\ref{ttmr}). To understand the proof of this result we need
the local results of Colding and Minicozzi which describe the
structure of compact, embedded minimal disks as essentially being
modeled by either a plane or a helicoid, and their one-sided
curvature estimate together with other results of a global nature
such as their limit lamination theorem for disks, see
Theorem~\ref{thmlimitlaminCM} below.

The remainder of the article, except for the last section, focuses on the case in
Theorem~\ref{classthm} where the surface has infinitely many ends.
Crucial in this discussion is the Ordering Theorem by Frohman and
Meeks (Theorem~\ref{thordering}) as well as two topological results,
the first on the non-existence of middle limit ends due to Collin,
Kusner, Meeks and Rosenberg (Theorem~\ref{thmckmr}) and the second
on the non-existence of properly embedded minimal planar domains
with just one limit end by Meeks, P\'erez and Ros
(Theorem~\ref{thmno1limitend}). It follows from these two
non-existence results that a properly embedded minimal planar domain
must have exactly two limit ends. These key ingredients represent
the content of Section~\ref{reduction2limitends}.

As a consequence of the results in Section \ref{reduction2limitends},
in Section \ref{sec2limitends} we obtain strong control on the
conformal structure and height differential of a properly embedded,
minimal planar domain $M$ with infinitely many ends. In this section
we also explain how Colding-Minicozzi theory can be applied to
obtain a curvature bound for such an $M$, after it has been the
normalized by a homothety so that its vertical flux is $1$,
which only depends on
the length of the horizontal component of its flux
vector (Theorem~\ref{thm1}). This bound leads to a quasi-periodicity
property of $M$, which is a cornerstone to finishing the
classification problem. Also in Section~\ref{sec2limitends} we
introduce the Shiffman Jacobi function $S_M$ and explain how the
existence of a related holomorphic deformation of $M$ preserving its
flux is sufficient to reduce the proof of Theorem \ref{classthm} in
the case with infinitely many ends to the singly-periodic case,
which was solved earlier by Meeks, P\'erez and Ros \cite{mpr1}.

In Section~\ref{secKdV} we explain how the existence of the desired
holomorphic deformation of $M$ follows from the integration of an
evolution equation for its Gauss map. The Shiffman function $S_M$
will be crucial at this point by enabling us to reduce this
evolution equation to an equation of type Korteweg de Vries (KdV). A
classical condition that implies global
 integrability of the KdV equation, i.e.
existence of (globally defined) meromorphic solutions $u(z,t)$ of the Cauchy problem associated
to the KdV equation, is that the initial condition $u(z)$
is an {\it algebro-geometric} meromorphic function, a concept
related to the hierarchy of the KdV equation. In our
setting, the final step of the classification of the properly
embedded minimal planar domains consists of proving that the initial
condition $u(z)$ for the Cauchy problem of the KdV equation
naturally associated to any quasiperiodic, possibly immersed,
minimal planar domain $M$ with two limit ends is algebro-geometric
(Section~\ref{subsec7.4}), which in turn is a consequence of the
fact that the space of bounded Jacobi functions on such a surface
$M$ is finite dimensional (Theorem~\ref{bounded}).

An important consequence of the proof of the classification of
properly embedded minimal planar domains is the characterization of the
asymptotic behavior of the ends of any properly embedded minimal
surface $M\subset \rth$ with finite genus. If $e\in {\cal E} (M)$ is
an end, then there exists properly embedded domain
$E(e)\subset M$ with compact boundary which represents $e$ and such that in a natural
sense $E(e)$ converges to the end of a plane, a catenoid, a helicoid
or to one of the limit ends of a Riemann minimal example. This asymptotic
characterization is due to Schoen~\cite{sc1} and Collin \cite{col1}
when $M$ has a finite number of ends greater
than one, to Meeks, P\'{e}rez and Ros in the case $M$ has an
infinite number of ends, and to Meeks and Rosenberg~\cite{mr8},
Bernstein and Breiner~\cite{bb2},
and Meeks and P\'erez~\cite{mpe3,mpe5} in the
case of just one end. These
asymptotic characterization results
will be explained  in Section \ref{secas} of these notes.

\section{Background.}
\label{Background}

An isometric immersion $X=(x_1,x_2,x_3)\colon M\to \R^3$ of a
Riemannian surface into Euclidean space is said to be {\it minimal}
if $x_i$ is a harmonic function on~$M$ for each $i$ (in the sequel,
we will identify the Riemannian surface $M$ with its image under the
isometric embedding). Since harmonicity is a local property, we can
extend the notion of minimality to immersed surfaces $M\subset
\R^3$. We will always assume all surfaces under consideration are
orientable. If $M\subset \R^3$ is an immersed oriented surface, we
will denote by $H$ the {\it mean curvature} function of $X$ (average
normal curvature) and by $N \colon M \to \esf^2(1) \subset \R^3$ its
Gauss map. Since $\Delta X= 2HN$ (here $\Delta $ is the
Laplace-Beltrami operator on $M$), we have that $M$ is minimal if
and only if $H=0$ identically. Expressing locally $M$ as the graph
of a function $u=u(x,y)$ (after a rotation), the last equation can
be written as the following quasilinear second order elliptic PDE:
\begin{equation}
\label{eq:ecminsurf}
\Div _0\left( \frac{\nabla _0u}{\sqrt{1+|\nabla _0u|^2}}\right) =0,
\end{equation}
where the subscript $\bullet _0$ indicates that the corresponding
object is computed with respect to the flat metric in the plane.

Let $\Omega $ be a subdomain with compact closure in a surface
$M\subset \R^3$ and let $u\in C^{\infty }_0(\Omega )$ be a compactly
supported smooth function. The first variation of the area functional
$A(t)=\mbox{Area}((X+tuN)(\Omega ))$ for the normally perturbed
immersions $X+utN$ (with $|t|$ sufficiently small) gives
\begin{equation}
\label{eq:1vararea} A'(0)=-2\int _{\Omega }uH\, dA,
\end{equation}
where $dA$ stands for the area element of $M$. Therefore,
$M$ is minimal if and only if it is
a critical point of the area functional for all compactly
supported variations. The second variation of area implies
that any point in a minimal surface has a
neighborhood with least-area relative to its boundary\footnote{This property
justifies the word ``minimal'' for these surfaces.}, and thus
$M\subset \R^3$ is minimal if and only
  if every point $p\in M$ has a
 neighborhood with least-area relative to its boundary.
 If one exchanges the area functional~$A$ by the {\it Dirichlet energy}
$E=\int _{\Omega }|\nabla X|^2 dA$, then
the two functionals are related by $E\geq 2A$, with equality
if and only if the immersion $X\colon M\to \R^3$ is conformal.
This relation between area and energy together with
the existence of isothermal coordinates on every Riemannian surface
allow us to state that a conformal immersion $X\colon M\to \R^3$ is
minimal if and only if it is a critical point of the Dirichlet
 energy for all compactly supported variations (or
 equivalently,
every point $p\in M$ has a neighborhood with least energy relative
to its boundary). Finally, the relation $A_p=-dN_p$ between the
differential of the Gauss map and the
{\it shape operator}%(a symmetric endomorphism of the tangent space
%$T_pM$ of $M$ at $p\in M$ whose orthogonal eigenvectors are
%{\it principal directions} and the corresponding eigenvalues
%are the {\it principal curvatures} of $M$ at $p$)
, together with the Cauchy-Riemann equations give that $M$ is {\it
minimal} if and only if its stereographically projected Gauss map
$g\colon M\to \C \cup \{ \infty \} $ is a holomorphic function. All
these equivalent definitions of minimality illustrate the wide
variety of branches of mathematics that appear in its study:
Differential Geometry, PDE, Calculus of Variations, Geometric
Measure Theory, Complex Analysis, etc.

The Gaussian curvature function\footnote{If needed, we will use the
notation $K_M$ to highlight the surface $M$ of which $K$ is the
Gaussian curvature function.} $K$ of a surface $M\subset \R^3$ can
be written as $K=\kappa _1\kappa _2=\det A$, where $\kappa _i$ are
the principal curvatures and $A$ the shape operator. Thus $|K|$ is
the absolute value of the Jacobian of the Gauss map~$N$. If $M$ is
minimal, then $\kappa _1=-\kappa _2$ and $K\leq 0$, hence the
{\it total curvature} $C(M)$ of $M$ is the
negative of the spherical area of $M$ through its Gauss map,
counting multiplicities:
\begin{equation}
\label{eq:curvtot} C(M)=\int_MK\, dA = -{\rm Area}(N\colon
M\rightarrow \esf ^2(1))\in [-\infty ,0].
\end{equation}
In the sequel, we will denote by $\B (p,R)=
\{ x\in \R^3 \ | \ |x-p|<R\} $ and $\B (R)=\B (\vec{0},R)$.

\subsection{Weierstrass representation and the definition of flux.}
Let $M\subset \R^3$ be a possibly immersed minimal surface, with
stereographically projected Gauss map $g\colon M\to \C \cup \{
\infty \} $. Since the third coordinate function $x_3$ of $M$ is
harmonic, it admits a locally well-defined harmonic conjugate
function $x_3^*$. The {\it height differential of $M$} is the
holomorphic 1-form $dh=dx_3+idx_3^*$ (note that $dh$ is not
necessarily exact on $M$). The minimal immersion $X\colon M\to \R^3$
can be written up to translation by the  vector $X(p_0)$, $p_0\in M$,
solely in terms of the {\it Weierstrass data} $(g,dh)$ as
\begin{equation}
\label{eq:repW}
 X(p)=\Re \int _{p_0}^p\left( \frac{1}{2}\left(
 \frac{1}{g}-g\right)
,\frac{i}{2}\left( \frac{1}{g}+g\right) ,1\right) dh,
\end{equation}
where $\Re $ stands for real part. The positive-definiteness of the
induced metric and the independence of (\ref{eq:repW}) with respect
to the integration path give rise to certain compatibility
conditions on the meromorphic data $(g,dh)$ for analytically
defining a minimal surface (Osserman~\cite{os3}); namely, if we
start with a meromorphic function $g$ and a holomorphic one-form
$dh$ on an abstract Riemann surface $M$, then the map $X\colon M\to
\R^3$ given by {\rm (\ref{eq:repW})} is a conformal minimal
immersion with Weierstrass data $(g,dh)$ provided that two
conditions hold:
\begin{equation}
\label{eqWeiers1}
\mbox{The zeros of $dh$ coincide with the poles and zeros of $g$,
with the same order.}
\end{equation}
\begin{equation}
\label{eqWeiers2}
\overline{\int _{\g }g\, dh}=\int _{\g }\frac{dh}{g},\quad
    \Re \int _{\g }dh=0 \ \mbox{ for any closed curve $\g \subset M$
    (period problem).}
\end{equation}
The {\it flux vector} of $M$ along a closed curve $\g \subset M$ is defined as
\begin{equation}
\label{eq:flux}
F( \gamma) = \int_\gamma \mbox{Rot}_{90^\circ }(\g ')
= \Im \int _{\g }\left( \frac{1}{2}\left(
\frac{1}{g}-g\right), \frac{i}{2}\left( \frac{1}{g}+g\right)
,1\right) dh\in \R^3,
\end{equation}
where $\mbox{Rot}_{90^\circ }$ denotes the counterclockwise rotation by angle $\pi /2$ in
 the tangent plane of $M$ at any point, and $\Im $ stands for imaginary
  part. Both the period condition (\ref{eqWeiers2}) and the flux vector (\ref{eq:flux})
only depend on the homology class of $\g $ in $M$.

\subsection{Maximum principles.}
Since minimal surfaces can be written locally as solutions of the
PDE~(\ref{eq:ecminsurf}),
they satisfy certain maximum principles.

\begin{theorem}[Interior and boundary maximum principles~\cite{sc1}]
\label{thmintmaxprin}
 Let $M_1,M_2$ be connected minimal surfaces in $\R^3$ and
 $p$ an interior point to both surfaces, such that $T_pM_1=T_pM_2=
 \{ x_3=0\} $. If
$M_1,M_2$ are locally expressed as the graphs of functions $u_1,u_2$
around $p$ and $u_1\leq u_2$ in a neighborhood of $p$, then
$M_1=M_2$ in a neighborhood of~$p$. The same conclusion holds if
$p$ is a boundary point of both surfaces and additionally,
$T_p\partial M_1=T_p\partial M_2$.
\end{theorem}
We also dispose of more sophisticated versions of the maximum
principle, where a first contact point of two minimal surfaces does
not occur at a finite point but at infinity, amongst which we state
two. The first one (whose proof we sketch for later purposes) was
proved by Hoffman and Meeks, and the second one is due to Meeks and
Rosenberg.

\begin{theorem}[Half-space Theorem~\cite{hm10}]
\label{thmhalf}
A properly immersed, non-planar minimal surface
without boundary cannot be contained in a half-space.
\end{theorem}
{\it Sketch of proof.} Arguing by contradiction, suppose that a
surface $M\subset \R^3$ as in the hypotheses is contained in $\{
x_3\geq 0\} $ (and so, $M\subset \{ x_3>0\} $ by
Theorem~\ref{thmintmaxprin}) but is not contained in $\{ x_3>c\} $
for any $c>0$. Since $M$ is proper, we can find a ball $\B (p,r)$
centered at a point $p\in \{ x_3=0\} $ such that $M\cap \B
(p,r)=\mbox{\O }$. Consider a vertical half-catenoid $C$ with
negative logarithmic growth, completely contained in $\{ x_3\leq 0\}
\cup \B (p,r)$, whose waist circle $\G $ is centered at $q= p+\ve
(0,0,1)$, $\ve >0$ small. Then $M\cap C=\mbox{\O }$. Now deform
$C=C(1)$ by a one-parameter family of non-compact annular pieces of
vertical catenoids $\{ C(r)\} _{r\in (0,1]}$, all having the same
boundary $\G $ as $C$, with negative logarithmic growths converging
to zero as $r\to 0$ and whose Gaussian curvatures blow up at a the
limit point of the waist circles of $C(r)$, which is the point $q$.
Then, the surfaces $C(r)$ converge on compact subsets of $\R^3-\{
q\} $ to the plane $\{ x_3=\ve \} $, and so, $M$ achieves a first
contact point with one of the catenoids, say $C(r_0)$, in this
family; the usual maximum principle for $M$ and $C(r_0)$ gives a
contradiction.
{\hfill\penalty10000\raisebox{-.09em}{$\Box$}\par\medskip}

\begin{theorem}[Maximum Principle at Infinity~\cite{mr7}]
\label{thmMaxprin}
Let $M_1,M_2\subset \R^3$ be disjoint, connected,
properly immersed minimal surfaces
with (possibly empty) boundary.
\begin{itemize}
\item[{\it i)}] If $\partial M_1\neq \mbox{\rm \O }$ or
$\partial M_2\neq \mbox{\rm \O }$, then after possibly re-indexing,
\[
{\rm dist} (M_1,M_2) = \inf \{ {\rm dist} (p, q) \mid
p \in \partial M_1, \, q \in M_2 \} .
\]
\item[{\it ii)}] If $\partial M_1=\partial M_2=\mbox{\rm \O}$,
 then $M_1$ and $M_2$
are flat.
\end{itemize}
\end{theorem}
\subsection{Monotonicity formula.}
Monotonicity formulas, as well as maximum principles, play a crucial
role in many areas of Geometric Analysis. For instance, we will see
in Theorem~\ref{thmckmr} how the monotonicity formula can be used to
discard middle limit ends for a properly embedded minimal surface.
We will state this basic result without proof here; it is a
consequence of the classical coarea formula applied to the distance
function to a point $p\in \R^3$, see for instance Corollary 4.2
in~\cite{cm34} for a detailed proof.
\begin{theorem}[Monotonicity Formula~\cite{cmCourant,kks1}]
\label{monotoform} Let $X\colon M\to \R^3$ be a connected, properly
immersed minimal surface. Given $p\in \R^3 $, let $A(R)=\mbox{\rm
Area}(X(M)\cap \B (p,R))$. Then, $A(R)R^{-2}$ is non-decreasing. In
particular,  $\lim _{R\to \infty}A(R)R^{-2} \geq \pi$ with equality
if and only if $M$ is a plane.
\end{theorem}

\subsection{Stability, Plateau problem and barrier constructions.}
\label{subsecbarrier} Recall that every (orientable) minimal surface
$M\subset \R^3$ is a critical point of the area functional for
compactly supported normal variations. $M$ is said to be {\it
stable} if it is a local minimum for such a variational problem. The
following well-known result indicates how restrictive is stability
for complete minimal surfaces. It was proved independently by
Fischer-Colbrie and Schoen~\cite{fs1}, do Carmo and Peng~\cite{cp1},
and Pogorelov \cite{po1} for orientable surfaces and more recently
by Ros \cite{ros9} in the case of non-orientable surfaces.
\begin{theorem}
\label{thmstablecompleteplane}
If $M\subset\R^3$ is a complete, immersed, stable minimal surface, then $M$ is a plane.
\end{theorem}
The {\it Plateau Problem} consists of finding a {\it compact}
surface of least area spanning a given boundary. This problem can be
solved under certain circumstances; this existence of compact
solutions together with a taking limits procedure leads to construct
non-compact {\it stable} minimal surfaces in $\rth$ that are
constrained to lie in regions of space whose boundaries have
non-negative mean curvature.  A huge amount of literature is devoted
to this procedure, but we will state here only a particular version.

Let $W$ be a compact Riemannian three-manifold with boundary which embeds in the
interior of another Riemannian three-manifold. $W$ is said to have
{\it piecewise smooth, mean convex boundary} if $\partial W$ is a
two-dimensional complex consisting of a finite number of smooth,
two-dimensional compact simplices with interior angles less than or
equal to $\pi $, each one with non-negative mean curvature with
respect to the inward pointing normal. In this situation, the
boundary of $W$ is a good barrier for solving Plateau problems in
the following sense:

\begin{theorem}[\cite{Me2,my1,my2,si1}]
 \label{Plateau}
Let $W$ be a compact Riemannian three-manifold with piecewise smooth
mean convex boundary. Let $\G $ be a smooth collection of pairwise
disjoint closed curves in $\partial W$, which bounds a compact
orientable surface in $W$. Then, there exists an embedded orientable
surface $\Sigma \subset W$ with $\partial \Sigma =\G $ that
minimizes area among all orientable surfaces with the same boundary.
\end{theorem}

Instead of giving an idea of the proof of Theorem~\ref{Plateau}, we
will illustrate it together with a taking limits procedure to obtain
a particular result in the non-compact setting. Consider two
disjoint, connected, properly embedded minimal surfaces $M_1, M_2$
in $\rth$ and let $W$ be the closed complement of $M_1\cup M_2$ in
$\rth$ that has both $M_1$ and $M_2$ on its boundary.
\begin{enumerate}
\item We first show how to produce compact, least-area surfaces in $W$ with prescribed
boundary lying in the boundary $\partial W$. Note that $W$ is a
complete flat three-manifold with boundary, and $\partial W$ has
mean curvature zero.  Meeks and Yau~\cite{my2} proved that $W$
embeds isometrically in a homogeneously regular\footnote{A
Riemannian three-manifold $N$ is {\it homogeneously regular}
if given $\ve >0$ there exists $\de > 0$ such that $\de$-balls in $N$ are $\ve $-uniformly
close to $\de$-balls in $\rth$ in the $C^2$-norm. In particular, if
$N$ is compact, then $N$ is homogeneously regular.} Riemannian
three-manifold $\widetilde{W}$ diffeomorphic to the interior of $W$
and with metric $\widetilde{g}$.  Morrey~\cite{mor1} proved that in
a homogeneously regular manifold, one can solve the classical
Plateau problem.  In particular, if $\Gamma$ is an embedded
$1$-cycle in $\widetilde{W}$ which bounds an orientable finite sum
of differentiable simplices in $\widetilde{W}$, then by standard
results in geometric measure theory~\cite{fed1}, $\Gamma$ is the
boundary of a compact, least-area embedded surface $\Sigma_\Gamma
(\widetilde{g})\subset \widetilde{W}$. Meeks and Yau proved that the
metric $\widetilde{g}$ on $\widetilde{W}$ can be approximated by a
family of homogeneously regular metrics $\{ g_n\}_{n \in
\mathbb{N}}$ on $\widetilde{W}$, which converges smoothly on compact
subsets of $\widetilde{W}$   to $\widetilde{g}$, and each $g_n$
satisfies a convexity condition outside of $W \subset
\widetilde{W}$, which forces the least-area surface
$\Sigma_{\Gamma}(g_n)$ to lie in $W$ if $\Gamma$ lies in $W$.  A
subsequence of the $\Sigma_{\Gamma}(g_n)$ converges to a smooth
minimal surface $\Sigma_\Gamma$ of least-area in $W$ with respect to
the original flat metric, thereby finishing our description of how
to solve (compact) Plateau-type problems in $W$.
\item We now describe the limit procedure to construct a non-compact, stable minimal
surface with prescribed boundary lying in $M_1$.
Let $M_1(1) \subset \hdots \subset M_1(n) \subset \hdots$ be a
compact exhaustion of $M_1$, and let $\Sigma_1(n)$ be
a least-area surface in $W$ with boundary $\partial M_1(n)$, constructed
as in the last paragraph. Let $\alpha$ be a compact arc in $W$ which joins
a point in  $M_1(1)$ to a point in $\partial W \cap M_2$.  By
elementary intersection theory, $\alpha$ intersects every least-area
surface $\Sigma_1(n)$.  By compactness of least-area surfaces, a
subsequence of the surfaces $\Sigma_1(n)$ converges to a
properly embedded area-minimizing surface $\Sigma$ in $W$
with a component $\Sigma_0$ which intersects $\alpha$;
this proof of the existence of $\Sigma$ is due to Meeks, Simon and
Yau~\cite{msy1}.
\item We finish this application of the barrier construction,
as follows. Since $\Sigma_0$ separates $\rth$, $\Sigma_0$ is orientable and
so, Theorem~\ref{thmstablecompleteplane} insures that $\Sigma_0$
is a plane.  Hence, $M_1$ and $M_2$ lie in closed half-spaces of $\rth$,
and the Half-space Theorem (Theorem~\ref{thmhalf}) implies that both
$M_1,M_2$ are planes.
\end{enumerate}
The above items 1--3 give the following generalization of Theorem~\ref{thmhalf},
also due to Hoffman and Meeks.
\begin{theorem}[Strong Half-space Theorem~\cite{hm10}]
If $M_1$ and $M_2$ are two disjoint, properly immersed minimal surfaces
in $\rth$, then $M_1$ and $M_2$ are parallel planes.
\end{theorem}

\subsection{The examples that appear in Theorem~\ref{classthm}.}
\label{subsecexamples}
We will now use
the Weierstrass representation for introducing the minimal planar
domains characterized in Theorem~\ref{classthm}. \vspace{.15cm}
\par
\noindent{\bf  The plane. }$M=\C $, $g(z)=1$, $dh=dz$. It is the
only complete, flat minimal surface in $\R^3$. \vspace{.15cm}
\par
\begin{figure}
\begin{center}
\includegraphics[width=16cm]{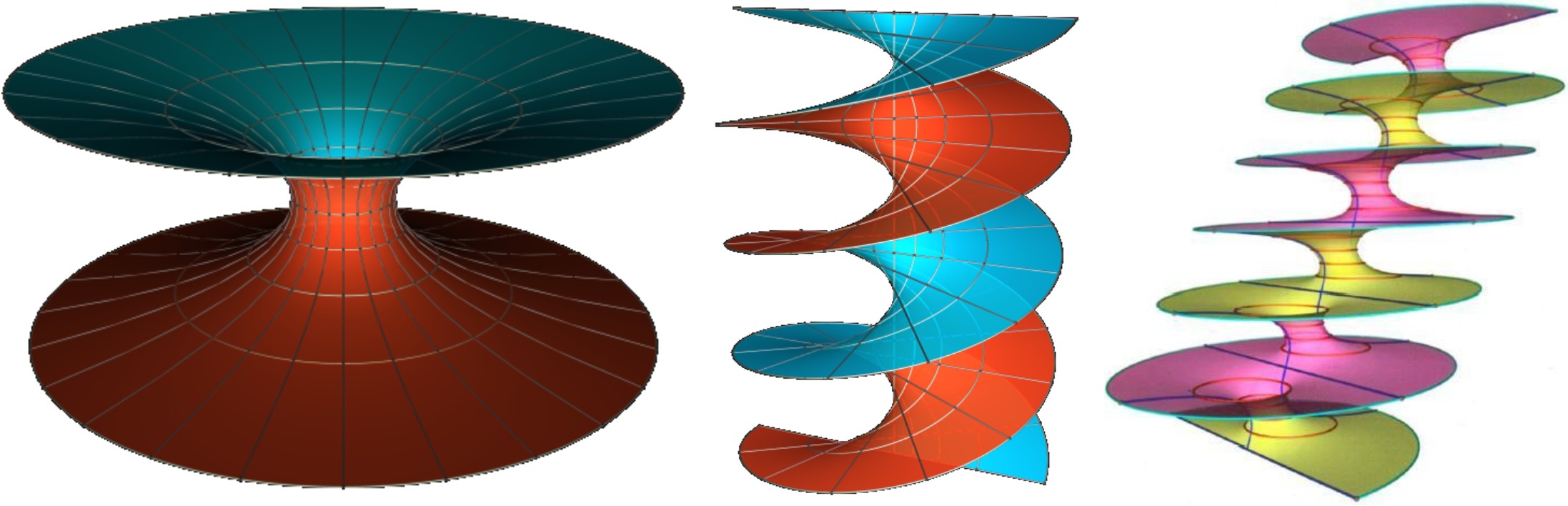}
\caption{Left: Catenoid. Center: Helicoid. Right: One of the Riemann minimal
examples. Figures courtesy of Matthias Weber.}
\label{cathelRiem}
\end{center}
\end{figure}

\noindent
{\bf The catenoid.} $M=\C -\{ 0\} $, $g(z)=z$,
 $dh=\frac{dz}{z}$ (Figure~\ref{cathelRiem} left).  This surface has
genus zero, two ends and total curvature $-4\pi $. Together with the
plane, the catenoid is the only minimal surface of revolution
(Bonnet~\cite{Bonnet2}). As we will see in Theorem~\ref{thmLR}, the
catenoid and the plane are the unique complete, embedded minimal
surfaces with genus zero and finite total curvature. Also, the
catenoid was characterized by Schoen~\cite{sc1} as being the unique
complete, immersed minimal surface with finite total curvature and
two embedded ends.

\par
\vspace{.2cm}
\noindent
{\bf The helicoid.} $M=\C $, $g(z)=e^z$,
 $dh=i\, dz$ (Figure~\ref{cathelRiem} center).
When viewed in $\R^3$, the helicoid has genus zero, one end and
infinite total curvature. Together with the plane, the helicoid is
the only ruled minimal surface (Catalan~\cite{catalan1}), and we
will see in Theorem~\ref{ttmr} below that it is the unique properly
embedded, simply-connected minimal surface. The vertical helicoid is
invariant by a vertical translation $T$ and by a 1-parameter family
of screw motions\footnote{A {\it screw motion} $S_{\t }$ is the
composition of a rotation of angle $\t $ around the $x_3$-axis with
a translation in the direction of this axis.} $S_{\theta }$, $\theta
>0$. Viewed in $\R^3/T$ or in $\R^3/S_{\t }$, the helicoid is a
properly embedded minimal surface with genus zero, two ends and
finite total curvature. The catenoid and the helicoid are {\it
conjugate} minimal surfaces, in the sense that the coordinate
functions of one of these surfaces are the harmonic conjugates of
the coordinate functions of the other one; in this case, we consider
the catenoid to be defined on its universal cover $e^z \colon
\mathbb{C} \rightarrow \mathbb{C}-\{0\}$ in order for the harmonic
conjugate of $x_3$ to be well-defined. \vspace{.15cm}
\par
\noindent {\bf The Riemann minimal examples.} They form a
one-parameter family, with Weierstrass data $M_\lambda =\{ (z,w)\in
(\C\cup \{ \infty \} )^2\ | \ w^2=z(z-\lambda )(\lambda z+1)\} -\{
(0,0),(\infty ,\infty )\} $, $g(z,w)=z$, $dh=A_{\lambda
}\frac{dz}{w}$, for each $\l >0$, where $A_{\lambda }$ is a non-zero
complex number satisfying $A_{\lambda }^2\in \R $ (one of these
surfaces is represented in Figure~\ref{cathelRiem} right). Together with the
plane, catenoid and helicoid, these examples were characterized by
Riemann~\cite{ri1} as the unique minimal surfaces which are foliated
by circles and lines in parallel planes. Each Riemann minimal
example $M_{\l }$ is topologically a cylinder minus an infinite set
of points which accumulates at infinity to the top or bottom ends of
the cylinder; hence, $M_{\l}$ is topologically the unique planar
domain with two limit ends. Furthermore, $M_{\l }$ is invariant
under reflection in the $(x_1,x_3)$-plane and by a translation
$T_{\l }$; the quotient surface $M_{\l }/T_{\l } \subset \R^3/T_{\l
}$ has genus one and two planar ends, provided that $T_{\l }$ is the
generator of the orientation preserving translations of $M_{\lambda }$.
The conjugate
minimal surface of $M_{\l }$ is $M_{1/\l }$ (the case $\l =1$ gives
the only self-conjugate surface in the family). See~\cite{mpr6} for
a more precise description of these surfaces.

\section{The case with $r$ ends, $2\leq r<\infty $.}
\label{ftc}
Complete minimal surfaces with {\it finite total curvature}
can be naturally thought of as compact algebraic objects,
which explains why these surfaces form
the most extensively studied family among complete minimal surfaces.

\begin{theorem}[Huber~\cite{hu1}, Osserman~\cite{os1}]
\label{thmftcOsserman} Let $M\subset \R^3$ be a complete (oriented),
immersed minimal surface with finite total curvature. Then, $M$ is
conformally a compact Riemann surface $\overline{M}$ minus a finite
number of points, and the Weierstrass representation $(g,dh)$ of $M$
extends meromorphically to $\overline{M}$.  In particular, the total
curvature of $M$ is a multiple of $-4 \pi$.
\end{theorem}
Under the hypotheses of the last theorem, the Gauss map $g$ has a
well-defined finite degree on $\overline{M}$, and equation (\ref{eq:curvtot})
implies that the total curvature of $M$ is $-4\pi $ times the degree
of $g$. The Gauss-Bonnet formula relates the degree of $g$ with the genus of
$\overline{M}$ and the number of ends (Jorge and Meeks~\cite{jm1}); although this
formula can be stated in the more general immersed case, we will only consider it
when all the ends of $M$ are embedded:
\begin{equation}
\mbox{deg}(g)=\mbox{genus}(\overline{M})+\# (\mbox{ends})-1.
\label{eq:JorgeMeeks}
\end{equation}
The asymptotics of a complete, embedded minimal surface in $\R^3$
with finite total curvature are also well-understood: after a
rotation, each embedded end of such a surface is a graph over the
exterior of a disk in the $(x_1,x_2)$-plane with height function
\begin{equation}
\label{eq:grafofinalctfemb}
x_3(x_1,x_2)=a \log
r+b+\frac{c_1x_1+c_2x_2}{r^2}+{\mathcal O}(r^{-2}),
\end{equation}
where $r=\sqrt{x_1^2+x_2^2}$, $a,b,c_1,c_2\in \R $ and ${\mathcal
O}(r^{-2})$ denotes a function such that $r^2{\mathcal O}(r^{-2})$
is bounded as $r\to \infty $ (Schoen~\cite{sc1}). When the {\it
logarithmic growth} $a$ in (\ref{eq:grafofinalctfemb}) is not zero,
the end is called a {\it catenoidal end} (and the surface is
asymptotic to a half-catenoid); if $a=0$, we have a {\it planar end}
(and the surface is asymptotic to a plane). A consequence of the
asymptotics (\ref{eq:grafofinalctfemb}) is that for minimal surfaces
with finite total curvature, completeness is
equivalent to properness (this is also true for
immersed surfaces).

The classification of the complete embedded minimal surfaces
with genus zero and finite total curvature in $\R^3$ was solved in 1991 by
L\'opez and Ros~\cite{lor1}. Their result is based on the fact
that every surface in this family can be deformed through minimal surfaces of the
same type;
this is a strong property which we will encounter in more general situations,
as in Theorem~\ref{thmintegrShiffman} below.

In the finite total curvature setting, the deformation is explicitly
given in terms of the Weierstrass representation: If $(g,dh)$ is the
Weierstrass pair of a minimal surface $M$, then for each $\l >0$ the
pair $(\l g,dh)$ satisfies condition (\ref{eqWeiers1}) and the
second equation in (\ref{eqWeiers2}). The first equation in
(\ref{eqWeiers2}) holds for $(\l g,dh)$ provided that the flux
vector $F(\g )$ of $M$ along every closed curve $\g \subset M$ is
vertical. If $M$ is assumed to have genus zero and finite total
curvature, then the homology classes of $M$ are generated by loops
around its planar and/or catenoidal ends. It is easy to check that
the flux vector of a catenoidal (resp. planar) end along a
non-trivial loop is $(0,0,\pm 2\pi a)$ after assuming that the
limiting normal vector at the end is $(0,0,\pm 1)$. Since
embeddedness implies that all the ends are parallel, then all the
flux vectors of our complete embedded minimal surface $M$ with genus
zero and finite total curvature are vertical; hence $(\l g,dh)$
defines a minimal immersion $X_{\l }$ by the formula (\ref{eq:repW})
for each $\l >0$; note that for $\l =1$ we obtain the starting
surface $M$.

A direct consequence of the maximum principle is that smooth
deformations of compact minimal surfaces remain embedded away from
their boundaries. The strong control on the asymptotics for complete
embedded minimal surfaces of finite total curvature implies
embeddedness throughout the entire deformation $\{ X_{\l }\} _{\l
>0}$. This last property excludes both points in $M$ with horizontal
tangent plane and planar ends (a local analysis of the deformation
around such points and ends produce self-intersections in $X_{\l }$
for values of the parameter $\l $ close to zero or infinity), which
in turn implies that the height function $x_3$ of $M$ is proper
without critical points. A simple application of Morse theory gives
that $M$ has just two ends, in which case the characterization of
the catenoid was previously solved by Schoen~\cite{sc1}. This is a
sketch of the proof of the following result.
\begin{theorem} [L\' opez, Ros~\cite{lor1}]
\label{thmLR}
The plane and the catenoid are the only complete,
embedded minimal surfaces in $\R^3$ with genus zero and finite
 total curvature.
\end{theorem}

The next step in our classification of all properly embedded minimal
planar domains is to understand the case when the number of ends is
finite and at least two. In 1993, Meeks and Rosenberg~\cite{mr6}
showed that if a properly embedded minimal surface $M\subset \R^3$
has at least two ends, then every annular end $E\subset M$ either
has finite total curvature or it satisfies the hypotheses of the
following conjecture, which was solved by Collin in 1997.
\begin{conjecture}[Generalized Nitsche Conjecture, Collin's Theorem~\cite{col1}]
\label{conjNitsche} \mbox{}\newline Let $E\subset \{ x_3\geq 0\} $
be a properly embedded minimal annulus with $\partial E\subset \{
x_3=0\} $, such that $E$ intersects each plane $\{ x_3=t\} $, $t>0$,
in a simple closed curve. Then, $E$ has finite total curvature.
\end{conjecture}
As a direct consequence of the last result and Theorem~\ref{thmLR},
we have that the plane and the catenoid are the unique properly
embedded minimal surfaces in $\R^3$ with genus zero and $r$ ends,
$2\leq r<\infty $.

Collin's original proof of Conjecture~\ref{conjNitsche} is a
beautiful and long argument based on the construction of auxiliary
minimal graphs which serve as guide posts to elucidate the shape of
$E$ in space. For later purposes, it will be more useful for us to
briefly explain a later proof due to Colding and Minicozzi, which is
based on the following scale invariant bound for the Gaussian
curvature of any embedded minimal disk in a half-space.
\begin{theorem}[One-sided curvature estimates, Colding, Minicozzi~\cite{cm23}]
\label{thmcurvestimCM}\mbox{}\newline There exists $\ve >0$ such
that the following holds. Given $r>0$ and an embedded minimal disk
$M\subset \B(2r)\cap \{ x_3>0\} $ with $\partial M\subset \partial
\B(2r)$, then for any component $M'$ of $M\cap \B(r)$
 which intersects $\B(\ve r)$,
\[
\sup _{M'}|K_{M}|\leq r^{-2}.
\]
\end{theorem}
Before sketching the proof of the Nitsche Conjecture, we will make a
few comments about the one-side curvature estimates. The catenoid
shows that the hypothesis in Theorem \ref{thmcurvestimCM} on $M$ to
be simply-connected is necessary. Theorem~\ref{thmcurvestimCM}
implies that if an embedded minimal disk is close enough to (and
lies at one side of) a plane, then reasonably large components of it
are graphs over this plane. The proof of
Theorem~\ref{thmcurvestimCM} is long and delicate,
see~\cite{cm21,cm22,cm23}.

Returning to the Nitsche Conjecture, we see that  it follows
directly from the next result. Given $\ve \in \R $, we denote by
${\mathcal C}_{\ve }$ the conical region $\{ x_3>\ve
\sqrt{x_1^2+x_2^2}\} $.
\begin{theorem}[Colding, Minicozzi~\cite{cm6}]
\label{thmCMNitsche}
There exists $\de >0$ such that any properly embedded
minimal annular end $E\subset {\mathcal C}_{-\de }$ has finite total curvature.
\end{theorem}
{\it Sketch of proof of Theorem~\ref{thmCMNitsche}.}
The argument starts by showing, for each $\de
>0$, the existence of a sequence $\{ y_j\} _j\subset E-{\mathcal C}_{\de }$ with $\| y_j\| \to
\infty $ (this is done by contradiction: if for a given $\de >0$
this property fails, then one use $E$ together with the boundary of
${\mathcal C}_{\de }$ as barriers to construct an end of finite
total curvature contained in ${\mathcal C}_{\de }$, which is clearly
impossible by the controlled asymptotics of catenoidal and planar
ends). The next step consists of choosing suitable radii $r_j>0$
such that the connected component $M_j$ of $E\cap \B(y_j,2r_j)$
which contains $y_j$ is a disk. Now if $\de
>0$ is sufficiently small in terms of the $\ve $ appearing in the one-sided
curvature estimates, we can apply Theorem~\ref{thmcurvestimCM} and
conclude a bound for the supremum of the absolute Gaussian curvature
of the component $M_j^1$ of $M_j\cap \B(y_j,r_j)$ which contains
$y_j$. A Harnack type inequality together with this curvature bound
gives a bound for the length of the intrinsic gradient of $x_3$ in
the intrinsic ball ${\mathcal B}_j$ in $M_j^1$ centered at $y_j$
with radius $5r_j/8$, which in turn implies (by choosing $\ve $
sufficiently small) that ${\mathcal B}_j$ is a graph with small
gradient over $x_3=0$, and one can control a bound by below of the
diameter of this graph. This allows to repeat the above argument
exchanging $y_j$ by a point $y_j^1$ in ${\mathcal B}_j^1$ at certain
distance from $y_j$, and the estimates are carefully done so that
the procedure can be iterated to go entirely around a curve $\g
_j\subset E$ whose projection to the $(x_1,x_2)$-plane links once
around the $x_3$-axis. The graphical property of $\g _j$ implies
that either $\g _j$ can be continued inside $E$ to spiral
indefinitely or it closes up with linking number one with the
$x_3$-axis. The first possibility contradicts that $E$ is properly
embedded, and in the second case the topology of $E$ implies that
$\partial E\cup \g _j$ bounds an annulus $E_j$. The above gradient
estimate gives a linear growth estimate for the length of $\g _j$ in
terms of $\| y_j\| $, from where the isoperimetric inequality for
doubly connected minimal surfaces by Osserman and
Schiffer~\cite{osSchi1} gives a quadratic growth estimate for the
area of $E_j$. Finally, this quadratic area growth property together
with the finiteness of the topology of $E$ imply that $E$ has finite
total curvature by the Gauss-Bonnet formula, finishing the outline
of proof.

\section{The one-ended case.}
\label{sec1conn} In our goal of classifying the properly embedded
minimal surfaces with genus zero, the simplest topology occurs when
the number of ends is one, and the surface is simply-connected. In
spite of this apparent simplicity, this problem remained open until
2005, when Meeks and Rosenberg gave a complete solution by using the
one-side curvature estimates (Theorem~\ref{thmcurvestimCM}) and
other aspects of Colding-Minicozzi theory that we comment on in this
section.

Classical minimal surface theory allows us to understand the
structure of limits of sequences of embedded minimal surfaces with
fixed genus, when the sequence has uniform local area and curvature
bounds, see for instance the survey by P\'erez and Ros~\cite{pro2}.
Colding and Minicozzi faced the same problem in the absence of such
uniform local bounds in a series of papers starting in
2004~\cite{cm25,cm21,cm22,cm24,cm23}. Their most important structure
theorem deals with the case in which all minimal surfaces in the
sequence are disks whose Gaussian curvature blows up near the
origin. To understand this phenomenon, one should think of a
sequence of rescaled helicoids $M_n=\l_nH = \{ \l_nx \ | \ x\in H
\}$, where $H$ is a fixed vertical helicoid with axis the $x_3$-axis
and $\l_n\in \R^+$, $\l_n \searrow 0$. The curvature of the sequence
$\{ M_n\} _n$ blows up along the $x_3$-axis and the $M_n$ converge
away from the axis to the foliation ${\mathcal L}$ of $\R^3$ by
horizontal planes. The $x_3$-axis is the singular set of
$C^1$-convergence $S({\mathcal L})$ of $M_n$ to ${\mathcal L}$, and each leaf $L$ of ${\mathcal L}$ extends smoothly across $L\cap
S({\mathcal L})$ (i.e. $S({\mathcal L})$ consists of removable
singularities of ${\mathcal L}$). The same behavior is mimicked by
any sequence of embedded minimal disks in balls centered at the
origin with radii tending to infinity:

\begin{theorem}[Limit Lamination Theorem for Disks, Colding, Minicozzi~\cite{cm23}]
\label{thmlimitlaminCM}\mbox{}\newline
 Let $M_n\subset \B(R_n)$ be a sequence of embedded minimal
disks with $\partial M_n\subset \partial \B(R_n)$ and $R_n\to \infty $.
If $\sup
|K_{M_n\cap \B(1)}|\to \infty $, then there exists a subsequence of
the $M_n$
(denoted in the same way) and a Lipschitz curve $S\colon \R \to \R^3$
 such that up to a rotation of $\R^3$,
\begin{enumerate}
    \item $x_3(S(t))=t$ for all $t\in \R $.
    \item Each $M_n$ consists of exactly two
    multigraphs$^{\ref{Ngraph}}$ away from $S(\R )$ which  spiral together.
    \item For each $\a \in (0,1)$, the surfaces $M_n-S(\R )$ converge
    in the $C^{\a }$-topology to the foliation
${\mathcal L}=\{ x_3=t\} _{t\in \R }$  by horizontal planes.
    \item $\sup |K_{M_n\cap \B(S(t),r)}|\to
    \infty $ as $n\to \infty $, for any $t\in \R $ and $r>0$.
\end{enumerate}
\end{theorem}
{\it Sketch of proof.} Similar as in the proof of the one-sided
curvature estimates (Theorem~\ref{thmcurvestimCM}), the proof of
this theorem is involved and runs through various
papers~\cite{cm21,cm22,cm23}
(references~\cite{cm27,cm28,cm34,cm33,cm37} by Colding and Minicozzi
are reading guides for the complete proofs of these results). We
will content ourselves with a rough idea of the argument. The first
step consists of showing that the embedded minimal disk $M_n$ with
large curvature at some interior point can be divided into
multivalued graphical building blocks $u_n(\rho ,\t )$ defined on
annuli\footnote{\label{Ngraph}In polar coordinates $(\rho, \t )$
with $\rho >0$ and $\t \in \R $, a {\it $k$-valued graph on an
annulus of inner radius $r$ and outer radius $R$,} is a
single-valued graph of a function $u(\rho ,\t )$ defined over $\{
(\rho ,\t )\ | \ r\leq \rho \leq R,\ |\t |\leq k\pi \} $, $k$ being
a positive integer. The {\it separation} between consecutive sheets
is $w(\rho ,\t )=u(\rho ,\t +2\pi )-u(\rho ,\t )\in \R$.}, and that
these basic pieces fit together properly, in the sense that the
number of sheets of $u_n(\rho ,\t )$ rapidly grows as the curvature
blows up and at the same time, the sheets do not accumulate in a
half-space. This is obtained by means of sublinear and logarithmic
bounds for the {\it separation}$^{\ref{Ngraph}}$ $w_n(\rho ,\t )$ as
a function of $\rho \to \infty$. Another consequence of these bounds
is that by allowing the {\it inner radius}$^{\ref{Ngraph}}$ of the
annulus where the multigraph is defined to go to zero, the sheets of
this multigraph collapse (i.e. $|w_n(\rho ,\t )|\to 0$ as $n\to
\infty $ for $\rho ,\t $ fixed); thus a subsequence of the $u_n$
converges to a smooth minimal graph through $\rho =0$. The fact that
the $R_n$ go to $\infty $ then implies this limit graph is entire
and, by the classical Bernstein's Theorem~\cite{bern}, it is a
plane.

The second step in the proof uses the one-sided curvature estimates
in the following manner: once it has been proven that an embedded
minimal disk $M$ contains a highly sheeted double multigraph
$\widetilde{M}$, then $\widetilde{M}$ plays the role of the plane in
the one-sided curvature estimate, which implies that reasonably
large pieces of $M$ consist of multigraphs away from a cone with
axis ``orthogonal'' to the double multigraph. The fact that the
singular set of convergence is a Lipschitz curve follows because the
aperture of this cone is universal (another consequence of
Theorem~\ref{thmcurvestimCM}).
{\hfill\penalty10000\raisebox{-.09em}{$\Box$}\par\medskip}

With the above discussion in mind, we can now state the main result of
this section.
\begin{theorem}[Meeks, Rosenberg~\cite{mr8}]
\label{ttmr} If $M\subset \R^3$ is a properly embedded,
simply-connected minimal surface, then $M$ is a plane or a helicoid.
\end{theorem}
{\it Sketch of Proof.} Take a sequence $\{ \lambda_n \}_{n}\subset
\R^+$ with $\l _n \to 0$ as $n\to \infty $, and consider the
rescaled surface $\lambda_n M$. Since $M$ is simply-connected,
Theorem~\ref{thmlimitlaminCM} gives that a subsequence of $\l _nM$
converges on compact subsets of $\rth$ to a minimal foliation $\lc$
of $\rth$ by parallel planes, with singular set of convergence
$S(\lc)$ being a Lipschitz curve that can be parameterized by the
height over the planes in ${\cal L}$. Furthermore, a consequence of
the proof of Theorem~\ref{thmlimitlaminCM} in our case is that for
$n$ large, the  almost flat multigraph which starts to form on $\l
_nM$ near the origin extends all the way to infinity. From here in
can be shown that the limit foliation $\lc$ is independent of the
sequence $\{\lambda_n\}_n$. After a rotation of $M$ and replacement
of the $\l _nM$ by a subsequence, we can suppose that the $\l _nM$
converge to the foliation $\lc$ of $\rth$ by  horizontal planes, on
compact subsets outside of the singular set of convergence given by
a Lipschitz curve $S(\lc)$ parameterized by its $x_3$-coordinate. In
particular, $S({\cal L})$ intersects each horizontal plane exactly
once.

The next step consists of proving that $M$ intersects transversely
each of the planes in ${\cal L}$. The idea now is to consider the
solid vertical cylinder $E=\{ x_1^2+x_2^2\leq 1, -1\leq x_3\leq 1\}
$. After a homothety and translation, we can assume that $S({\cal
L})\cap \{ x_3=0\} =\{ \vec{0}\} $ and $S({\cal L})\cap
\{ -1\leq x_3\leq 1\} $ is contained in the convex component of the solid
cone whose boundary has the origin as vertex and that passes through
the circles in $\partial E$. The Colding-Minicozzi picture of
Theorem~\ref{thmlimitlaminCM} implies that for $n$ large, $\l _nM$
intersects $\partial E\cap \{ -1<x_3< 1\} $ in a finite number of
spiraling curves. For simplicity, we will suppose additionally that
the foliation of $\partial E\cap \{ -1<x_3< 1\} $ by horizontal
circles is transversal to $(\l _nM)\cap [\partial E\cap
\{ -1<x_3< 1\} ]$ (in general, one needs to deform slightly these
horizontal circles to almost horizontal Jordan curves to have this
transversality property), and consider the foliation of $E$ given by
the flat disks $D(t)$ bounded by these circles (here $t\in [-1,1]$
denotes height; in general, $D(t)$ is a minimal almost flat disk
constructed by Rado's theorem). Since at a point of tangency, the
minimal surfaces $\l _nM$ and $D(t)$ intersect hyperbolically
(negative index), Morse theory implies that each minimal disk $D(t)$
intersects $\l _nM$ transversely in a simple arc for all $n$ large.
This property together with the openness of the Gauss map of the
original surface $M$, implies that $M$ is transverse to $\lc$, as
desired. In terms of the Weierstrass representation, we now know
that the stereographical projection of the Gauss map $g \colon M \to
\C \cup \{\infty \}$ can be expressed as $g(z)=e^{H(z)}$ for some
holomorphic function $H\colon M \to \C$.

The next goal is to demonstrate that $M$ is conformally $\C $, $M$
intersects every horizontal plane in just one arc and its height
function can be written $x_3=\Re (z)$, $z\in \C $.
In the original proof by Meeks and Rosenberg,
all of these
properties can be deduced from the non-existence of {\it asymptotic
curves} in $M$, a concept that we now explain\footnote{\label{foootn}
For an
alternative short argument, see Remark~\ref{remnew} below.}. Note that the
non-existence of points in $M$ with vertical normal vector implies
that the intrinsic gradient $\nabla x_3$ of the third coordinate
function does not vanish on $M$. An integral curve $\g :[0,\infty
)\to M$ of $\nabla x_3$ is called an {\it asymptotic curve} if $\g $
limits to a finite height as its parameter goes to $\infty $.
Suppose for the moment that $M$ does not admit asymptotic curves.
Consider a component $\G $ of $M\cap \{ x_3=0\} $, which we know it
is smooth. The mapping $F:\G \times \R \to M$ given by $F(p,t)=\g
_p(t)$ where $\g _p$ is the unique integral curve of $\nabla x_3$
with $\g _p(0)=p$, is a local diffeomorphism. Using that $M$ does
not have asymptotic curves, it can be shown that $\partial F(\G
\times \R )$ is empty, hence $F(\G \times \R )=M$ since $M$ is
connected. Now consider the holomorphic function $h=x_3+ix_3^*:M\to
\C $, where $x_3^*$ is the (globally well-defined) harmonic
conjugate function of $x_3$. Again the transversality of $M$ to
every horizontal plane implies that $h$ is a local biholomorphism.
Since $\G \subset \{ x_3=0\} $, $h$ maps $\G $ diffeomorphically
onto an interval $I\subset i\R \subset \C $. As $M$ has no
asymptotic curves, $h$ maps any integral curve $\g _p$ of $\nabla
x_3$ onto a complete horizontal line in $\C $. Thus $h(M)=\R \times
I$ and $h$ is a biholomorphism between these two surfaces. The first
sentence of this paragraph will be proved provided that $I=i\R $.
Otherwise, $[0,\infty )\times I$ is conformally the closed unit disk
minus a closed interval in its boundary, which is not parabolic as a
Riemann surface (i.e. bounded harmonic functions on it are not
determined by their boundary values). Therefore $M\cap \{ x_3\geq
0\}$ is not parabolic, which contradicts Theorem~\ref{thmCKMR}
below. It remains to prove that $M$ does not admit asymptotic
curves. The argument is by contradiction: if $\g \subset M $ is an
asymptotic curve, then one can find a piece of $M$ which is a graph
with infinitely many connected components above a certain horizontal
plane, with zero boundary values. This contradicts the  existence of
an upper bound for the number of components of a minimal graph over
a possibly disconnected, proper domain in $\R^2$ with zero boundary
values (Meeks and Rosenberg proved their own version of this bound
in~\cite{mr8} following previous arguments of Colding and Minicozzi
for harmonic functions; later and sharper versions of this bound can
be found in the papers by Li and Wang~\cite{liw1} and Tkachev~\cite{tk1}).

At this point, we know that the Weierstrass pair of $M$ is
$(g,dh)=(e^{H(z)},dz)$, $z\in \C $, where $H$ is an entire function.
The last step in the proof is to show that $H(z)$ is a linear
function of the form $az+b$ (because in that case $M$ is an
associate surface to a vertical helicoid; but such a surface is
embedded only if it is actually a helicoid). Assuming that
$H$ is a polynomial, the explicit expression
of the Gaussian curvature  $K$ in terms of the Weierstrass  data
implies that $H(z)$ is linear if and only if  $M$ has bounded curvature.
This fact completes the proof of Theorem~\ref{ttmr} provided that
$H$ is a polynomial and $K$ is bounded.
On the other hand, Theorem~\ref{thmlimitlaminCM} and a clever blow-up argument
on the scale of curvature allows us to argue in the bounded curvature
setting, and we then are left with ruling out the case that $H$ has an
essential singularity at $\infty $. This is done by analyzing the
inverse image of a latitude by the Gauss map of the surface. This concludes
our sketch of proof. {\hfill\penalty10000\raisebox{-.09em}{$\Box$}\par\medskip}

In the last proof we mentioned a result on parabolicity for minimal
surfaces with boundary, which we next state.
The proof of this auxiliary result uses the harmonic measure
and universal harmonic functions (see for instance~\cite{mpe1} for these
concepts), and we will skip its proof here.

\begin{theorem}[Collin, Kusner, Meeks, Rosenberg~\cite{ckmr1}]
\label{thmCKMR} Let $M$ be a connected, properly immersed minimal
surface in $\R^3$, possibly with boundary. Then, every component of
the intersection
 of $M$ with a closed half-space is a parabolic surface with boundary.
\end{theorem}

\section{Infinitely many ends I: one limit end is not possible.}
\label{reduction2limitends}

In the sequel, we will consider the case of $M$ being a properly
embedded minimal surface, whose topology is that of a sphere
$\esf^2$ minus an infinite, compact, totally disconnected subset
${\cal E}(M)$. Viewed as a subset ${\cal E}(M)$ of $\esf^2$, the set
${\cal E}(M)$ of ends of $M$ must have accumulation points, which
are called {\it limit ends}\footnote{See Section~2.7 of~\cite{mpe1}
for a generalization of the notion of limit end to a non-compact
connected $n$-manifold.} of $M$. The isolated points in ${\cal
E}(M)$ are called {\it simple ends.}

Next we explain the first two ingredients needed to understand the
geometry of properly embedded minimal surfaces with more than one
end: the notion of limit tangent plane at infinity and the Ordering
Theorem. Every properly embedded minimal surface $M\subset \R^3$
with more than one end admits in one of its two closed complements a
properly embedded minimal surface $\Sigma $ with finite total curvature and
compact boundary (produced via the barrier construction method, see
Section~\ref{subsecbarrier}). By the discussion in
Section~\ref{ftc}, the ends of $\Sigma $ are of catenoidal or planar
type with parallel normal vectors at infinity since $\Sigma $ is
embedded. The plane passing through the origin which is orthogonal
to the limiting normal vectors at the ends of $\Sigma $ does not
depend on $\Sigma$, and it is called the {\it limit tangent plane at
infinity} of $M$ (for details, see Callahan, Hoffman and
Meeks~\cite{chm3}).

\begin{theorem}[Ordering Theorem, Frohman, Meeks~\cite{fme2}]
\label{thordering} Let $M\subset \R^3$ be a properly embedded
minimal surface with more than one end and horizontal limit tangent
plane at infinity. Then, the space ${\cal E}(M)$ of ends of $M$ is
linearly ordered geometrically by the relative heights of the ends
over the $(x_1,x_2)$-plane, and ${\cal E}(M)$ embeds topologically
as a compact totally disconnected subspace of $[0,1]$ in an ordering preserving way.
\end{theorem}
{\it Proof that the subset $X\subset {\cal E}(M)$ of
 all the ends of $M$ with proper
annular representatives has a natural linear ordering.} Suppose
$M\subset \R^3$ is a properly embedded minimal surface and $X$ is
the set of ends which have proper annular representatives. By
Collin's Theorem (Conjecture~\ref{conjNitsche}), every proper
annular representative of an end in $X$ has finite total curvature
and thus, it is asymptotic to a horizontal plane or to a
half-catenoid (recall that $M$ has horizontal limit tangent plane at
infinity). Since these ends of $M$ are all graphs over complements
of compact subdomains in the $(x_1,x_2)$-plane as in
equation~(\ref{eq:grafofinalctfemb}), we have that the set of ends
in $X$ has a natural linear ordering by relative heights over the
$(x_1, x_2)$-plane, and the Ordering Theorem is proved for this
restricted collection of ends.
{\hfill\penalty10000\raisebox{-.09em}{$\Box$}\par\medskip}
\begin{remark}
\label{remark5.2}
%{\rm
The proof of Theorem~\ref{thordering} in the general case is more
involved. For later purposes, we will only indicate what is done.
One starts the proof by using the barrier construction to find ends
of finite total curvature in one of the closed complements of $M$ in
space, adapted to each of the non-annular ends of $M$; more
precisely, one separates each non-annular end representative $E$ of
$M$ with compact boundary from the non-compact domain $M-E$ by a
properly embedded, orientable least area surface $\Sigma _1\subset
\R^3-M$ with $\partial \Sigma _1=\partial E=\partial (M-E)$
constructed using $M$ as a barrier against itself. The asymptotics
of such a  $\Sigma _1$ consists of a positive number of graphical
ends of planar or catenoidal type, with vertical limiting normal
vector. Then one uses these naturally ordered surfaces of the type
$\Sigma _1$ to extend the linear ordering to the entire set of ends
of $M$, see~\cite{fme2} for further details.
%}
\end{remark}
Since ${\mathcal E}(M)\subset [0,1]$ is a compact subspace, the
above linear ordering on ${\mathcal E}(M)$ lets us define the {\it
top} (resp. {\it bottom}) end $e_T$ (resp. $e_B$) of $M$ as the
unique maximal (resp. minimal) element in ${\mathcal E}(M)$. If
$e\in {\mathcal E}(M)$ is neither the top nor the bottom end of $M$,
then it is called a {\it middle} end of $M$. Another key result,
related to conformal properties and area growth, is the following
non-existence result for middle limit ends for a properly embedded
minimal surface.

\begin{theorem}[Collin, Kusner, Meeks, Rosenberg~\cite{ckmr1}]
\label{thmckmr} Let $M\subset \R^3$ be a properly embedded minimal
surface with more than one end and horizontal limit tangent plane at
infinity. Then, any limit end of $M$ must be a top or bottom end. In
particular, $M$ can have at most two limit ends, each middle end is
simple and the number of ends of $M$ is countable.
\end{theorem}
{\it Sketch of proof.} The arguments in Remark~\ref{remark5.2}
insure that every middle end of a surface $M$ as described in
Theorem~\ref{thmckmr} can be represented by a proper subdomain $E
\subset M$ with compact boundary such that $E$ ``lies between two
half-catenoids''. This means that $E$ is contained in a neighborhood
$W$ of the $(x_1,x_2)$-plane, $W$ being topologically a slab, whose
width grows at most logarithmically with the distance from the
origin. This constraint on a middle end representative can be used
in the following way to deduce that the area of this end grows at
most quadratically in terms of the distance to the origin.

For simplicity, we will assume that $E$ is trapped between two
horizontal planes, rather than between two
half-catenoids\footnote{The general case can be treated in a similar
way, although the auxiliary function $f$ is more complicated.}, i.e.
$E\subset W: = \{(x_1, x_2, x_3)$ $| \ r\geq 1, \ 0\leq x_3\leq
1\}$, where $r=\sqrt{x_1^2+x_2^2}$. We claim that both $|\nabla
x_3|^2$ and $\Delta \ln r$ are in $L^1(E)$, where $\nabla ,\Delta $
are the intrinsic gradient and laplacian on $M$. It is not hard to
check using the inequality in equation (\ref{eq:estimlaplaclogr})
below that the function $\ln r- x_3^2$ restricts to every minimal
surface lying in $W$ as a superharmonic function\footnote{This is
called a {\it universal superharmonic function} for the region $W$;
for instance, $x_1$ or $-x_1^2$ are universal superharmonic
functions for all of $\R^3$.}. In particular, the restriction $f
\colon E \rightarrow \re$ is superharmonic and proper.  Suppose
$f(\partial E) \subset [-1, c]$ for some $c>0$.  Replacing $E$ by
$f^{-1}[c, \infty)$ and taking $t>c$, the Divergence Theorem gives
(we can also assume that both $c,t$ are regular values of $f$):
\[
\int_{f^{-1}[c, t]} \Delta f \, dA=-\int_{f^{-1}(c)}|\nabla f|\, ds
+ \int_{f^{-1}(t)}|\nabla f|\, ds,
\]
where $dA,ds$ are the corresponding area and length elements. As $f$
is superharmonic, the function $t\mapsto \int _{f^{-1}[c, t]} \Delta
f\, dA$ is monotonically decreasing and bounded from below by
$-\int_{f^{-1}(c)}|\nabla f|\, ds$.  In particular, $\Delta f$ lies
in $L^1(E)$. Furthermore, $|\Delta f|=|\Delta \ln r - 2 |\nabla
x_3|^2| \geq -|\Delta \ln r|+2|\nabla x_3|^2$.

At this point we need a useful inequality also due to Collin,
Kusner, Meeks and Rosenberg~\cite{ckmr1}, valid for every immersed
minimal surface in $\R^3$:
\begin{equation}
\label{eq:estimlaplaclogr} |\Delta \ln r|\leq
 \frac{| \nabla x_3|^2}{r^2} \qquad \mbox{ in $M-(x_3$-axis)}.
\end{equation}
By the estimate (\ref{eq:estimlaplaclogr}), we have
$|\Delta f|\geq
 \left(2-\frac{1}{r^2}\right) |\nabla x_3|^2$.
Since $r^2\geq 1$ in $W$, it follows $|\Delta f|\geq |\nabla x_3|^2$
and thus, both $|\nabla x_3|^2$ and $\Delta \ln r$ are in $L^1(E)$,
as desired. Geometrically, this means that outside of a subdomain of
$E$ of finite area, $E$ can be assumed to be as close to being
horizontal as one desires, and in particular, for the radial
function $r$ on this horizontal part of $E$, $|\nabla r|$ is almost
equal to 1.

Let $r_0=\max r|_{\partial E}$ and $E(t)$ be the subdomain
of $E$ that lies inside the region $\{ r_0^2\leq x_1^2 +x_2^2
 \leq t^2\}$.
Since
\[
\int_{E(t)} \Delta \ln r \, dA= -\int_{r=r_0}
\frac{|\nabla r|}{r}ds + \int_{r=t} \frac{|\nabla r|}{r}ds
 = \mbox{ const. } +\frac{1}{t}\int_{r=t}|\nabla r|\, ds
\]
and
 $\Delta \ln r\in L^1(E)$, then the following limit exists:
\begin{equation}
\label{eq:ckmr1}
\lim_{t\to\infty} \frac{1}{t}\int_{r=t}|\nabla r|\, ds=C
\end{equation}
for some positive constant $C$. Thus, $t\mapsto \int _{r=t}| \nabla
r|\, ds$ grows at most linearly as $t\to \infty $. By the coarea
formula, for $t_1$ fixed and large,
\begin{equation}
\label{eq:ckmr2}
\int_{E\cap \{ t_1\leq r\leq t\} }|\nabla r|^2\, dA=
\int _{t_1}^t\left( \int _{r=\tau }|\nabla r|\, ds\right) d\tau ;
\end{equation}
hence, $t\mapsto \int _{E\cap \{ t_1\leq r\leq t\} }|\nabla r|^2\, dA$
grows at most quadratically as $t\to \infty $.
Finally, since outside of a domain of finite area $E$ is arbitrarily
close to horizontal and $|\nabla r|$ is
almost equal to one, we conclude that the area of $E\cap \{ r\leq t\} $
grows at most quadratically as $t\to \infty $.
In fact, from (\ref{eq:ckmr1}) and (\ref{eq:ckmr2}) it follows that
\[
\int _{E\cap \{ r\leq t\} }dA=\frac{C}{2}t^2+o(t^2),
\]
where $t^{-2}o(t^2)\to 0$ as $t\to \infty $. Furthermore, it can be
proved that the constant $C$ must be an integer multiple of $2\pi$
(using the quadratic area growth property, the homothetically shrunk
surfaces $\frac{1}{n}E$ converge as $n \rightarrow \infty$ in the
sense of geometric measure theory to a locally finite, minimal
integral varifold $V$ with empty boundary, and $V$ is supported on
the limit of $\frac{1}{n}W$, which is the $(x_1,x_2)$-plane; thus
$V$ is an integer multiple of the $(x_1,x_2)$-plane, which implies
that $C$ must be an integer multiple of $2\pi$).

Finally, every end representative of a minimal surface must have
asymptotic area growth at least equal to half of the area growth of
a plane (as follows from the monotonicity formula,
Theorem~\ref{monotoform}). Since we have checked that each middle
end $e$ of a properly embedded minimal surface has a representative
with at most quadratic area growth, then $e$ admits a representative
which have exactly one end, and this means that $e$ is never a limit
end. {\hfill\penalty10000\raisebox{-.09em}{$\Box$}\par\medskip}

We remark that Theorem~\ref{thmckmr} does not make any assumption on
the genus of the minimal surface. Our next goal is to discard the
possibility of just one limit end for properly embedded minimal
surfaces with {\it finite} genus (in particular, this result holds
in our search of the examples with genus zero), which is the content
of Theorem~\ref{thmno1limitend} below. In order to understand the
proof of this theorem, we will need the following notions and
results.

\begin{definition}
{\rm Suppose that $\{ M_n\} _n$ is a sequence of
connected, properly embedded minimal surfaces in an open set
$A\subset \R^3$. Given $p\in A$ and $n\in \N $, let $r_n(p)>0$ be
the largest radius of an extrinsic open ball $\B \subset A$ centered
at $p$ such that $\B $ intersects $M_n$ in simply-connected
components. If for every $p\in A$ the sequence $\{ r_n(p)\} _n$ is
bounded away from zero, we say that $\{ M_n\} _n$ is  {\it locally
simply-connected in $A$}. If $A=\R^3$ and for all $p\in \R^3$, the
radius $r_n(p)$ is bounded from below by a fixed positive constant
for all $n$ large, we say that $\{ M_n\} _n$ is {\it uniformly
locally simply-connected}. }
\end{definition}

\begin{definition}
{\rm A {\it lamination} of an open subset $U\subset \R^3$ is the
union of a collection of pairwise disjoint, connected, injectively
immersed surfaces with a certain local product structure. More
precisely, it is a pair $({\mathcal L},{\mathcal A})$ where
${\mathcal L}$ is a closed subset of $U$ and ${\mathcal A}=\{
\varphi _{\be }\colon \D \times (0,1)\to U_{\be }\} _{\be }$ is a
collection of coordinate charts of $\R^3$ (here $\D $ is the open
unit disk, $(0,1)$ the open unit interval and $U_{\be }$ an open
subset of $U$); the local product structure is described by the
property that for each $\be $, there exists a closed subset $C_{\be
}$ of $(0,1)$ such that $\varphi _{\be }^{-1}(U_{\be }\cap {\mathcal
L})=\D \times C_{\be }$. It is customary to denote a lamination only
by ${\mathcal L}$, omitting the charts in ${\mathcal A}$. A
lamination ${\mathcal L}$ is a  {\it foliation of $U$} if ${\mathcal
L}=U$. Every lamination ${\mathcal L}$ naturally decomposes into a
union of disjoint connected surfaces, called the {\it leaves} of
${\mathcal L}$. A lamination is {\it minimal} if all its leaves are
minimal surfaces.
 }
\end{definition}
The simplest examples of minimal laminations of $\R^3$ are a closed
family of parallel planes, and ${\cal L}=\{ M\} $, where $M$ is a
properly embedded minimal surface. A crucial ingredient of the proof
of the key topological result described in
Theorem~\ref{thmno1limitend} will be the following structure theorem
of minimal laminations in $\R^3$:

\begin{thm}[Structure of minimal laminations in $\R^3$]
\label{thm2.12}
Let $\lc\subset \R^3$ be a minimal lamination. Then, one of the following possibilities hold.
\begin{enumerate}
\item $\lc=\{ L\} $ where $L$ is a properly embedded minimal surface in $\rth$.
\item $\lc$ has more than one leaf. In this case, $\lc={\cal P}\cup {\cal L}_1$ where ${\cal P}$
consists of the disjoint union of a non-empty closed set of parallel
planes, and ${\cal L}_1$ is a (possibly empty) collection of
pairwise disjoint, complete embedded minimal surfaces. Each $L\in
{\cal L}_1$ has infinite genus, unbounded Gaussian curvature, and is
properly embedded in one of the open slabs and half-spaces
components of $\R^3-{\cal P}$, which we will call $C(L)$. If
$L,L'\in {\cal L}_1$ are distinct, then $C(L)\cap C(L')=\mbox{\O }$.
Finally, each plane $\Pi $ contained in $C(L)$ divides $L$ into
exactly two components.
\end{enumerate}
\end{thm}
The proof of Theorem~\ref{thm2.12} can be found in two papers, due
to Meeks and Rosenberg~\cite{mr8} and Meeks, P\'erez and
Ros~\cite{mpr3}. Although the proof of the above theorem is a bit
long, we next devote some paragraphs to comment some aspects of this
proof since many of the arguments that follow will be used somewhere
else in these notes.
\begin{enumerate}
\item First note that the local structure of a lamination
${\cal L}$ of $\R^3$ %and the one-sided curvature estimates in Theorem \ref{thmcurvestimCM}
implies that the Gaussian curvature
function of the leaves of ${\cal L}$ is locally bounded (in bounded
extrinsic balls). Reciprocally, if a complete embedded minimal
surface $M\subset \R^3$ has locally bounded Gaussian curvature
(bounded in extrinsic balls), then its closure $\overline{M}$ has
the structure of a minimal lamination of $\R^3$. This holds because if $\{ p_n\}
_n\subset M$ is a sequence of points that converges to some $p\in
\overline{M}$, then the local boundedness property of $K_M$ implies
that there exists $\de =\delta (p)>0$ such that for $n$ sufficiently
large, $M$ is a graph $G_n$ over the disk $D(p_n,\delta)\subset
T_{p_n}M$ of radius $\delta$  and center $p_n$, and we have
$C^k$-bounds for these graphs for all $k$. Up to a subsequence, the
planes $T_{p_n}M$ converge to some plane $P$ with $p\in P$, and  the
$G_n|_{D(p_n,\de /2)}$ converge to a minimal graph  $G_\infty$ over
$D(p,\frac{\de}{2})\subset P$; the absence of self-intersections in
$M$ and the maximum principle insure that both $P$ and $G_{\infty }$
are independent of the sequence $\{ p_n\} _n$, and that each $G_n$
is disjoint from $G_\infty$; thus we have constructed a local
structure of lamination around $p\in \overline{M}$, which can be
continued analytically along $\overline{M}$.
\item Given a minimal lamination  ${\cal L}$ of an open set $A\subset \R^3$ and
a point $p\in {\cal L}$, we say that $p$ is a {\it limit point} if for all
$\ve >0$ small enough, the ball $\B (p,\ve )$ intersects ${\cal L}$ in an
infinite number of (disk) components. If a leaf $L\in{\cal L}$ contains a
limit point, then $L$ consists entirely of limit points, and $L$ is then
called a {\it limit leaf} of ${\cal L}$.

\item Let ${\cal L}$ be a minimal lamination of an open set $A\subset \R^3$
and $L$ a limit leaf of ${\cal L}$. Then, the universal covering
space $\widetilde{L}$ of $L$ is stable: By lifting arguments, this property
can be reduced to the case $L$ is simply-connected, and this
particular case follows by expressing every compact disk $\Omega
\subset L$ as the uniform limit of disjoint compact disks $\Omega
(n)$ in leaves $L(n)\in {\cal L}$, and then writing the $\Omega (n)$ as
normal graphs over $\Omega $ of functions $u_n\in C^{\infty
}(\Omega )$ such that $u_n\to 0$ as $n\to \infty $
in the $C^k$-topology for each $k$.
The lamination structure of ${\cal L}$ allows us assume that
$0<u_n<u_{n+1}$ in $\Omega $ for all $n$. After normalizing suitably
$u_{n+1}-u_n$, we produce a positive limit $v\in C^{\infty }(\Omega
)$, which satisfies the linearized version of the minimal surface
equation (the so-called {\it Jacobi equation}), $\Delta v-2K_Lv=0$
in $\Omega $. Since $v>0$ in $\Omega $, it is a standard fact that
the first eigenvalue of the Jacobi operator $\Delta -2K_L$ in
$\Omega $ is positive. Since $\Omega $ is an arbitrary compact
disk in $L$, we deduce that $L$ is stable (in fact, a recent
general result of Meeks, P\'erez and Ros implies that
any limit leaf $L$  of ${\cal L}$ is stable without assuming
it is simply-connected, see \cite{mpr18, mpr19}).
\item A direct consequence of items 1 and 3 above together with
Theorem~\ref{thmstablecompleteplane} is that if $M\subset \R^3$ is a connected,
complete embedded minimal surface with locally bounded curvature, then:
\begin{enumerate}
\item $M$ is proper in $\rth$, or
\item $M$ is proper in an open half-space (resp. slab) of
$\rth$ with limit set the boundary plane
(resp. planes) of this half-space (resp. slab).
\end{enumerate}
In particular, in order to understand the structure of minimal
laminations of $\R^3$, we only need to analyze the behavior in a
neighborhood of a limit plane.
\item In both this item and item 6 below, $M\subset \R^3$
will denote a connected, embedded minimal surface
with $K_M$ locally bounded, such that $M$ is not proper in $\R^3$
and $P$ is a limit plane of $M$; after a rotation, we can assume
$P=\{ x_3=0\} $ and $M$ limits to $P$ from above. We claim that {\it
for any $\ve>0$, the surface $M\cap \{ 0<x_3\leq  \ve \}  $ has unbounded
curvature.} This follows because otherwise, we can choose a smaller
$\ve >0$ so that the vertical projection $\pi \colon M\cap
\{ 0\leq x_3\leq \ve \} \to P$ is a submersion. Let $\Delta$  be a
component of $M\cap \{ 0\leq x_3\leq \ve \} $. Since $\Delta
-\partial \Delta$ is properly embedded in the simply-connected slab
$\{ 0<x_3<\ve \} $, it separates this slab. It follows that each
vertical line in $\rth$ intersects $\Delta$ transversally in at most
one point. This means that $\Delta $ is a
graph over its orthogonal projection in $\{ x_3=0\} $. In
particular, $\Delta$ is proper in $\{ 0\leq x_3\leq \ve \}$. A
straightforward application of the proof of the Halfspace
Theorem~\ref{thmhalf} now gives a contradiction.

\item {\it For every $\ve >0$, the surface $M\cap \{ 0<x_3\leq \ve \}$ is connected.}
The failure of this property produces, using two components
$M(1),M(2)$ of $M\cap \{ 0<x_3\leq \ve \} $ as barriers, a properly
embedded minimal surface $\Sigma \subset \{ 0<x_3\leq \ve \} $ which is
stable, orientable and has the same boundary as, say, $M(1)$. Thus
$\Sigma $ satisfies curvature estimates away from its boundary
(Schoen~\cite{sc3}), which allows us to replace $M$ by $\Sigma $ in
the arguments of item 5 to obtain a contradiction.
\end{enumerate}
By items 1--6 above, Theorem~\ref{thm2.12} will be proved
provided that we check that every $L\in
{\cal L}_1$ has infinite genus, which in turn is a consequence of
the next result.
\begin{proposition}[Meeks, P\'erez and Ros~\cite{mpr3}]
\label{propos2.11}
If $M\subset \R^3$ is a complete embedded minimal surface
with finite genus and locally bounded Gaussian curvature, then $M$ is proper.
\end{proposition}
{\it Sketch of Proof.} One argues by contradiction assuming that $M$
is not proper in $\R^3$. By item~4 above, $M$ is proper in an open
region $W\subset \R^3$ which is (up to a rotation and finite
rescaling) the slab $\{ 0<x_3<1\} $ or the halfspace $\{ x_3>0\} $.
Since $M$ is not proper in $\R^3$, it cannot have finite topology
(by Theorem~\ref{thmCMCalabiYau} below). As $M$ has
finite genus, then it has infinitely many ends. Using similar
arguments as in the proof of the Ordering Theorem
(Theorem~\ref{thordering}), every pair of ends of $M$ can be
separated by a stable, properly embedded minimal surface with
compact boundary. This gives a linear ordering of the ends of $M$ by
relative heights over the $(x_1,x_2)$-plane, in spite of $M$ not
being proper in $\R^3$ but only being proper in $W$. In this new
setting, the arguments in the proof of Theorem~\ref{thmckmr} apply
to give that the middle ends of $M$ are simple, hence topologically
annuli and asymptotically planar or catenoidal.

If some middle end $E$ of $M$ is planar, then all of the middle ends
of $M$ below $E$ are clearly planar. Furthermore, $M$ has unbounded
curvature in every slab $\{ 0<x_3\leq \ve \} $, $\ve >0$ small (by
item 5 before this proposition). The contradiction now follows from
the proof of the curvature estimates of the two-limit-ended case
(Theorem~\ref{thm1} below), which can be extended to this
situation\footnote{We remark that Theorem~\ref{thm1} assumes that
$M$ is proper in $\R^3$; nevertheless, here we need to apply the
arguments in its proof of this situation in which $M$ is proper in
$W$.}.

Hence all the annular ends of $M$ are catenoidal (in particular,
$W=\{ x_3>0\} $). By flux arguments one can show that $M$ cannot
have a top limit end, so it has exactly one limit end which is its
bottom end, that limits to $\{ x_3=0\} $, with all the annular
catenoidal ends of positive logarithmic growth. The rest of the
proof only uses the (connected) portion $M(\ve )$ of $M$ in a slab
of the form $\{ 0<x_3\leq \ve \} $ with $\ve >0$ small enough so that
$M(\ve )$ is a planar domain. Since $K_{M(\ve )}$ is not bounded but
is locally bounded, we can find a divergent sequence $p_n\in M(\ve
)$ with $x_3(p_n)\to 0$ and $|K_{M(\ve )}|(p_n)\to \infty $. The
one-sided curvature estimate (Theorem~\ref{thmcurvestimCM}) implies
that there exists a sequence $r_n\searrow 0$ such that $M\cap \B
(p_n,r_n)$ contains some  component which is not a disk. Then the
desired contradiction will follow after analyzing a new sequence of
minimal surfaces $\widehat{M}(n)$, obtained by {\it blowing-up
$M(\ve )$ around $p_n$ on the scale of topology.} This new notion
deserves some brief explanation, since it will be useful in other
applications (for instance, in the proof of
Theorem~\ref{thmno1limitend} below). It consists of homothetically
expanding $M(\ve )$ so that in the new scale, $p_n$ becomes the
origin and the surface $\widehat{M}(n)$ obtained from expansion of
$M(\ve )$ in the $n$-th step, intersects every ball in $\R^3$ of
radius less than $1/2$ in simply-connected components for $n$ large,
but $\widehat{M}(n)\cap \overline{\B }(1)$ contains a component
which is not a disk (roughly speaking, this is achieved by using the
inverse of the injectivity radius function of $M(\ve )$ as the ratio
of expansion). In the new scale, the contradiction will appear after
consideration of two cases, depending on whether the sequence $\{
K_{\widehat{M}(n)}\} _n$ is locally bounded or unbounded in $\R^3$.

If the sequence $\{ K_{\widehat{M}(n)}\} _n$ is not locally bounded,
then one applies to $\{ \widehat{M}(n)\} _n$ the generalization by
Colding-Minicozzi of Theorem~\ref{thmlimitlaminCM} to planar domains
(see Theorem~\ref{t:t5.1CM} below), and concludes that after
extracting a subsequence, the $\widehat{M}(n)$ converge to a minimal
foliation ${\cal F}$ of $\R^3$ by parallel planes, with singular set
of convergence $S({\cal F})$ consisting of two Lipschitz curves that
intersect the planes in ${\cal F}$ exactly once each. By Meeks' regularity
theorem for the singular set (see Theorem~\ref{thmregular} below),
$S({\cal F})$ consists of two straight lines orthogonal to the
planes in ${\cal F}$. Furthermore, the distance between these two
lines is at most 2 (this comes from the non-simply-connected
property of some component of $\widehat{M}(n)\cap \overline{\B
}(1)$), and this limit picture allows to find a closed curve $\de
_n\subset \widehat{M}(n)$ arbitrarily close to a doubly covered
straight line segment $l$ contained in one of the planes of ${\cal
F}$, such that the each of the extrema of $l$ lies in one of the
singular lines in $S({\cal F})$. Since $\widehat{M}(n)$ has genus
zero, $\de _n$ separates $\widehat{M}(n)$. Therefore, $\de _n$ is
the boundary of a non-compact subdomain $\Omega _n\subset
\widehat{M}(n)$ with a finite number of vertical catenoidal ends,
all with positive logarithmic growth. By the Divergence Theorem, the
flux of $\widehat{M}(n)$ along $\de _n$ is vertical. But this flux
converges to a vector parallel to the planes of ${\cal F}$, of
length at most 4. It follows that the planes of ${\cal F}$ are
vertical. This fact leads to a contradiction with the maximum
principle applied to the function $x_3$ on the domain $\Delta
(n)\subset \widehat{M}(n)$ with finite topology and boundary $\de
_n$ (note that $\Delta (n)$ must lie above the height $\min
x_3|_{\de _n}$ because the ends of $\Delta (n)$ are catenoidal with
positive logarithmic growth). Therefore, $\{ K_{\widehat{M}(n)}\}
_n$ must be locally bounded.

Finally if $\{ K_{\widehat{M}(n)}\} _n$ is locally bounded, then the
arguments in item 1 before Proposition~\ref{propos2.11} (extended to
a sequence of embedded surfaces instead of a single surface) imply
that after extracting a subsequence, the $\widehat{M}(n)$ converge
to a minimal lamination $\widehat{\cal L}$ of $\R^3$. The
non-simply-connected property of some component of
$\widehat{M}(n)\cap \overline{\B }(1)$ is then used to prove that
$\widehat{\cal L}$ contains a non-flat leaf $\widehat{L}$. Since
$\widehat{L}$ is not a plane, it is not stable, which in turns
implies that the multiplicity of the convergence of portions of
$\widehat{M}(n)$ to $\widehat{L}$ is 1 (the arguments for this
property are similar to those applied in item 3 before
Proposition~\ref{propos2.11}). This fact and the verticality of the
flux of the $\widehat{M}(n)$ imply that $\widehat{L}$ has vertical
flux as well. Then one finishes this case by discarding all
possibilities for such an $\widehat{L}$ as a limit of the
$\widehat{M}(n)$; the verticality of the flux of $\widehat{L}$
together with the L\'opez-Ros deformation argument explained in
Section~\ref{ftc} are crucial here (for instance, one of the
possibilities for $\widehat{L}$ is being a properly embedded minimal
planar domain with two limit ends; this case is discarded by using
Theorem~\ref{thm8} below).
 See Theorem~7 in~\cite{mpr3} for further details.
{\hfill\penalty10000\raisebox{-.09em}{$\Box$}\par\medskip}
\par
\vspace{.2cm} In the last proof we used three auxiliary results,
which we next state for future reference.
\begin{theorem}[Colding-Minicozzi~\cite{cm35}]
\label{thmCMCalabiYau} If $M\subset \R^3$ is a complete, embedded
minimal surface with finite topology, then $M$ is proper.
\end{theorem}

\begin{theorem}[Limit Lamination Theorem for Planar Domains, Colding, Minicozzi~\cite{cm25}]
\label{t:t5.1CM} Let $M_n\subset \B(R_n)$ be a locally
simply-connected sequence of embedded minimal planar domains with
$\partial M_n\subset \partial \B(R_n)$, $R_n\to \infty $, such that
$M_n\cap \B(2)$ contains a component which is not a disk for any
$n$. If $\sup| K_{M_n\cap \B(1)}|\to \infty \, $, then there exists
a subsequence of the $M_n$ (denoted in the same way) and two
vertical lines $S_1,S_2$, such that
\begin{description}
\item[{\it (a)}] $M_n$ converges away from $S_1\cup S_2$
to the foliation ${\cal F}$ of $\R^3$ by horizontal planes.
\item[{\it (b)}] Away from $S_1\cup S_2$, each $M_n$
consists of exactly two multivalued graphs spiraling together. Near
$S_1$ and $S_2$, the pair of multivalued graphs form double spiral
staircases with opposite handedness at $S_1$ and $S_2$. Thus,
circling only $S_1$ or only $S_2$ results in going either up or
down, while a path circling both $S_1$
and $S_2$ closes up.%, see Figure~\ref{ULSC1}.
%\begin{figure}
%\centerbmp{14.7cm}{5.3cm}{ULSC1.bmp x=14.7cm y=5.3cm}
%\caption{Left: Two oppositely handed double spiral staircases. Right: The limit foliation by parallel
% planes and the singular set of convergence $S_1\cup S_2$.}
%\label{ULSC1}
%\end{figure}
\end{description}
\end{theorem}

\begin{theorem}[Regularity of $S({\cal L})$, Meeks~\cite{me25}]
\label{thmregular} Suppose $\{ M_n\}_n$ is a locally
simply-connected sequence of properly embedded minimal surfaces in a
three-manifold, that converges $C^{\a}$ to a minimal foliation
${\cal L}$  outside a locally finite collection of Lipschitz curves
$S({\cal L})$ transverse to ${\cal L}$. Then, $S(\mathcal{L})$
consists of a locally finite collection of integral curves of the
unit Lipschitz normal vector field to $\mathcal{L}$. In particular,
the curves in $S({\cal L})$ are $C^{1,1}$ and orthogonal to the
leaves of ${\cal L}$.
\end{theorem}

After all these preliminaries, we are ready to discard
the one-limit-ended case in our search of all properly
embedded minimal surfaces with finite genus.
\begin{theorem}[Meeks, P\'erez, Ros \cite{mpr4}]
\label{thmno1limitend}
 If $M\subset \R^3$ is a properly embedded minimal surface
with finite genus, then $M$ cannot have exactly one limit end.
\end{theorem}
{\it Sketch of proof.} After a rotation, we can suppose that $M$ has
a horizontal limit tangent plane at infinity and by
Theorems~\ref{thordering} and~\ref{thmckmr}, its ends are linearly
ordered by relative increasing heights, ${\mathcal E}(M)=\{
e_1,e_2,\ldots ,e_{\infty }\} $ with $e_{\infty }$ being the limit
end of $M$ and its top end. Since for each finite $n$ the end $e_n$
has a annular representative, Collin's Theorem implies that $e_n$
has a representative with finite total curvature, which is therefore
asymptotic to a graphical annular end of a vertical catenoid or
plane. By embeddedness, this catenoidal or planar end has vertical
limit normal vector, and its logarithmic growth $a_n$ satisfies
$a_n\leq a_{n+1}$ for all $n$. Note that $a_1<0$ (otherwise we
contradict Theorem~\ref{thmhalf}). By flux reasons, $a_n<0$ for all
$n$ (one $a_k$ being positive implies that $a_n\geq a_k>0$ for all
$k\geq n$, which cannot be balanced by finitely many negative
logarithmic growths; a similar argument works if one $a_k$ is zero).

The next step consists of analyzing
the limits of $M$ under homothetic shrinkings:
\begin{claim}
\label{claim5.12}
Given a sequence $\l _n\searrow 0$,
the sequence of surfaces $\{ \l _nM\} _n$ is locally
simply-connected in $\R^3-\{ \vec{0}\} $.
\end{claim}
The proof of this property is by contradiction: if it fails around a
point $p\in \R^3-\{ \vec{0}\} $, then one blows-up $\l _nM$ on the
scale of topology around $p$, as we did in the proof of
Proposition~\ref{propos2.11}. Thus, we produce expanded versions
$\widehat{M}(n)$ of $\l _nM$ so that $p$ becomes the origin and
$\widehat{M}(n)$ intersects any ball in $\R^3$  of radius less than
$1/2$ in simply-connected components for $n$ large, but
$\widehat{M}(n)\cap \overline{\B}(1)$ contains a component which is
not simply-connected. We will now try to mimic the arguments in the
last two paragraphs of the proof of Proposition~\ref{propos2.11},
commenting only the differences between the two situations.

If $\{ K_{\widehat{M}(n)}\} _n$ is not locally bounded in $\R^3$,
then Theorems~\ref{t:t5.1CM} and~\ref{thmregular} produce a limit
picture for the $\widehat{M}(n)$ which is a foliation ${\cal F}$ of
$\R^3$ by parallel planes, with two parallel straight lines as
singular set of convergence. Then one finds connection loops $\de
_n\subset \widehat{M}(n)$ as in the proof of
Proposition~\ref{propos2.11}, which together with a flux argument,
imply that the foliation ${\cal F}$ is by vertical planes. Now the
contradiction comes from the maximum principle applied to the
function $x_3$ on the domain $\Delta (n)\subset \widehat{M}(n)$ with
finite topology and boundary $\de _n$ (note that $\Delta (n)$ now
lies below the height $\max x_3|_{\de _n}$ because the ends of
$\Delta (n)$ are catenoidal with negative logarithmic growth,
compare with the situation in the next to the last paragraph of the
proof of Proposition~\ref{propos2.11}).

If $\{ K_{\widehat{M}(n)}\} _n$ is locally bounded in $\R^3$, then
after replacing by a subsequence, the $\widehat{M}(n)$
converge to a minimal lamination $\widehat{\cal L}$ of $\R^3$ with
at least one non-flat leaf $\widehat{L}$. The difference now with the
last paragraph in the proof of Proposition~\ref{propos2.11}) is that
we know that $\widehat{L}$ is proper in $\R^3$ (since it has
genus zero and by Proposition~\ref{propos2.11}). Therefore,
this current situation is even simpler than the one in the
last paragraph of the proof of Proposition~\ref{propos2.11},
where we discarded all the possibilities for $\widehat{L}$
by using that its flux is vertical (for instance, the case of
$\widehat{L}$ being a two-limit ended planar domain is
discarded using Theorem~\ref{thm8} below).
This finishes the sketch of the proof of Claim~\ref{claim5.12}.

\par
\vspace{.2cm} Once we know that $\{ \l _nM\} _n$ is locally
simply-connected in $\R^3-\{ \vec{0}\} $, it can be proved that the
limits of certain subsequences of $\{ \l _nM\} _n$ consist of
(possibly singular) minimal laminations ${\cal L}$ of $H(*)=\{
x_3\geq 0\} -\{ 0\} \subset \R^3$ containing $\partial H(*)$ as a
leaf. Subsequently, one checks that every such a limit lamination
${\cal L}$ has no singular points and that the singular set of
convergence $S({\mathcal L})$ of $\l _nM$ to ${\mathcal L}$ is
empty. In particular, taking $\l _n=\| p_n\| ^{-1}$ where
 $p_n$ is any divergent sequence on
$M$, the fact that $S({\mathcal L})=\mbox{\O }$ for
the corresponding limit minimal
lamination ${\mathcal L}$ insures that the Gaussian
 curvature of $M$ decays at least
quadratically in terms of the distance function to
 the origin. In this situation, the Quadratic Curvature Decay
Theorem (Theorem~\ref{thm1introd} below), insures that $M$ has
finite total curvature, which is impossible since $M$ has an
infinite number of ends.
{\hfill\penalty10000\raisebox{-.09em}{$\Box$}\par\medskip}

We finish this section by stating an auxiliary result that was used in
the last paragraph of the sketch of proof of Theorem~\ref{thmno1limitend}.
\begin{theorem}[Quadratic Curvature Decay Theorem, Meeks, P\'erez,
Ros~\cite{mpr10}] \label{thm1introd}\mbox{}\newline Let $M\subset
\R^3-\{ \vec{0}\} $ be an embedded minimal surface with compact
boundary (possibly empty), which is complete outside the origin
$\vec{0}$; i.e. all divergent paths of finite length on~$M$ limit
to~$\vec{0}$. Then, $M$
%A complete, embedded minimal
%surface in $\R^3$ with compact boundary
%   (possibly empty)
has quadratic decay of curvature if and only if
its closure in $\R^3$ has finite total curvature.
\end{theorem}

\section{Infinitely many ends II: two limit ends.}
\label{sec2limitends}

Consider a properly embedded minimal planar domain $M\subset \R^3$
with two limit ends and horizontal limit tangent plane at infinity.
Our goal in Sections~\ref{sec2limitends} and~\ref{secKdV} is to
prove that $M$ is one of the Riemann minimal examples introduced in
Section~\ref{subsecexamples}.

From the preceding sections we know that the limit ends of $M$ are
its top and bottom ends and each middle end of $M$ has a
representative which is either planar or catenoidal. Indeed, all
middle ends of $M$ are planar: otherwise $M$ contains a catenoidal
end $e$ with, say, logarithmic growth $a>0$. Consider an embedded
closed curve $\g \subset M$ separating the two limit ends of $M$
(for this is useful to think of $M$ topologically as a sphere minus
a closed set ${\cal E}$ which consists of an infinite sequence of
points that accumulates at two different points). Let $M_1$ be the
closure of the component of $M-\g $ which contains $e$. By
embeddedness, all the annular ends of $M_1$ above $e$ must have
logarithmic growth at least $a$, which contradicts that the flux of
$M$ along $\g \subset M$ is finite.

Once we know that the middle ends of $M$ are all planar, we can
separate each pair of consecutive planar ends of $M$ by a horizontal
plane which intersects $M$ in a compact set. An elementary analysis
of the third coordinate function restricted to the subdomain
$M(-)=M\cap \{ -\infty <x_3\leq 0 \} $ (which is a parabolic Riemann
surface with boundary by Theorem~\ref{thmCKMR}) shows that $M(-)$ is
conformally equivalent to the closed unit disk $\overline{\D }$
minus a closed subset $\{ \xi _k\} _k\cup \{ 0\} $ where $\xi _k\in
\overline{\D }-\{ 0\} $ converges to zero as $k\to \infty $, and
that the third coordinate function of $M$ can be written as $x_3(\xi
)=\a \log |\xi |$, $\xi \in \overline{\D }$, where $\a >0$ (in
particular, there are no points in $M$ where the tangent plane is
horizontal, or equivalently $M$ intersects every horizontal plane of
$\R^3$ in a Jordan curve or an open arc). After a suitable homothety
so that $M$ has vertical component of its flux vector along a
compact horizontal section equal to $1$ (i.e. $\a =\frac{1}{2\pi
}$), we deduce that the following properties hold:
\begin{enumerate}
\item $M$ can be conformally parameterized by the cylinder
$\C /\langle i\rangle $ (here $i=\sqrt{-1}$) punctured in an
infinite discrete set of points $\{ p_j,q_j\} _{j\in \Z }$ which
correspond to the planar ends of $M$.
\item The stereographically projected Gauss map of $M$
extends through the planar ends to a meromorphic function $g$ on $\C
/\langle i\rangle $ which has double zeros at the points $p_j$ and
double poles at the $q_j$ (otherwise we contradict that $M$
intersects every horizontal plane asymptotic to a planar end in an
open arc).
\item The height differential of $M$ is $dh=dz$ with $z$ being the
usual conformal coordinate on~$\C $ (equivalently,
the third coordinate function of $M$ is $x_3(z)=\Re (z)$).
\item The planar ends are ordered by their heights so
that $\Re (p_j)<\Re (q_j)<\Re (p_{j+1})$ for all
$j$ with $\Re (p_j)\to \infty $ (resp. $\Re (p_j) \to -\infty $)
when $j\to \infty $ (resp. $j\to -\infty $).
\end{enumerate}

\subsection{Curvature estimates and quasiperiodicity.}
The next step consists of proving that every surface $M$ as before
admits an estimate for its absolute Gaussian curvature that depends
solely on an  upper bound for the horizontal component of the flux
vector $F(M)$ of $M$ along a compact horizontal section. Note that
$F(M)$ does not depend on the height of the horizontal plane which
produces the compact section since the flux of a planar end is zero;
for this reason, we will simply call $F(M)$ the {\it flux vector} of
$M$.

\begin{theorem}[Meeks, P\'erez, Ros~\cite{mpr3}]
\label{teorKestim} Let $\{ M(i)\} _i\subset {\cal M}$ be a sequence
of properly embedded minimal surfaces with genus zero, two limit
ends, horizontal limit tangent plane at infinity and flux vector
$F(M(i))=(H(i),1)\in \R^2\times \R $. If $\{ H(i)\} _i$ is bounded
from above in $\R^2 $, then the Gaussian curvature function
$K_{M(i)}$ of the $M(i)$ is uniformly bounded.
\end{theorem}
{\it Sketch of proof.} Arguing by contradiction, assume $\{
|K_{M(i)}|\} _i$ is not uniformly bounded. Then, we blow up $M(i)$
{\it on the scale of curvature,} which means that we choose suitable
points $p(i)\in M(i)$ such that after translating $M(i)$ by $-p(i)$
and expanding by $\l (i)=\sqrt{|K_{M(i)}(p(i))|}\to +\infty $, we
produce new properly embedded minimal surfaces $M'(i)\subset \R^3$
which, after passing to a subsequence, converge uniformly on compact
subsets of $\R^3$ with multiplicity 1 to a properly embedded minimal
surface $H\subset \R^3$ with $\vec{0}\in H$, $|K_H(\vec{0})|=1$ and
$|K_H|\leq 1$ in $H$. Since $H$ is a limit with multiplicity 1 of
surfaces of genus zero, then $H$ has genus zero.

If $H$ is not simply-connected, then the discussion in previous
sections shows that one can find an embedded closed curve $\g
\subset H$ such that the flux of $H$ along $\g $ is finite and
non-zero. This is a contradiction, since $\g $ produces related
non-trivial loops $\g '(i)\subset M'(i)$ converging to $\g $ as
$i\to \infty $; if we call $\g (i)$ to the loop in $M(i)$ which
corresponds to $\g'(i)$ in the original scale, then the third
component of the flux of $M(i)$ along $\g (i)$ times $\l (i)$
converges as $i\to \infty $ to the third component of the flux of
$H$ along $\g $, which is finite. Since $\l (i)\to \infty $, we
deduce that the third component of the flux of $M(i)$ along $\g (i)$
tends to zero, which is impossible by our normalization on $F(i)$.
Therefore, $H$ is simply-connected.

Since $H$ is a non-flat, properly embedded minimal surface  which is
simply-connected, Theorem~\ref{ttmr} implies that $H$ is a helicoid.
Furthermore, $H$ is a vertical helicoid, since its Gauss map does
not take vertical directions (because the Gauss maps of the surfaces $M'(i)$
share the same property). Once this first helicoidal limit $H$ of
rescalings of the $M(i)$ has been found, one rescales and rotates
again the $M(i)$ in a rather delicate way:
\begin{claim}
\label{claim6.2} There exist a universal $\tau _0>0$ and angles $\t
(i)\in [0,2\pi )$ such that for any $\tau >\tau _0$, one can find
numbers $\mu (\tau ,i)>0$ and embedded closed curves $\de (\tau
,i)\subset \mu (\tau ,i) \mbox{\rm Rot}_{\t (i)}(M(i)-p(i))$ (here
$\mbox{\rm Rot}_{\t}$ denotes the rotation of angle $\t $ around the
$x_3$-axis) so that the flux of the rotated and rescaled surface
$\mu (\tau ,i) \mbox{\rm Rot}_{\t (i)}(M(i)-p(i))$ along $\delta
(\tau ,i)$ decomposes as
\begin{equation}
\label{eq:teorKestim1}
\mbox{\rm Flux}\left( \mu (\tau ,i) \mbox{\rm Rot}_{\t (i)}(M(i)-p(i)),\de (\tau
,i)\right) =V(\tau ,i)+W(\tau ,i)
\end{equation}
where $V(\tau ,i),W(\tau ,i)\in \R^3$ are vectors such that $\lim
_{i\rightarrow \infty }V(\tau, i)=(12\tau ,0,0)$
and $\{ \| W(r,i)\| \} _i$ is bounded by a constant independent of $\tau $.
\end{claim}
Assuming this technical property, the proof of
Theorem~\ref{teorKestim} finishes as follows. First note that for
any properly embedded minimal surface $M\subset \R^3$ with genus
zero, two limit ends and horizontal limit tangent plane at infinity,
the angle between the flux vector $F(M)$ and its horizontal
component $(H(M),0)$ is invariant under translations, homotheties
and rotations around the $x_3$-axis. By (\ref{eq:teorKestim1}), the
corresponding angles for the flux vectors of the surfaces $ \mu
(\tau ,i) \mbox{Rot}_{\t (i)}(M(i)-p(i))$ tend to zero as
$i\rightarrow \infty $ and $\tau \rightarrow \infty $. But those
angles are nothing but the angles for $M(i)$, which are bounded away
from zero because of the hypothesis of $\{ H(M(i))\} _i$ being
bounded above. This contradiction proves the theorem, modulo the
above claim.

The proof of Claim~\ref{claim6.2} is delicate, and we will only
mention that it uses some results of  Colding-Minicozzi
(Theorems~\ref{thmcurvestimCM} and~\ref{t:t5.1CM}) and the rescaling
on the scale of topology that we explained in the proof of
Proposition~\ref{propos2.11}. For details, see~\cite{mpr3}.
{\hfill\penalty10000\raisebox{-.09em}{$\Box$}\par\medskip}

With the curvature bound given in Theorem~\ref{teorKestim} in hand,
one can use standard arguments based on the maximum principle to
find an embedded regular neighborhood of constant positive radius
bounded from below by a constant which only depends on the curvature
bound of the minimal surface (see for instance Lemma 1
in~\cite{mpr1}). The existence of such a regular neighborhood
implies uniform local area bounds for a sequence of surfaces $\{
M(i)\} _i$ under the hypotheses of Theorem~\ref{teorKestim}.
Finally, these curvature and area bounds allow one to apply
classical results for taking limits (after extracting a subsequence)
of suitable translations of the $M(i)$. Summarizing, we have the
next statement.

\begin{theorem} \label{thm1}
Suppose $M$ is a properly embedded minimal surface in $\R^3$ with
genus zero and two limit ends. Assume that $M$ is normalized by a
rotation and homothety so that it has horizontal limit tangent plane
at infinity and the vertical component of its flux equals~1.  Then:
\begin{enumerate}
\item The middle ends $\{e_n \mid n \in \Z \}$ of $M$ are
planar and have heights ${\cal H} = \{x_3(e_n) \mid n \in \mathbb{Z}
\}$ such that $x_3(e_n) < x_3 (e_{n+1})$ for all $n\in \Z $;
\item ${\displaystyle \lim_{n \rightarrow
\infty} x_3(e_n) = \infty}$ and ${\displaystyle \lim_{n \rightarrow
-\infty} x_3(e_n) = - \infty}$;
\item Every horizontal plane intersects $M$
in a simple closed curve when its height is not in ${\cal H}$ and in a
single properly embedded arc when its height is in ${\cal H}$;
\item $M$ has bounded Gaussian curvature, with the
bound of its curvature depending only on an upper bound of the
horizontal component of the flux of $M$.

\item If the Gaussian curvature of $M$ is bounded from below in absolute
value by $\ve ^2$, then $M$ has a regular neighborhood of radius
$1/\ve $ and so, the spacings $S(n) = x_3 (e_{n+1}) - x_3(e_n)$
between consecutive ends are bounded from below by $2/\ve $.
Furthermore, these spacings are also bounded by above.
\item $M$ is quasiperiodic in the following
sense.  There exists a divergent sequence $V(n)\in \R ^3$ such that
the translated surfaces $M + V(n)$ converge to a properly embedded
minimal surface $M(\infty )$ of genus zero, two limit ends, horizontal limit tangent
plane at infinity and with the same flux
as $M$.
\end{enumerate}
\end{theorem}
{\it Sketch of proof.} Everything has been proved
(at least, outlined with some details) except for item~{\it 6}
and the last sentence of item~{\it 5}. To see that item~{\it 6}
holds, one starts by noticing that for any two consecutive planar
ends $e_n, e_{n+1}$ of $M$ there is a point $p(n)$ with
$x_3(e_n)<x_3(p(n))<x_3(e_{n+1})$ where the tangent plane of $M$ is
vertical (by continuity, since the Gauss map alternates from the
north pole to the south pole or vice versa at $e_n,e_{n+1}$).
Extrinsically close to $p(n)$ we can find a point $q(n)\in M$ such
that $|K_M(q(n)|\geq \ve >0$, for $\ve$ fixed and small (otherwise
we produce relatively large flat vertical regions in $M$ around
$p(n)$, which contradicts that the third component of the flux
vector of $M$ is one). Now one considers the translated surfaces
$M-q(n)$, which have local area and curvature estimates. The
arguments in the paragraph just before the statement of
Theorem~\ref{thm1} give that a subsequence of the $M-q(n)$ converges
with multiplicity 1 to a connected, non-flat, properly embedded
minimal surface $M(\infty)$ of genus zero, with length of the
vertical component of its flux vector at most 1. Since $M(\infty )$
is not flat, the open mapping property for the Gauss map of minimal
surfaces implies that the Gauss map of $M(\infty )$ omits the
vertical directions. Since the vertical component of the flux of
$M(\infty )$ is finite, $M(\infty )$ cannot be a helicoid. By
Theorem~\ref{ttmr}, $M(\infty )$ is not simply-connected. Since it
has genus zero, $M(\infty )$ has at least two ends, and thus it must
have a planar or a catenoid type end  (which must be horizontal). In
particular, $M(\infty)$ has a horizontal limit tangent plane at
infinity.

Next we check that $M(\infty )$ cannot have catenoidal ends.
Otherwise, $M(\infty)$ contains a horizontal, strictly convex Jordan
curve $\G (\infty )$, along which this limit surface has vertical
flux. Since $M(\infty )$ is the limit of (a subsequence of)
$M-q(n)$, the curve $\G (\infty )$ produces horizontal, strictly
convex Jordan curves $\G (n)\subset M$, with $\G (n)-q(n)\to \G
(\infty )$ as $n\to \infty $. As the ends of $M$ are planar and the
flux of $M(\infty )$ along $\G (\infty )$ is vertical, we deduce
that the flux of $M$ along $\G(n)$ is also vertical for $n$ large
enough. Now one applies a variant of the L\' opez-Ros deformation
argument to the portion of $M$ bounded by $\G (n)\cup \G (n+k)$,
with $n,k$ large, to find a contradiction. Therefore, $M(\infty )$
cannot have catenoidal ends.

Finally, the non-existence of catenoidal ends on $M(\infty )$
implies that $M(\infty )$ neither has finitely many ends (recall it
is not simply-connected) nor does it have exactly one limit
end\footnote{Although $M(\infty )$ is known to have genus zero, we
cannot apply Theorem~\ref{thmno1limitend} here to avoid the
possibility that $M(\infty )$ has just one limit end, since the end
of the proof of Claim~\ref{claim5.12} uses this Theorem~\ref{thm1}
to discard the existence of a certain minimal surface $\widehat{L}$
with two limit ends.}. Therefore, $M(\infty )$ has two limit ends
and item {\it 6} of Theorem~\ref{thm1} is proved. The last sentence
in item~{\it 5} is also proved, since if the spacing between the
ends of $M$ were unbounded, then a variation of the above arguments
would yield a limit surface of translations of $M$ with a catenoidal
end, which is impossible.
{\hfill\penalty10000\raisebox{-.09em}{$\Box$}\par\medskip}

For the sake of completeness, we next state a technical result which
was used at the end of the proofs of Proposition~\ref{propos2.11}
and Claim~\ref{claim5.12}. Although this proposition and this claim
were stated in earlier sections of this article, chronologically
they were obtained later than the next property.

\begin{theorem}[Meeks, P\'erez, Ros~\cite{mpr3}]
\label{thm8} Let $M\subset \R^3$ be a properly embedded minimal
surface with genus zero, two limit ends and horizontal limit tangent
plane at infinity. Then, $M$ has non-zero horizontal flux along any
compact horizontal section.
\end{theorem}
{\it Sketch of Proof.} Arguing by contradiction, suppose that $M$
has vertical flux. The main tool used here is again a variant of the
L\'opez-Ros deformation (see Section~\ref{ftc}). Since $M$ has
vertical flux, the L\'opez-Ros deformation produces a one-parameter
family $\{ X_{\l }\colon M\to \R^3\} _{\l }$ of complete minimal
immersions with $X_{1}(M)=M$. Since the middle ends of $M$ are
planar, for $\l $ large $X_{\l }$ has self-intersections. As $X_{1}$
is an embedding, it makes sense to consider the largest $\l_0 \in
[1,\infty)$ such that $X_{\l }$ is injective for all $\l \in [1,\l
_0)$. The usual maximum principle for minimal surfaces implies that
limits of embedded minimal surfaces are embedded, from where we
conclude that $X_{\l _0}$ is also an embedding. The contradiction
will follow by proving that, whenever $X_{\l }$ is injective, then
for nearby parameter values $\l '$, $X_{\l '}$ is also injective
(just apply this property to $\l =\l _0$). First note that the
height differential is preserved along the L\'opez-Ros deformation,
and thus the heights of the planar ends of $X_{\l }$ remain constant
in $\l $. From here is not difficult to prove the desired openness
property of the injective parameter values, provided that the
surface $M$ is periodic (invariant by a translation in $\R^3$),
since one can then argue with the quotient surface under the
translation, which has finitely many ends. The key observation here
is that the quasiperiodicity property in item~{\it 6} of
Theorem~\ref{thm1} applied to $X_{\l _0}(M)$ is enough to preserve
the embeddedness of $X_{\l }(M)$ for nearby values of $\l $.
{\hfill\penalty10000\raisebox{-.09em}{$\Box$}\par\medskip}

\subsection{The Shiffman function.}
Next we explain the last tool necessary to finish our classification
of properly embedded minimal planar domains with two limit ends: the
{\it Shiffman function.} This is a particular example of a {\it
Jacobi function,} which can be considered to be an infinitesimal
deformation for a minimal surface $M$ of the kind we are
considering.

In the sequel we will denote by ${\cal M}$ the space of properly
embedded, minimal planar domains ${M\subset\rth}$ with two limit
ends, horizontal limit tangent plane at infinity and flux vector
$F(M)=(h,0,1)$ for some $h=h(M)>0$, identified up to translations.
Recall that we proved that for every $M\in {\cal M}$, properties
1--4 stated before Theorem~\ref{teorKestim} hold. In particular, the
planar sections $M\cap \{ x_3=c\} $, $c\in \R $, are either Jordan
curves or Jordan arcs. In the natural conformal coordinate $z=x+iy$
such that the height differential of $M$ is $dz$, the level curve
$M\cap x_3^{-1}(c)$ corresponds to $y\mapsto z_c(y)=c+iy$ in the
$z$-plane. The planar curvature of this level curve is
\begin{equation}
\label{eq:curvplanarsect} \kappa _c(y)=\left. \left[
\frac{|g|}{1+|g|^2}\Re \left( \frac{g'}{g}\right) \right] \right|
_{z=z_c(y)},
\end{equation}
where $g$ is the stereographically projected Gauss map of $M$ and the
prime stands for derivative with respect to $z$.

\begin{definition}
\label{defShiffman}
{\rm In the above setting, the {\it Shiffman function} of $M$ is
\begin{equation}
\label{eq:uShiffman}
S_M=\Lambda \frac{\partial \kappa
 _c}{\partial y}=\Im \left[ \frac{3}{2}\left( \frac{g'}{g}\right) ^2-\frac{g''}{g}
 -\frac{1}{1+|g|^2}\left( \frac{g'}{g}\right) ^2\right],
 \end{equation}
 where the induced metric $ds^2$ on $M$ by the inner product of
$\R^3$ is $ds^2=\Lambda ^2|dz|^2$ (i.e. $\Lambda =\frac{1}{2}(|g|+|g|^{-1})$) and
$\Im $ stands for imaginary part.
 }
\end{definition}
\begin{remark}
{\it The Shiffman function can also be defined for surfaces not in
${\cal M}$: we only need $(g(z),dh=dz)$ to be the Weierstrass pair
of a minimal surface with respect to a local conformal coordinate.
Every minimal surface admits such a local representation around a
point with non-vertical normal vector.}
\end{remark}
Perhaps the most remarkable property of $S_M$ is that, from a
variational point of view, it is an infinitesimal deformation of $M$
by minimal surfaces, i.e. $S_M$ lies in the kernel of the {\it
Jacobi operator} $L=\Delta -2K$ of $M$ (here $K$ stands for the
Gaussian curvature function of $M$). Functions in the kernel of $L$
are called {\it Jacobi functions.} In his pioneering
paper~\cite{sh1}, Shiffman himself exploited this property when he
proved that if a minimal annulus $M$ is bounded by two circles in
parallel planes, then $M$ is foliated by circles in the intermediate
planes. For a surface $M\in {\cal M}$ and with the notation of items
1--4 before the statement of Theorem~\ref{teorKestim}, its Shiffman
function $S_M$ can be defined globally on $M=(\C /\langle i\rangle
)-\{ p_j,q_j\} _j$. By writing the local expression of $g$ around a
zero or pole, it is not difficult to check that $S_M$ is bounded
around the middle ends of $M$. Hence $S_M$ has a continuous
extension to the cylinder $\C /\langle i\rangle $ obtained after
attaching to $M$ its planar ends. Since $S_M$ solves the Jacobi
equation $LS_M=0$ and $L$ is of the form $\Delta +q$ for $q$ smooth
on $\C /\langle i\rangle $ (now $\Delta $ refers to the laplacian in
the flat metric on $\C /\langle i\rangle $),  then elliptic
regularity implies that $S_M$ extends smoothly to $\C /\langle
i\rangle $. Another key property of $S_M$ for every $M\in {\cal M}$
is that $S_M$ is bounded on $\C /\langle i\rangle $, as we next
explain.

\begin{definition}
\label{defgquasi} {\rm A meromorphic function $g\colon \C/\langle
i\rangle \to \C \cup \{ \infty \} $ is said to be {\it
quasiperiodic} if it satisfies the following two conditions:
\begin{enumerate}
\item There exists a constant $C>0$ such that the distance between
any two distinct points in $g^{-1}(\{ 0,\infty \} )\subset \C / \langle i\rangle $ is at least $C$
and given any $p\in g^{-1}(\{ 0,\infty \} )$, there exists at least one point in
$g^{-1}(\{ 0,\infty \} )-\{ p\} $ of distance less than $1/C$ from $p$.
\item For any divergent sequence $\{ z_k\} _k\subset
\C /\langle i\rangle $, there exists a subsequence of the
meromorphic functions $g_k(z)=g(z+z_k)$ which converges uniformly on
compact subsets of $\C /\langle i\rangle $ to a non-constant
meromorphic function $g_{\infty }\colon \C/\langle i\rangle \to \C
\cup \{ \infty \} $ (thus $g_{\infty }$ is also quasiperiodic).
\end{enumerate}
}
\end{definition}
Note that items {\it 5, 6} of Theorem~\ref{thm1} imply that the
(stereographically projected) Gauss map of a surface $M\in {\cal M}$
extends across the planar ends to a quasiperiodic meromorphic
function. This quasiperiodicity of $g$ together with equation
(\ref{eq:uShiffman}) clearly imply that $S_M$ is bounded on $\C
/\langle i\rangle $.

Let us return to the definition of the Shiffman function. Since the
conformal factor $\Lambda $ in (\ref{eq:uShiffman}) is a positive
function, we deduce that $S_M=0$ on $M$ if and only if $M$ is
foliated by pieces of circles and straight lines in horizontal
planes. In the nineteenth century, Riemann~\cite{ri2,ri1} classified
all minimal surfaces with such a foliation property: they reduce to
the examples appearing in Theorem~\ref{classthm}. Therefore, a
possible approach to proving Theorem~\ref{classthm} is to verify
that $S_M=0$. Instead of doing this directly, we will show that
given $M\in {\cal M}$, the Shiffman function $S_M$ can be
integrated, in the sense that there exists a one-parameter family
$\{ M_t\} _t\subset {\cal M}$ such that $M_0=M$ and the normal
component of the variational vector of this variation, when
restricted to each $M_t$, is exactly the Shiffman Jacobi function
$S_{M_{t}}$. In fact, in our arguments below we will need to
complexify our framework by considering complex-valued
Jacobi functions of the type $v+iv^*$, where $v^*$ is the Jacobi
conjugate function of a Jacobi function $v$ (see the definition
of ${\cal J}_{\C }(g)$ below) and {\it holomorphic} deformations
$t\mapsto (g_t,dz)$ of the Weierstrass data $(g,dz)$ of $M
\in {\cal M}$; this holomorphicity of the deformation will be crucial
in the proof of the linearity of the complex-valued Shiffman
function $S_M+iS_M^*\in {\cal J}_{\C }(g)$ of $M\in {\cal M}$
(Proposition~\ref{propos6.13}); in turn, the linearity of $S_M$
will imply that $M$ is periodic (Lemma~\ref{lema6.15}) and subsequently,
that $M$ is a Riemann minimal example (Proposition~\ref{propos7.3}).
In order for this framework to make sense, we
need discuss some technical issues.

First note that the Weierstrass pair of every $M\in {\cal M}$
depends solely on the Gauss map $g$ (not on the height differential
$dh=dz$), hence one can think of a deformation of $M$ inside ${\cal
M}$ as a certain kind of deformation of $g$ in an appropriate space
of {\it allowed} Gauss maps; namely, we consider the space of
meromorphic functions
\[
{\cal W}=\left\{ g\colon \C/\langle i\rangle  \to \C \cup \{ \infty \} \ \mbox{quasiperiodic\ } | \
(g)=\prod _{j\in \Z }p_j^2 q_j^{-2}%,\ \Re (p_j)<\Re (q_j)<\Re (p_{j+1})\  \forall j
\right\} ,
\]
where $(g)$ denotes the divisor of zeros and poles of $g$ on $\C
/\langle i\rangle $. By Theorem~\ref{thm1}, the Gauss map of every
$M\in {\cal M}$ lies in ${\cal W}$. Note that by
Definition~\ref{defgquasi}, any limit $g_{\infty }$ of a convergent
subsequence of $g_k(z)=g(z+z_k)$, with $\{ z_k\} _k\subset \C
/\langle i\rangle $ being a divergent sequence, satisfies that
$g_{\infty }$ lies in~${\cal W}$. The divisor of zeros $Z=\prod
_jp_j^2$ of a function $g\in {\cal W}$ is quasiperiodic, in the
sense that for every divergent sequence $\{ z_k\} _k\subset \C
/\langle i\rangle $, there exists a subsequence of $\{ Z+z_k\} _k$
which converges in the Hausdorff distance on compact subsets of $\C
/\langle i\rangle $ to a divisor $Z_{\infty }$ in $\C /\langle
i\rangle $, and a similar property holds for poles. Reciprocally,
two disjoint quasiperiodic divisors $Z=\prod _jp_j^2$, $P=\prod
_jq_j^2$ in $\C /\langle i\rangle $ define a unique quasiperiodic
meromorphic function $g$ (up to multiplicative non-zero constants)
whose principal divisor is $(g)=Z/P$. Existence is given in Douady
and Douady~\cite{Douady1}, while uniqueness follows from the
Liouville theorem together with fact that the function $f=g_1/g_2$,
where $g_1,g_2\in {\cal W}$ have $(g_1)=(g_2)$, is bounded on $\C
/\langle i\rangle $ (if $f$ were unbounded on $\C /\langle i\rangle
$, then we could find a divergent sequence $\{ z_k\} _k\subset \C
/\langle i\rangle $ such that $f(z_k)$ diverges; but the
quasiperiodicity of $g_1$ and $g_2$ implies that, after extracting a
subsequence, $f_k(z)=f(z+z_k)$ converges uniformly on compact
subsets of $\C $ to a meromorphic function $f_{\infty }\colon
\C/\langle i \rangle \to \C \cup \{ \infty \}$ which is not constant
infinity, thus $f_{\infty }$ has no poles by the Hurwitz theorem,
which contradicts that $f_k(0)=f(z_k)\to \infty $ as $k\to \infty
$). This association, to each $g\in {\cal W}$ of the quasiperiodic
set of its zeros and poles in $\C /\langle i\rangle $ together with
the value of $g$ at a prescribed point $z_0\in (\C /\langle i\rangle
)-g^{-1}(\{ 0,\infty \} )$, leads us to define the notion of
holomorphic dependence of a curve $t\mapsto g_t\in {\cal W}$:
\begin{definition}
{\rm Given $\ve >0$, we say that a curve $t\in \D (\ve ):=\{ t\in \C
\ | \ |t|<\ve \}\to g_t\in {\cal W}$ is {\it holomorphic} if the
corresponding functions $p_j(t),q_j(t),g_t(z_0)$ depend
holomorphically on~$t$. In this case, the function $\dot{g}\colon \C
/\langle i\rangle \to \C \cup \{ \infty \} $ given by $z\in \C
/\langle i\rangle \mapsto \left. \frac{d}{dt}\right| _{t=0}g_t(z)$
is meromorphic on $\C /\langle i\rangle $. We will call $\dot{g}$
the {\it infinitesimal deformation of $g$ associated to the curve
$t\mapsto g_t$}. }
\end{definition}

If $\dot{g}$ is the infinitesimal deformation of $g=g_0\in {\cal W}$
associated to the curve $t\mapsto g_t$ and $g$ has principal divisor
$(g)=\prod _jp_j^2q_j^{-2}$, then the principal divisor of $\dot{g}$
clearly satisfies $(\dot{g})\geq \prod _jp_jq_j^{-3}$.
%\end{equation}
%In particular, if $\dot{g}$ is constant, then $\dot{g}=0$.
The converse also holds: if $g\in {\cal W}$ and $f$ is a meromorphic
function on $\C /\langle i\rangle $ whose principal divisor verifies
$(f)\geq \prod _{j}p_jq_j^{-3}$, then $f$ can be proved to be the
infinitesimal deformation of $g$ associated to a holomorphic curve $
t\mapsto g_t\in {\cal W}$ with $g_0=g$. This leads to define a sort
of {\it tangent space at $g$ to ${\cal W}$,} which is formally the
complex linear space of infinitesimal deformations of $g$ associated
to holomorphic curves:
\begin{equation}
\label{eq:TgW}
T_g{\cal W}=\left\{ f\colon \C/\langle i\rangle \to \C \cup \{ \infty\} \mbox{ meromorphic }\ | \
(f)\geq \prod _jp_jq_j^{-3} \right\} .
\end{equation}
\begin{remark}
\label{remarkABC} Given $g\in {\cal W}$, the functions $g,g'\in
T_g{\cal W}$ are respectively the infinitesimal deformations at
$t=0$ associated to the holomorphic curves $t\mapsto (t+1)g(z)$,
$t\mapsto g(z+t)$.
\end{remark}

The usual period problem can be reformulated in ${\cal W}$ as
follows. Let $\g =\{ it\ | \ t\in [0,1]\}$ be the generator of the
homology of the cylinder $\C /\langle i\rangle $. Given $g\in {\cal
W}$, the list $\left( (\C /\langle i\rangle \right) -g^{-1}(\{
0,\infty \} ),g,dh=dz)$ is the Weierstrass data of a complete,
immersed minimal surface in $\R^3$ with individually embedded
horizontal planar ends at the zeros and poles of $g$ if and only if
the corresponding period problem~(\ref{eqWeiers2}) can be solved
(note that condition~(\ref{eqWeiers1}) holds automatically).
Equation (\ref{eqWeiers2}) can be encoded in the {\it period map}
$\mbox{Per}\colon {\cal W}\to \C ^2\times \C^{\Z}\times \C ^{\Z }$,
defined as
\begin{equation}
\label{eq:periodmap}
\mbox{Per}(g)=\left( \int _{\g }\frac{dz}{g}, \int _{\g }g\, dz,
\{ \mbox{Res}_{p_j}\left( \frac{dz}{g}\right) \} _j,
\{ \mbox{Res}_{q_j}(g\, dz)\} _j\right) .
\end{equation}
The subset of $g\in {\cal W}$ such that $(g,dz)$ solves the period
problem can be written as
\begin{equation}
\label{eq:Mimm}
{\cal M}_{\mbox{\footnotesize imm}}=\mbox{Per}^{-1}\{ (a,\overline{a},0,0)\ | \ a\in \C \} .
\end{equation}

\begin{definition}
{\rm
A {\it quasiperiodic, immersed minimal surface of Riemann type}
is a minimal surface $M\subset \R^3$ which admits a Weierstrass
pair of the form $(g,dz)$ on $(\C /\langle i\rangle )-g^{-1}(\{ 0,\infty \} )$
where $g$ lies in ${\cal M}_{\mbox{\rm \footnotesize imm}}$.
}
\end{definition}
Given $g\in {\cal W}$ (not necessarily in ${\cal M}_{\mbox{\rm
\footnotesize imm}}$), we denote by $N=\left(
\frac{2g}{|g|^2+1},\frac{|g|^2-1}{|g|^2+1}\right) \in \C \times \R
\equiv \R^3$ its related spherical ``Gauss map''. The Jacobi
equation $Lv=\Delta v+|\nabla N|^2v=0$ can be formally written as
\begin{equation}
\label{eq:Jacobi}
v_{z\overline{z}}+2\frac{|g'|^2}{(1+|g|^2)^2}v=0.
\end{equation}
A function $v\colon (\C /\langle i\rangle )-g^{-1}(\{ 0,\infty \}
)\to \R $ which satisfies equation~(\ref{eq:Jacobi}) is called a
{\it Jacobi function associated to $g$.} We will denote by ${\cal
J}(g)$ to the space of Jacobi functions associated to $g$. Since
$N:(\C /\langle i\rangle )-g^{-1}(\{ 0,\infty \} )\to \esf^2$ is
harmonic, then $\Delta N+|\nabla N|^2N=0$, which implies that the
functions $v=\langle N,a\rangle $ where $a\in \R^3 $, are always
Jacobi functions associated to $g$. Such functions will be called
{\it linear Jacobi functions.}

If $M\subset \R^3$ is a quasiperiodic, immersed minimal surface of
Riemann type (i.e. $g\in {\cal M}_{\mbox{\rm \footnotesize imm}}$),
then we define its {\it Shiffman function} $S_M$ by equation
(\ref{eq:uShiffman}). Note that the formula (\ref{eq:uShiffman})
makes sense for every $g\in {\cal W}$, but in order $S_M$ to be
bounded at each point in $g^{-1}(\{ 0,\infty \} )$ we need $g$ to
lie in ${\cal M}_{\mbox{\rm \footnotesize imm}}$.

Given a minimal surface $M$, let $B(N)\subset M$ denote the set of
branch points of the Gauss map $N$ of $M$. Given a Jacobi function
$v\colon M\to \R $, the formula
\begin{equation}
\label{eq:MontielRos}
X_v:=vN+\frac{1}{|N_z|^2}\{ v_zN_{\overline{z}}+v_{\overline{z}}N_z\}\colon
M-B(N)\rightarrow \R^3
\end{equation}
defines a branched minimal immersion (possibly constant) with the
same Gauss map $N$ as $M$, where $z$ is any local conformal
coordinate on $M$. The {\it support function}\footnote{Since
dilations of a minimal immersion $X\colon M\to \R^3$ are also
minimal and have the same Gauss map $N$ as $X$, the normal component
of the variational part of the deformation $X_t=\langle tX,N\rangle
$ is a Jacobi function. This function $\langle X,N \rangle $ is
called the {\it support function} of $X$.} of $X_v$ is $\langle
X_v,N\rangle =v$. The correspondence $v\mapsto X_v$ was studied by
Montiel and Ros~\cite{mro1} (see also Ejiri and Kotani~\cite{ek2}).
Linear Jacobi functions $v=\langle N,a\rangle $, $a\in \R^3$,
produce constant maps $X_v$ (and vice versa), and the law $v\mapsto
X_v$ can be considered to be a linear isomorphism from the linear
space of Jacobi functions on $M$ modulo the subspace of linear
Jacobi functions, onto the linear space of all branched minimal
immersions $X\colon M-B(N)\rightarrow \R^3$ with Gauss map $N$
modulo the constant maps.

The {\it conjugate Jacobi function} $v^*$ of a Jacobi function
$v\colon M\to \R $ is defined (locally) as the support function
$\langle (X_v)^*,N\rangle $ of the conjugate minimal immersion
$(X_v)^*$ (recall that such a conjugate minimal immersion is an
isometric minimal immersion of the underlying Riemannian surface,
whose coordinate functions are the harmonic conjugates to the ones
of $X_v$ and whose Gauss map is the same as for $M$). Since the
conjugate minimal immersion is defined modulo additive vector-valued
constants, then $v^*$ is defined up to additive scalar constants.
Furthermore, $v^*$ is globally well-defined precisely when $(X_v)^*$
is globally well-defined.

We consider the complex linear space
\[
{\cal J}_{\C }(g)=\{ v+iv^*\ | \ v\in {\cal J}(g) \mbox{ and $v^*$ is globally defined}\} .
\]
Thus, ${\cal J}_{\C }(g)$ is the space of support functions $\langle
X,N\rangle $ of holomorphic maps $X\colon (\C /\langle i\rangle
)-B(N)\to \C^3$ whose real and imaginary parts have differentials
orthogonal to $N=\left(
\frac{2g}{|g|^2+1},\frac{|g|^2-1}{|g|^2+1}\right) \in \C \times \R
\equiv \R^3$. The (complex) linear functions of $g$, $L_{\C }(g):=\{
\langle N,a\rangle \ | \ a\in \C^3\} $, form a complex linear
subspace of ${\cal J}_{\C }(g)$,

With the above framework, we are ready to state a precise
formulation for the integration of the Shiffman function of every
surface $M\in {\cal M}$, which will be proved in
Section~\ref{secKdV}.
\begin{theorem}
\label{thmintegrShiffman} Let $M$ be a quasiperiodic, immersed
minimal surface of Riemann type, with Gauss map $g:\C /\langle
i\rangle \to \C \cup \{ \infty \} $ and height differential $dz$.
Let $(g)=\prod _{j\in \Z }p_j^2q_j^{-2}$ be the principal divisor of
$g$ and let $z_0\in (\C /\langle i\rangle )-g^{-1}(\{ 0,\infty \}
)$. Then, the Shiffman function $S_M$ of $M$ can be {\it
holomorphically integrated} in the following sense: There exist
$\varepsilon >0$ and families $\{ p_j(t)\} _j$, $\{q_j(t)\}
_j\subset \C /\langle i\rangle $, $a(t)\in \C -\{ 0\} $  such that
\begin{description}
\item[{\it i)}] For each $j\in \Z $, the functions $t\in \D (\ve )\mapsto p_j(t)$,
$t\mapsto q_j(t)\in \C /\langle i\rangle $ are holomorphic with
$p_j(0)=p_j$, $q_j(0)=q_j$. Also, the function $t\mapsto a(t)$ is
holomorphic as well.
\item[{\it ii)}] For any $t\in \D (\varepsilon )$, the divisor
$\prod _jp_j(t)^2q_j(t)^{-2}$ defines an element $g_t\in {\cal
M}_{\mbox{\rm \footnotesize imm}}$ with $g_0=g$ and $g_t(z_0)=a(t)$.
Let $M_t$ be the quasiperiodic, immersed minimal surface of Riemann
type with Weierstrass pair  $(g_t,dz)$.
\item[{\it iii)}] For every $t\in \D (\ve )$, the derivative of $t\mapsto g_t$ with respect to $t$ equals
\[
\frac{d}{dt}g_t=\frac{i}{2}\left( g_t'''-3\frac{g_t'g_t''}{g_t}+
\frac{3}{2}\frac{(g_t')^3}{g_t^2}\right)  \quad \mbox{on }\C /\langle i\rangle .
\]
\end{description}
Furthermore if $M$ is embedded, then the surfaces $M_t$ are also
embedded for $|t|$ sufficiently small.
\end{theorem}

The condition {\it iii)} in Theorem~\ref{thmintegrShiffman} implies
that if $\Psi _t=\left( \frac{1}{2}(\frac{1}{g_t}-g_t),
\frac{i}{2}(\frac{1}{g_t}+g_t),1\right) dz$ is the usual Weierstrass
form associated to the pair $(g_t,z)$, then
\begin{equation}
\label{eqnueva}
\langle \left. \frac{d}{dt}\right| _t\int ^z\Psi _t,N_t\rangle
=-\frac{1}{2}(S_{M_t}+iS_{M_t}^*)+\langle a(t),N_t\rangle ,
\end{equation}
where $N_t=\left( \frac{2\Re (g_t)}{|g_t|^2+1},\frac{2\Im (g_t)}
{|g_t|^2+1},\frac{|g_t|^2-1}{|g_t|^2+1}\right)$ is the spherical
Gauss map of $M_t$ and $a(t)\in \C ^3$. In other words, the
normal component of the complex-valued variational field of $t
\mapsto (g_t,dz)$ is proportional to the (complex) Shiffman function $S_{M_t}+iS_{M_t}^*$
of $M_t$ modulo an additive linear function of $N_t$, for all $t$.
To understand why this is
true, we need two auxiliary results which will be also used in
Section~\ref{secKdV}.

\begin{proposition}
\label{propos4.8}
Given $g\in {\cal M}_{\mbox{\rm \footnotesize imm}}$, we have:
\begin{enumerate}
\item Let $h\colon \C /\langle i\rangle \to \C \cup \{ \infty \} $
be a meromorphic function which is a rational expression of $g$ and
its derivatives with respect to $z$ up to some order, such that
\begin{equation}
\label{eq:gpuntodeh}
\dot{g}(h)=\left(\frac{g^3h'}{2g'}\right) '
\end{equation}
belongs to $T_g{\cal W}$.
Then, the complex valued function
\begin{equation}
\label{eq:Jacobideh}
f(h)=\frac{g^2h'}{g'}+\frac{2gh}{1+|g|^2}
\end{equation}
lies in ${\cal J}_{\C }(g)$, is quasiperiodic\footnote{Given $g\in
{\cal W}$, a Jacobi function $v\in {\cal J}(g)$ is said to be {\it
quasiperiodic} if for every divergent sequence $\{ z_k\} _k\subset
\C /\langle i\rangle $, there exists a subsequence of the functions
$v_k(z)=v(z+z_k)$ which converges uniformly on compact subsets of
$(\C /\langle i\rangle )-g_{\infty }^{-1}(\{ 0,\infty \} )$ to a
function $v_{\infty }$, where $g_{\infty }\in {\cal W}$ is the limit
of (a subsequence of) $\{ g_k(z)=g(z+z_k)\} _k$, which exists since
$g$ is quasiperiodic. Note that $v_{\infty }\in {\cal J}(g_{\infty
})$ and that if $v_{\infty }$ is constant, then $v_{\infty }=0$. A
similar definition can be made exchanging ${\cal J}(g)$ by ${\cal
J}_{\C }(g)$.} and bounded on $\C /\langle i\rangle $. Furthermore,
for every closed curve $\G \subset \C /\langle i\rangle $,
\begin{equation}
\label{eq:periodconstant}
\int _{\G }\frac{\dot{g}(h)}{g^2}dz=\int _{\G }\dot{g}(h)\, dz=0.
\end{equation}
\item Reciprocally, if $\dot{g}\in T_g{\cal W}$ satisfies {\rm (\ref{eq:periodconstant}),}
then there exists a meromorphic function $h$ on $\C /\langle i\rangle $
such that (\ref{eq:gpuntodeh}) holds.
\end{enumerate}
\end{proposition}
{\it Sketch of proof.}
To prove item {\it 1,} assume $\dot{g}(h)\in T_g{\cal W}$. Then
there exists a holomorphic curve $t\mapsto g_t\in {\cal W}$
such that $g_0=g$ and $\left. \frac{d}{dt} \right| _{t=0}g_t=\dot{g}(h)$. Therefore
$\langle \left. \frac{d}{dt}\right| _0\int ^z\Psi _t,N\rangle
\in {\cal J}_{\C }(g)$, where
$\Psi _t=\left( \frac{1}{2}(\frac{1}{g_t}-g_t),\frac{i}{2}(\frac{1}{g_t}+g_t),1\right) dz$ and
$N=\left( \frac{2\Re (g)}{|g|^2+1},\frac{2\Im (g)}{|g|^2+1},\frac{|g|^2-1}{|g|^2+1}\right)$.
A simple calculation gives
\begin{equation}
\label{eq:primitivas}
\int ^z\frac{\dot{g}(h)}{g^2}dz=\frac{gh'}{2g'}+h,\qquad
\int ^z\dot{g}(h)\, dz=\frac{g^3h'}{2g'}
\end{equation}
up to additive complex numbers, from where it is straightforward to get
\begin{equation}
\label{eq:propos4.13new}
\langle \left. \frac{d}{dt}\right| _0\int ^z\Psi _t,N\rangle =
-\frac{1}{2}f(h)+\langle a,N\rangle ,
\end{equation}
for some $a\in \C^3$. Equation (\ref{eq:propos4.13new}) implies that
$f\in {\cal J}_{\C }(g)$. The quasiperiodicity of $f(h)$ follows
directly from the quasiperiodicity of $g$ (recall that $h$ is a
rational function of $g$ and its derivatives). In order to prove
that $f(h)$ is bounded on $\C /\langle i\rangle $, one first shows
that $f(h)$ is bounded around every zero and pole of $g$ and around
every zero of $g'$ (this is a direct computation using the behavior
of $\dot{g}(h)$ at zeros and poles of $g$, see the definition of
$T_g{\cal W}$). Then the boundedness of $f(h)$ in $\C /\langle
i\rangle $ can be reduced to work away from the discrete set
$g^{-1}(\{ 0,\infty \} )\cup (g')^{-1}(0)$, where it follows
directly because $g$ is quasiperiodic, $h$ is a rational expression
of $g$ and its derivatives, and $f$ is given in terms of $g,h$ by
the formula (\ref{eq:Jacobideh}). Therefore $f(h)$ is bounded on $\C
/\langle i\rangle $. Finally, (\ref{eq:periodconstant}) is a direct
consequence of (\ref{eq:primitivas}). This gives item {\it 1} of the
proposition. Concerning item {\it 2}, equation (\ref{eq:primitivas})
together with the hypothesis (\ref{eq:periodconstant}) allow us to
find a meromorphic function $h$ on $\C /\langle i \rangle $ such
that (\ref{eq:gpuntodeh}) holds.
{\hfill\penalty10000\raisebox{-.09em}{$\Box$}\par\medskip}

\begin{remark}
\label{remark6.12}
%{\rm
An important interpretation of equation (\ref{eq:periodconstant}) is
that $\dot{g}(h)$ lies in the kernel of the differential $d\,
\mbox{\rm Per}_g$ of the period map at $g\in {\cal M}_{\mbox{\rm
\footnotesize imm}}$, see equation (\ref{eq:periodmap}).
%}
\end{remark}
It is clarifying to analyze some examples of $\dot{g}(h)$ and $f(h)$
coming from particular choices of $h$ in
Proposition~\ref{propos4.8}. If one takes $h=c_1+\frac{c_2}{g^2}$ in
(\ref{eq:gpuntodeh}) with $c_1,c_2\in \C $, then $\dot{g}(h)=0$ (and
vice versa). In this case, $f(h)$ is a complex linear combination of
$\frac{g}{1+|g|^2},\frac{\overline{g}}{1+|g|^2}$, which can be
viewed as a horizontal linear function of the ``Gauss map'' $g$. If
one takes $h=\frac{1}{g}$ in (\ref{eq:gpuntodeh}), then
$\dot{g}(h)=-\frac{1}{2}g'$ and $f(h)=\frac{1-|g|^2}{1+|g|^2}$,
which is a vertical linear function of $g$. A slightly more subtle
choice of $h$ explains the nature of condition {\it iii)} in
Theorem~\ref{thmintegrShiffman}:

\begin{corollary}
\label{corol4.10}
Let $M$ be a quasiperiodic, immersed minimal surface of Riemann type.
Then, its Shiffman function $S_M$ given by (\ref{eq:uShiffman}) admits
a globally defined conjugate Jacobi $S_M^*$, and
$S_M+iS_M^*=f$ is given by equation (\ref{eq:Jacobideh})
for
\begin{equation}
\label{Shiffmandeh}
h=h_S=\frac{i}{2}\frac{(g')^2}{g^3}.
\end{equation}
In particular:
\begin{enumerate}
\item Both $S_M,S_M^*$ are bounded on the
cylinder $M\cup g^{-1}(\{ 0,\infty \} )$.
\item The corresponding infinitesimal deformation
$\dot{g}_S=\dot{g}(h_S)\in T_g{\cal W}$ is given by
\begin{equation}
\label{gpuntodeShiffman}
\dot{g}_S=\frac{i}{2}\left( g'''-3\frac{g'g''}{g}+\frac{3}{2}\frac{(g')^3}{g^2}\right) .
\end{equation}
\item If $\dot{g}_S=0$ on $M$, then both $S_M,S_M^*$ are linear.
\end{enumerate}
\end{corollary}
\begin{proof}
Note that $h$ defined by (\ref{Shiffmandeh}) is a rational
expression of $g$ and $g'$. A direct computation gives that plugging
(\ref{Shiffmandeh}) into (\ref{eq:gpuntodeh}) we obtain
(\ref{gpuntodeShiffman}), and that this last expression has the
correct behavior at $g^{-1}(\{ 0,\infty \})$ expressed in
(\ref{eq:TgW}), thus $\dot{g}_S\in T_g{\cal W}$. Now item {\it 1} of
this corollary follows from Proposition~\ref{propos4.8}. Regarding
item~{\it 3}, if $\dot{g}(h_S)=0$ then (\ref{eq:gpuntodeh}) gives
$h_S=b-\frac{c}{g^2}$ for $b,c\in \C $. Substituting in
(\ref{eq:Jacobideh}) we have
$S_M+iS_M^*=2c\frac{\overline{g}}{1+|g|^2}+2b\frac{g}{1+|g|^2}$.
Hence, both $S_M,S_M^*$ are linear.
\end{proof}

We said above that a way of finishing the proof of
Theorem~\ref{classthm} is to show that if $M\in {\cal M}$, then
$S_M=0$. Next we explain that it is enough to demonstrate a weaker
condition, namely that $S_M$ is linear. This is the goal of
Lemma~\ref{lema6.15} and Proposition~\ref{propos7.3} below, which do
not use Theorem~\ref{thmintegrShiffman} in their
proofs\footnote{Theorem~\ref{thmintegrShiffman} will be used when
proving the linearity of $S_M$ for every $M\in {\cal M}$, see
Proposition~\ref{propos6.13}.}.
\begin{lemma}
\label{lema6.15} Suppose that the Shiffman function $S_M$ of a
quasiperiodic, immersed minimal surface of Riemann type $M\subset
\R^3$, is linear. Then, $M$ is singly-periodic (invariant by a
translation), and its smallest orientable quotient surface
$M/\langle v\rangle\subset \R^3/\langle v\rangle$ is a properly
immersed minimal torus with two planar ends and total curvature
$-8\pi $, where $v\in \R^3-\{ 0\} $ is a translation vector of $M$.
\end{lemma}
{\it Sketch of proof.} The argument uses the notion of conjugate
Jacobi functions and the Montiel-Ros correspondence between Jacobi
functions and branched minimal immersions explained before the
statement of Theorem~\ref{thmintegrShiffman}. Suppose that $M$ is a
minimal surface verifying the hypotheses in the lemma, and let $N$
be its Gauss map. Since the Shiffman function $S_M$ of $M$ is
assumed to be linear, then its related branched minimal immersion
$X_v$ is constant, where $v=S_M$. Thus the conjugate minimal
immersion $(X_v)^*$ is also constant, which implies that the
conjugate Jacobi function $S_M^*$ of $S_M$ exists globally and it is
linear. The linearity of the complex valued Jacobi function
$S_M+iS_M^*$ together with equation~(\ref{eq:uShiffman}) imply that
there exists $a\in \C ^3$ such that
\begin{equation}
\label{eq:lema7.2A}
\frac{3}{2}\left( \frac{g'}{g}\right) ^2-\frac{g''}{g}
 -\frac{1}{1+|g|^2}\left( \frac{g'}{g}\right) ^2=\langle N,a\rangle .
\end{equation}
After writing the last right-hand-side in terms of $g$ and
manipulating algebraically the resulting equation, one arrives to
the following ODE for $g$:
\[
\overline{g}\left( \frac{3}{2}\frac{(g')^2}{g}-g''-B-a_3g\right) =
\frac{g''}{g}-\frac{1}{2}\left( \frac{g'}{g}\right) ^2+Ag-a_3,
\]
where $a=(a_1,a_2,a_3)$, $2a_1=A+B$ and $2a_2=i(A-B)$.
Since $g$ is holomorphic and not constant, we deduce that
\[
\frac{3}{2}\frac{(g')^2}{g}-g''-B-a_3g=0,\qquad
\frac{g''}{g}-\frac{1}{2}\left( \frac{g'}{g}\right) ^2+Ag-a_3=0.
\]
After elimination of $g''$ in both equations, we have $(g')^2=g(-A
g^2+2a_3 g+B)$. This equation implies that the map $\pi =(g,g')$ is
a possibly branched holomorphic covering from the cylinder $M\cup \{
\mbox{planar ends}\} \equiv \C /\langle i\rangle $ onto the compact
Riemann surface $\Sigma =\{ (\xi ,w)\in (\C\cup \{ \infty \} )^2\ |
\ w^2=\xi(-A\xi^2+2a_3\xi+B)\} $. Clearly, $\Sigma $ is either a
sphere or a torus.

Suppose for the moment that $\Sigma $ is a sphere. Then, consider
the meromorphic differential $\frac{d\xi }{w}$ on $\Sigma $, whose
pullback by $\pi $ is $\pi ^*(\frac{d\xi }{w})=\frac{dg}{g'}=dz$.
Given a pole $P\in \Sigma $ of $\frac{d\xi }{w}$, choose a point
$z_0\in \C /\langle i\rangle $ such that $\pi (z_0)=P$. The residue
of $\frac{d\xi }{w}$ at $P$ can be computed as the integral of
$\frac{d\xi }{w}$ along a small closed curve $\G _P\subset \Sigma $
which winds once around $P$. After lifting $\G _P$ through $\pi $
locally around $z_0$ we obtain a closed curve $\widetilde{\G
}_P\subset \C /\langle i\rangle $ which winds a positive integer
number of times  around $z_0$, depending on the branching order of
$\pi $ at $z_0$. Hence the residue of $\frac{d\xi }{w}$ at $P$
equals a positive integer multiple of the residue of $dz$ at $z_0$,
which is zero. Therefore, $\frac{d\xi }{w}$ has residue zero at all
its poles, and so, it is exact on $\Sigma $. This implies that $dz$
is also exact on $\C /\langle i\rangle $, which is impossible. This
contradiction proves that $\Sigma $ is a torus.

Finally, consider on $\Sigma $ the Weierstrass pair $\left( g_1(\xi
,w)=\xi ,dh_1=\frac{d\xi }{w}\right) $. The metric associated to
this pair is $\left( \frac{1}{2}(|\xi |+|\xi |^{-1})\frac{|d\xi
|}{|w|}\right) ^2$, which can be easily proven to be positive
definite and complete in $\Sigma -\{ (0,0),(\infty ,\infty )\} $.
Note that $g_1\circ \pi =g$ and $\pi ^*(\frac{d\xi }{w})=dz$. This
implies that the Weierstrass pair $(g,dz)$ of $M$ can be induced on
the twice punctured torus $\Sigma -\{ (0,0),(\infty ,\infty )\} $,
from where the lemma follows easily.
{\hfill\penalty10000\raisebox{-.09em}{$\Box$}\par\medskip}

\begin{proposition}
\label{propos7.3}
If the Shiffman function $S_M$ of an embedded surface $M\in {\cal M}$ is linear,
then $M$ is a Riemann minimal example.
\end{proposition}
{\it Sketch of proof.} Let ${\cal M}_1\subset {\cal M}$ be the
subset of surfaces which are singly-periodic and their smallest
orientable quotient is a properly embedded, twice-punctured minimal
torus in a quotient of $\R^3$ by a translation. By
Lemma~\ref{lema6.15}, our proposition reduces to proving that ${\cal
M}_1$ coincides with the family of Riemann minimal examples. This is
a particular case of the main theorem in the 1998 paper~\cite{mpr1}
by Meeks, P\'erez and Ros, and we now give some ideas of its proof.

Consider the flux map $h\colon {\cal M}_1\to (0,\infty )$ which
associates to each surface $M\in {\cal M}_1$ the positive number
$h(M)$ such that $F=(h(M),0,1)$ is the flux vector of $M$. The goal
is to show that $h$ is a homeomorphism if we endow ${\cal M}_1$ with
the uniform topology on compact sets. This is a consequence of the
following steps.
\begin{description}
\item[(A)] The family ${\cal R}\subset {\cal M}_1$ of Riemann minimal examples is a
connected component of ${\cal M}_1$.
\item[(B)] $h:{\cal M}_1\to (0,\infty )$ is a proper  map.
\item[(C)] $h:{\cal M}_1\to (0,\infty )$ is an open map.
\item[(D)] There exists $\ve >0$ such that if $M\in {\cal M}_1$ and $h(M)\in (0,\ve )$,
then $M\in {\cal R}$.
\end{description}
We briefly indicate how the above four steps can be proved. In order
to show {\bf (A)}, first note that ${\cal R}$ is a path-connected
set in ${\cal M}_1$ (by construction), and ${\cal R}$ is closed in
${\cal M}_1$ (since Riemann minimal examples are characterized in
${\cal M}_1$ by being foliated by circles and lines in parallel
planes, a condition which is preserved under limits in the uniform
topology on compact sets). Hence {\bf (A)} will follow if we check
that ${\cal R}$ is open in ${\cal M}_1$. We will give two different
reasons of why this openness holds. The first one is based on the
fact that each Riemann minimal example $R$ is non-degenerate in the
sense that its space of periodic, bounded Jacobi functions reduces
to the linear Jacobi functions (as follows from Montiel and
Ros~\cite{mro1} since all the branch values of the Gauss map of $R$
lie on a spherical equator), together with the property that being
non-degenerate is an open condition (see P\'erez~\cite{perez3}). The
second proof is by contradiction: otherwise we find surfaces
$\{\Sigma _n\} _n\subset {\cal M}_1- {\cal R}$ which converge on
compact subsets of $\R^3$ to some Riemann minimal example $R\in
{\cal R}$. Since $\Sigma _n\notin {\cal R}$, then the Shiffman
functions $S_{\Sigma _n}$ of the $\Sigma _n$ cannot be identically
zero. After normalizing $\widehat{S}_{\Sigma _n}=\frac{1}{\sup
_{\Sigma _n}|S_{\Sigma _n}|}S_{\Sigma _n}$ (recall that $S_{\Sigma
_n}$ is bounded on $\Sigma _n$, hence $\sup _{\Sigma _n}|S_{\Sigma
_n}|$ exists), we have a sequence of bounded Jacobi functions, for
which a subsequence converges to a periodic Jacobi function
$\widehat{S}_{\infty }$ on $R$. Recall that the zeros of the
Shiffman function of a surface are the critical points of the
curvature of the horizontal level sections of this surface. By the
Four Vertex Theorem, each $\widehat{S}_{\Sigma _n}$ has at least
four zeros on each compact horizontal section of ${\Sigma _n}$
(counted with multiplicity), and thus, the same holds for
$\widehat{S}_{\infty }$ on each compact horizontal section of $R$.
On the other hand, $\widehat{S}_{\infty }$ is a periodic, bounded
Jacobi function on $R$, hence it is linear. Now one finishes by
checking that a linear Jacobi function on a Riemann minimal example
has at most two zeros on each horizontal circle, and that these
zeros are simple (see Assertion~9.5 in Meeks, P\'erez and
Ros~\cite{mpr6}).

The properness of $h$ in step {\bf (B)} is a direct consequence of
the curvature estimates in Theorem~\ref{teorKestim} (actually one
only needs a weaker version of this curvature estimates, namely for
{\it singly-periodic} minimal planar domains, which was proved for the
first time in Theorem~4 of~\cite{mpr1}).

The openness property {\bf (C)} can be proved as follows. Consider
the space ${\cal W}_1=\{ (\Sigma ,g,[\a ])\} $ where $\Sigma $ is a
compact Riemann surface of genus one, $g\colon \Sigma \to \C \cup \{
\infty \} $ is a degree two meromorphic function with a double zero
$p$ and a double pole $q$, and $[\a ]$ is a homology class in
$H_1(\Sigma -\{ p,q\} ,\Z )$ which is non-trivial in $H_1(\Sigma ,\Z
)$. We denote the elements in ${\cal W}_1$ simply by $g$. There
exists a natural inclusion ${\cal W}_1\subset {\cal W}$ that unwraps
both $\Sigma $ and $g$ so that $[\a ]$ remains in the lifted
cylinder that cyclically covers $\Sigma $ (elements in ${\cal
W}_1$ can be viewed as singly-periodic elements in ${\cal W}$). The
space ${\cal W}_1$ is a two-dimensional complex manifold with local
charts given by $g\mapsto (a_1+a_2,a_1\cdot a_2)$, where $a_1,a_2\in
\C-\{ 0\} $ are the (possibly equal) branch values of $g\in {\cal
W}_1$ close to a given element $g_0\in {\cal W}_1$ (in a chart, we
can forget about the homology class $[\a ]$
associated to $g$ after
identification with that of $g_0$). Given $g\in {\cal W}_1$, we
associate a unique holomorphic differential $\phi $ on $\Sigma $ by
the equation $\int _{\a }\phi =i$. Consider the related {\it period
map} $\mbox{Per}_1\colon {\cal W}_1\to \C^2$ given by
$\mbox{Per}_1(g)=\left( \int _{\a }\frac{\phi }{g},\int _{\a }g\,
\phi \right) $ (compare with (\ref{eq:periodmap}) and recall that we
are considering singly-periodic minimal surfaces with only two ends,
so in order to express the period problem we do not need to impose
the vanishing of the residue of $\frac{\phi }{g}$ at $p$ or of $g\,
\phi $ at $q$ because the residues of a meromorphic differential on
a compact Riemann surface add up to zero). Then, the space of
elements $g\in {\cal W}_1$ such that $(g,\phi )$ is the Weierstrass
pair of an {\it immersed} minimal surface are ${\cal
M}_{1}^{\mbox{\footnotesize imm}}= \mbox{Per}_1^{-1}(\{
(a,\overline{a})\ | \ a\in \C \} )$. Since $\mbox{Per}_1$ is a
holomorphic map, for $a\in \C $ fixed the set ${\cal
M}_1^{\mbox{\footnotesize imm}}(a)=\mbox{Per}_1^{-1}(a)$ is a
complex analytic subvariety of ${\cal W}_1$. As the limit of
embedded surfaces is embedded, the subset ${\cal M}_1 \subset {\cal
M}_{1}^{\mbox{\footnotesize imm}}$ of {\it embedded} surfaces is
closed in ${\cal M}_{1}^{\mbox{\footnotesize imm}}$. An application
the maximum principle at infinity (Theorem~\ref{thmMaxprin}) gives
that ${\cal M}_1$ is also open in ${\cal M}_{1}^{\mbox{\footnotesize
imm}}$. In particular, the set ${\cal M}_1(a)={\cal
M}_{1}^{\mbox{\footnotesize imm}} (a)\cap {\cal M}_1$ is a complex
analytic subvariety of ${\cal W}_1$. By the uniform curvature
estimates for {\it embedded} (possibly non-singly-periodic) surfaces
in Theorem~\ref{teorKestim} and subsequent uniform local area
estimates, we deduce that ${\cal M}_1(a)$ is compact. Since the only
compact, complex analytic subvarieties of ${\cal W}_1$ are finite
sets, we deduce that ${\cal M}_1(a)$ is finite. Thus, given $M\in
{\cal M}_1$, there exists an open neighborhood $U$ of $M$ in ${\cal
W}_1$ such that $U\cap {\cal M}_{1}(a)= U\cap {\cal
M}_{1}^{\mbox{\footnotesize imm}}(a)=\{ M\} $. In this setting, the
openness theorem for finite holomorphic maps (Chapter 5.2 of
Griffiths and Harris~\cite{GriHar}) gives that $\mbox{Per}_1$ is an
open map locally around~$M$. Finally, the formula
$\mbox{Per}_1(g)=(i\overline{h},-ih,0,0)$ where $h=h(M)$ and $g$ is
the Gauss map of $M\in {\cal M}_1$, relates the period map
$\mbox{Per}_1$ with the flux map $h\colon {\cal M}_1\to (0,\infty
)$, from where the desired openness for $h$ follows.

Finally, to prove property {\bf (D)} one firstly studies the
boundary points of ${\cal W}_1$ which correspond to limits of
surfaces $M\in {\cal M}_1$ such that $h(M)$ converges to zero (in
this case the surfaces $M\in {\cal M}_1$ degenerate into a singular
object consisting of a stack of two vertical catenoids). In second
place one parameterizes ${\cal W}_1$ around the related boundary
points, and notices that these singular objects form an analytic
subvariety of a polydisc in $\C^2$ (say centered at the origin,
which corresponds to the stack of two catenoids). Furthermore, not
only ${\cal W}_1$ but also the Period map Per$_1$ can be
holomorphically extended through these singular points, and the
Jacobian of the extended Period map at the origin is a bijection.
Now ${\bf (D)}$ follows from the Inverse Function Theorem\footnote{A
simpler proof of property {\bf (D)} which does not use the
periodicity of the surfaces can be found in Theorem~36 of Meeks and
P\'erez~\cite{mpe1}.} applied to the extension of Per$_1$ to the
polydisc. {\hfill\penalty10000\raisebox{-.09em}{$\Box$}\par\medskip}

To finish the classification of the properly embedded minimal
surfaces with two limit ends modulo Theorem~\ref{thmintegrShiffman},
we only need to show the following result.

\begin{proposition}
\label{propos6.13}
The Shiffman function of every surface $M\in {\cal M}$ is linear.
\end{proposition}
\begin{proof}
Consider a surface $M\in {\cal M}$.
According to the notation at the beginning of this Section~\ref{sec2limitends}, the
heights of the planar ends of $M\in {\cal M}$ are
\begin{equation}
\label{eq:heightends}
\ldots <\Re (p_{-1})<\Re (q_{-1})<\Re (p_0)<\Re (q_0)<\Re (p_1)<\Re (q_1)<\ldots
\end{equation}
Consider the positive functions $h_j\colon {\cal M}\to \R $ for $j\in \N$ given by
\[
h_1(M)=\Re (q_0-p_0),\quad h_2(M)=\Re (p_0-q_{-1}),\quad h_3(M)=\Re (p_1-p_0),\quad h_4(M)=\Re (p_0-p_{-1})\ldots
\]
thus the functions $h_j$ associate to each surface $M\in {\cal M}$
the absolute value of the relative height of its planar ends with
respect to one of these ends, namely $p_0$. First note that $h_j$ is
continuous (we endow ${\cal M}$ with the topology of the uniform
convergence on compact subsets of $\R^3$).

Let ${\cal M}_F$ be the set of surfaces in ${\cal M}$ with the same
flux vector $F$ as $M$. By the uniform curvature estimates in
Theorem~\ref{teorKestim} and subsequent uniform local area
estimates, ${\cal M}_F$ is compact. Thus, the set ${\cal M}_F(1)= \{
M'\in {\cal M}_F \ | \ h_1(M')=\max _{{\cal M}_F}h_1\} $ is
non-empty. Consider the restriction of $h_2$ to ${\cal M}_F(1)$, and
as before maximize $h_2$ on ${\cal M}_F(1)$, hence the set ${\cal
M}_F(2)=\{ M'\in {\cal M}_F(1)\ | \ h_2(M')=\max _{{\cal
M}_F(1)}(h_2)\} $ is non-empty. Repeating the argument, induction
lets us for each $j\in \N $ maximize $h_{j+1}$ in ${\cal M}_F(j)=\{
M'\in {\cal M}_F(j-1)\ | \ h_j(M')=\max _{{\cal M}_F(j-1)}(h_j)\}
\neq \mbox{\O}$. Since the compact subsets ${\cal M}_F(j)$ satisfy
${\cal M}_F(j)\supset {\cal M}_F(j+1)$ for all $j$, this collection
of closed sets satisfies the finite intersection property. By the
compactness of ${\cal M}_F$,  we conclude that $\bigcap _{j\in \N}
{\cal M}_F(j)\neq \mbox{\O }$. Thus there exists a surface $M_{\max
}\in {\cal M}_F$ that maximizes each of the functions $h_{j+1}$ in
${\cal M}_F(j)$ for all $j\geq 1$. In the same way, we find a
surface $M_{\min }\in {\cal M}_F$ that minimizes the functions
$h_{j+1}$ on ${\cal M}_F(j)$ for all $j\geq 1$.

Next we prove that if a surface $M_0\in {\cal M}_F$ maximizes all
the functions $h_j$ as in the previous paragraph, then its Shiffman
function $S_{M_0}$ is linear (for minimizing surfaces the argument is
similar). By our assumption, Theorem~\ref{thmintegrShiffman} holds
for $M_0$. Hence there exists a curve of functions $t\in \D (\ve
)\mapsto g_t\in {\cal M}_{\mbox{\rm \footnotesize imm}} \subset
{\cal W}$ such that the zeros  $p_j(t)$ and poles $q_j(t)$ of $g_t$
depend holomorphically on $t$, satisfying items {\it i), ii)} and
{\it iii)} of Theorem~\ref{thmintegrShiffman} for $M=M_0$. With the
notation of that theorem, let $\psi _t\colon (\C /\langle i\rangle
)-\{ p_j(t),q_j(t)\} _j\to \R^3$ be the parametrization of $M_t$
given by $\psi _t(z)=\Re \int _{z_0}^z\Psi _t$, where
$z_0\in (\C /\langle i\rangle )- \{ p_j(t),q_j(t)\ | \ j\in \Z, |t|<\ve \} $
and $\Psi _t=\left( \frac{1}{2}(\frac{1}{g_t}-g _t),
\frac{i}{2}(\frac{1}{g_t}+g_t),1\right) \, dz$.
Item {\it iii)} of Theorem~\ref{thmintegrShiffman} together with
Corollary~\ref{corol4.10} and equation~(\ref{eq:propos4.13new})
imply that equation (\ref{eqnueva}) holds, i.e.
$\langle \left. \frac{d}{dt}\right| _t\int ^z\Psi _t,N_t\rangle $
is, up to a multiplicative constant,
the complex valued Shiffman function $S_{M_t}+iS_{M_t}^*$ of $M_t$
plus a linear function of the Gauss map $N_t$ of $M_t$. Furthermore,
equation (\ref{eq:periodconstant}) applied to
$\dot{g}_t=\frac{d}{dt}g_t$ gives that $\dot{g}_t$ lies in the
kernel of the differential of the period map Per at $g_t$, for all
$t$, see Remark~\ref{remark6.12}. Hence, the (complex) period map
remains constant along $t\mapsto g_t$, which implies that $M_t\in
{\cal M}_F$ for all $t$. This condition gives that the harmonic
function $t\in \D (\ve  )\mapsto h_1(M_t)=\Re (q_0(t)-p_0(t))$
attains a maximum at $t=0$, so it is constant. Since the function
$t\in \D(\ve )\mapsto q_0(t)-p_0(t)$ is holomorphic with constant
real part, then $q_0(t)-p_0(t)$ does not depend on $t$. The same
argument applies to each function $t\mapsto h_j(M_t)$ with $j\in \N
$, hence for every $t$ all the planar ends $p_j(t),q_j(t)$ of $M_t$
are placed at
\[
p_j(t)=p_0(t)+p_j-p_0,\quad q_j(t)=p_0(t)+q_j-p_0.
\]
Geometrically, this means that the maps $\psi _t$ coincide with
$\psi _0$ up to translations in the parameter domain and in $\R^3$.
Therefore, the normal part of the variational field of $t\mapsto
\psi _t$ is linear, and thus the Shiffman function of $M_0$ is
linear.

Finally, we prove that the Shiffman function of every $M\in {\cal
M}$ is linear. Given $M\in {\cal M}$, let $F=(h,0,1)$ be its flux
vector. By the arguments above, we find embedded minimal surfaces
$M_{\max }, M_{\min }\in {\cal M}_F$ such that $M_{\max }$ (resp.
$M_{\min }$) maximizes (resp. minimizes) all the functions $h_j$ in
the above sense, $j\in \N $. Furthermore, the arguments in the last
paragraph imply that the Shiffman functions of $M_{\max }, M_{\min
}$ are linear. By Proposition~\ref{propos7.3}, both $M_{\max
},M_{\min }$ are Riemann minimal examples. But ${\cal M}_F$ contains
at most one Riemann minimal example, since the flux is a parameter
for the space of Riemann examples. This implies that $M_{\max
}=M_{\min }$. On the other hand, the vertical distance between the
ends $p_0,q_0$ of $M$ is bounded above (resp. below) by the distance
between the corresponding ends of $M_{\max }$ (resp. of $M_{\min
}$). So, the vertical distance between the ends $p_0,q_0$ of $M$ is
maximum, or equivalently, $M$ maximizes $h_1$ on ${\cal M}_F$.
Analogously, $M$ maximizes all the functions $h_j$ and so, its
Shiffman function $S_M$ is linear.
\end{proof}

\section{Infinitely many ends III: The KdV equation.}
\label{secKdV} In this section we will prove
Theorem~\ref{thmintegrShiffman}, which finishes our classification
of the properly embedded minimal planar domains in $\R^3$. The main
tool in this proof is the {\it Korteweg-de Vries equation} (KdV) and
its {\it hierarchy,} two cornerstones in integrable systems theory.

\subsection{Relationship between the KdV equation and the Shiffman function.}
Recall that by Corollary~\ref{corol4.10}, the Shiffman function
$S_M$ of a quasiperiodic immersed minimal surface of Riemann type
$M$ admits a globally defined conjugate Jacobi function $S_M^*$ on
$\C /\langle i\rangle $, and $f=f_S=S_M+iS_M^*\in {\cal J}_{\C }(g)$
is given by equation~(\ref{eq:Jacobideh}) with
$h=h_S=\frac{i}{2}\frac{(g')^2}{g^3}$. This function $h$, when
plugged into equation (\ref{eq:gpuntodeh}), produces the
quasiperiodic meromorphic function $\dot{g}_S\in T_g{\cal W}$ given
by equation (\ref{gpuntodeShiffman}). By definition of $T_g{\cal
W}$, $\dot{g}_S$ is the derivative \underline{at $t=0$} of a
holomorphic curve $t\in \D (\ve )\mapsto g_t\in {\cal W}$ with
$g_0=g$. What we want to prove is that such a holomorphic curve
$t\mapsto g_t$ can be chosen so that \underline{for all $t$}, the
pair $(g_t,dz)$ is the Weierstrass data of a minimal surface $M_t\in
{\cal M}$ and
\begin{equation}
\label{gpuntodeShiffmant}
\left. \frac{d}{dt}\right| _{t}g_t=\frac{i}{2}\left( g_t'''-3\frac{g_t'g_t''}{g_t}+
\frac{3}{2}\frac{(g_t')^3}{g_t^2}\right) .
\end{equation}
Therefore, one could think of the above problem as finding an
integral curve of a vector field in a manifold. Unfortunately, this
approach is unsatisfactory: on one hand, the right-hand-side of
(\ref{gpuntodeShiffman}) does not give an element of $T_g{\cal W}$
\underline{for all $g\in {\cal W}$} (it can be proved that
$g'''-3\frac{g'g''}{g}+\frac{3}{2}\frac{(g')^3}{g^2}\in T_g{\cal W}$
provided that $(g,dz)$ closes periods at the zeros and poles of $g$;
in particular, this holds when $g\in {\cal M}_{\mbox{\rm
\footnotesize imm}}$). On the other hand, ${\cal M}_{\mbox{\rm
\footnotesize imm}}$ is not known to have a manifold structure, so
the general theory of integral curves of vector fields does not
apply to ${\cal M}_{\mbox{\rm \footnotesize imm}}$.

A more satisfactory viewpoint  is to consider
(\ref{gpuntodeShiffmant}) as an evolution equation with respect to
the complex time $t$. Therefore, one could apply general PDE
techniques to find solutions $g_t=g_t(z)$ of this initial value
problem, only defined {\it a priori} locally around a point $z_0\in
(\C /\langle i \rangle )-g^{-1}(\{ 0,\infty \} )$ with the initial
condition $g_0=g$. Such solutions are not necessarily global on $\C
/ \langle i\rangle $, might develop essential singularities (we need
the $g_t$ to be {\it meromorphic} in $z$ on $\C / \langle i\rangle $
in order $g_t$ to lie in ${\cal W}$), and even if $g_t$ were
meromorphic on $\C /\langle i\rangle $, it is not clear that $g_t$
would have only double zeros and poles and other properties
necessary to give rise to minimal surfaces $M_t$ in ${\cal M}$.
(Meromorphic) KdV theory will be crucial to solve all these problems
as we next explain.

The change of variables
\begin{equation}
\label{u}
u =-\frac{3(g')^2}{4g^2}+\frac{g''}{2g}.
\end{equation}
transforms the evolution equation (\ref{gpuntodeShiffmant}) into one
of the standard forms of the KdV equation\footnote{In the literature
it is usual to find different other KdV equations, with different
coefficients for $u''', uu'$; all of them are equivalent after a
change of variables.}
\begin{equation}
\label{kdv}
\frac{\partial u}{\partial t} = -u'''-6 u u'.
\end{equation}
The change of variables (\ref{u}) has an easy explanation. The three
terms in the right-hand-side of (\ref{gpuntodeShiffman}) are
rational expressions of derivatives of $g$ with respect to $z$, with
some common homogeneity (order three in derivatives  and degree one
under multiplication $g\mapsto \l g$). One way of transforming
(\ref{gpuntodeShiffman}) into a polynomial expression is by mens of
the change of variables $x=g'/g$, which gives an equation of mKdV
type\footnote{mKdV is an abbreviation for
{\it modified Korteweg-de Vries.}}, namely $\dot{x}=\frac{i}{2}(x'''- \frac{3}{2}x^2x')$. It is a
standard fact that mKdV equations in $x$ can be transformed into KdV
equations in $u$ through the so called {\it Miura transformations,}
$x\mapsto u=ax'+bx^2$ with $a,b$ suitable constants (see for
example~\cite{gewe1} page~273). Equation~(\ref{u}) is just the
composition of $g\mapsto x$ with a Miura transformation. As KdV
theory is more standard than mKdV theory we have opted to deal with
the KdV, although it would have been possible to perform entirely
all what follows directly with the mKdV equation.

A well-known condition on the initial condition $u(z)$ which
guarantees that the Cauchy problem (\ref{kdv}) can be solved
globally producing a holomorphic curve $t\mapsto u_t$ of meromorphic
functions $u_t(z)$ on $\C /\langle i\rangle $, is that $u(z)$ is an
{\it algebro-geometric} potential of the KdV equation. Next we
describe this notion, and we will postpone to
section~\ref{subsec7.3} the property that the holomorphic
integration of (\ref{gpuntodeShiffman}) amounts to solving globally
in $\C /\langle i\rangle $ the Cauchy problem  for equation
(\ref{kdv}) with initial condition $u$ given by~(\ref{u}) (this
reduction is not direct, since the law $g\mapsto u$ in~(\ref{u})
might not be invertible).

\subsection{Algebro-geometric potentials for the KdV equation.}
The KdV equation (\ref{kdv}) is just one of the terms in a sequence of
evolution equations of $u$, called the {\it KdV hierarchy:}
\begin{equation}
\label{kdvn}
\left\{ \frac{\partial u}{\partial t_n} = -\partial_z{\cal P}_{n+1}(u)\right\} _{n\geq 0},
\end{equation}
where $\partial _z=\frac{\partial }{\partial z}$ and ${\cal
P}_{n+1}(u)$ is a differential operator given by a polynomial
expression of $u$ and its derivatives up to order $2n$, defined by
the recurrence law
\begin{eqnarray}
\label{law}
\left\{ \begin{array}{l}
\partial_z {\cal P}_{n+1}(u) = (\partial_{zzz} + 4u\,\partial_z+2u'){\cal P}_{n}(u), \\
\rule{0cm}{.5cm}{\cal P}_{0}(u)=\frac{1}{2}.
\end{array}\right.
\end{eqnarray}
The first operators ${\cal P}_j(u)$ and evolution equations
of the KdV hierarchy are given by
\begin{equation}
\label{eq:KdVhie}
\mbox{}\hspace{-.9cm}
\left.
\begin{array}{l}
{\cal P}_{1}(u)=u\\
\rule{0cm}{.5cm}{\cal P}_{2}(u)=u''+3u^2\qquad \mbox{(KdV)}\\
\rule{0cm}{.5cm}{\cal P}_{3}(u)=u^{(4)}+10 uu''+5(u')^2+10u^3\\
\mbox{}\hspace{.7cm}\vdots
\end{array}
\right|
\begin{array}{l}
\frac{\partial u}{\partial t_0} = -u'\\
\rule{0cm}{.5cm}\frac{\partial u}{\partial t_1} = -u'''- 6uu'\qquad \mbox{(KdV)}\\
\rule{0cm}{.5cm}\frac{\partial u}{\partial t_2} = -u^{(5)}- 10 uu''' -20 u'u''-30 u^2u'\\
\mbox{}\hspace{.7cm}\vdots
\end{array}
\end{equation}
The Cauchy problem for the $n$-th equation of the KdV hierarchy
consists of finding a solution $u(z,t)$ of $\frac{\partial
u}{\partial t_n}= -\partial _z{\cal P}_{n+1}(u)$ with prescribed
initial condition $u(z,0)=u(z)$.

\begin{definition}
\label{def7.1} {\rm Given a meromorphic function $u=u(z)$ defined on
an open set of $\C $, the right-hand-side of (\ref{kdvn}) gives a
sequence of functions of $z$, each of which is a polynomial
expression in $u$ and its derivatives, which we will call {\it
infinitesimal flows} of $u$. With an abuse of notation, we will
denote the $n$-th infinitesimal flow of $u$ by $\frac{\partial
u}{\partial t_n}$ (though $u$ only depends on $z$). The function
$u(z)$ is said to be an {\it algebro-geometric potential of the KdV
equation} (or simply {\it algebro-geometric}) if there exists an
infinitesimal flow $\frac{\partial u}{\partial t_n}$ which is a
complex linear combination of the lower order infinitesimal flows.
%\begin{equation}
%\label{eq:algcj}
%\frac{\partial u}{\partial t_n} =  c_0 \frac{\partial u}{\partial t_0}+ \ldots + c_{n-1}\frac{\partial u}{\partial t_{n-1}},
%\end{equation}
%with $c_0,\ldots,c_{n-1}\in \C $.
}
\end{definition}
Segal and Wilson \cite{SeWi} proved that if $u$ is
algebro-geometric, then it extends to a meromorphic function
$u\colon \C \to \C \cup \{\infty\}$ (see also Gesztesy and Weikard
\cite{gewe1}). We will use later two well-known properties of
algebro-geometric potentials, which can be found in~\cite{gewe1} and
Weikard~\cite{weik1}:
\begin{lemma}
\label{lema7.2}
Let $u(z)$ be an algebro-geometric potential. Then:
\begin{enumerate}
\item If $u$ has a pole at $z=z_0$, then there exists $k\in \Z$ such that around $z_0$,
\[
u(z)=\frac{-k(k+1)}{(z-z_0)^2} + \mbox{\rm holomorphic}(z),
\]
\item All the solutions of the linear Schr\"{o}dinger equation $y''+u\, y =0$
are meromorphic functions $y\colon \C \to \C \cup \{\infty\}$.
\end{enumerate}
\end{lemma}
For our purposes, the key property of algebro-geometric potentials
is that the Cauchy problem for any of the equations in the KdV
hierarchy can be (uniquely) solved if the initial condition is
algebro-geometric. To understand why this is true, suppose
$u=u(z):\C \to \C \cup\{\infty\}$ is algebro-geometric, with
$\frac{\partial u}{\partial t_n} =  c_0 \frac{\partial u}{\partial
t_0}+ \ldots + c_{n-1}\frac{\partial u}{\partial t_{n-1}}$,
$c_0,\ldots ,c_{n-1}\in \C $. Calling $\frac{\partial }{\partial
s}=\frac{\partial }{\partial t_n} - c_0 \frac{\partial }{\partial
t_0}- \ldots -c_{n-1}\frac{\partial }{\partial t_{n-1}}$, then the
next PDE system encodes solving the $k$-th equation in the KdV
hierarchy among functions $u(z,t)$ which are algebro-geometric with
the same coefficients $c_j\in \C$ as $u(z)$ (note that the (A-G) below is
an ODE in~$z$, and (KdV) is a PDE in $z,t$):
\begin{equation}
\label{eq:A-G,KdV}
\left. \begin{array}{l}
{\mbox{\rm (A-G)\hspace{1cm}}\displaystyle \frac{\partial u}{\partial s}=0,}
 \\
\rule{0cm}{.7cm}{\mbox{\rm (KdV)\hspace{0.9cm}}\displaystyle \frac{\partial u}{\partial t} =
-\partial _z{\cal P}_{k+1}(u)}
\end{array}\right\}
\end{equation}

According to the Frobenius Theorem, the integrability condition of
(\ref{eq:A-G,KdV}) is given by the commutativity $\frac{\partial
}{\partial s}\frac{\partial u}{\partial t}=\frac{\partial }{\partial
t} \frac{\partial u}{\partial s}$, which in turn reduces to
$\frac{\partial }{\partial t_j}\frac{\partial u}{\partial t_k}=
\frac{\partial }{\partial t_k} \frac{\partial u}{\partial t_j}$ for
all $j=0,\ldots ,n-1$. This commutativity is a well-known fact in
KdV theory. Therefore, given any $z_0\in \C$ which is not a pole of
$u(z)$, there exists a $\de >0$ and a unique solution $u(z,t)$,
$(z,t)\in \{ |z-z_0|<\delta)\} \times \D(\delta)$, of the
system~(\ref{eq:A-G,KdV}) with initial conditions
\begin{equation}
\label{jet}
\frac{\partial^j u}{\partial z^j }(z_0,0)= u^{(j)}(z_0), \qquad
j=0,\ldots,2n.
\end{equation}
(Note that (A-G) is an ODE of order $2n+1$). Since $u_t=u_t(z)$
is algebro-geometric by (A-G), then $u_t$ extends meromorphically
to the whole plane $\C $. As $u(z)$ satisfies both equations (A-G) and
(\ref{jet}), we have $u(z,0)=u(z)$. Consequently, we have solved the
Cauchy problem for the $k$-th equation of the KdV hierarchy with
algebro-geometric initial condition $u(z)$.

\subsection{Proof of Theorem~\ref{thmintegrShiffman} provided that $u$ is
algebro-geometric.}
\label{subsec7.3}

Coming back to our setting where $u(z)$ is given by (\ref{u}) for a
given $g\in {\cal M}_{\mbox{\rm \footnotesize imm}}$, {\it suppose
that $u$ is algebro-geometric.} By the arguments in the last
section, we can solve the Cauchy problem for the KdV equation (which
is the $k$-th equation in the KdV hierarchy for $k=1$), obtaining a
family $t\in \D (\de )\mapsto u_t(z)=u(z,t)$ of algebro-geometric
meromorphic functions on $\C $, with $u(z,0)=u(z)$. It is then
interesting to know how periodicity of the initial condition
propagates to $u_t(z)$ (note that $u(z)$ is defined on $\C /\langle
i\rangle $ since $g$ is). Clearly, the uniqueness of solution of
(\ref{eq:A-G,KdV}) with initial conditions (\ref{jet}) implies that
if $u(z)$ is invariant by the translation of a vector $\omega \in \C
$, then $u_t(z)$ has the same invariance for all $t$. In particular
for each $t\in \D(\delta)$, $u_t(z)$ descends to the cylinder $\C /
\langle i\rangle $ as a meromorphic function. It is now time to
produce $g_t(z)$ from $u_t(z)$ via equation~(\ref{u}), which is our
next goal.

To do this, we first choose a meromorphic function $y\colon \C \to
\C \cup\{\infty\}$ such that $g=1/y^2$ (the existence of $y$ is
guaranteed because $g$ has double zeroes and double poles without
residues since $g\in {\cal M}_{\mbox{\footnotesize imm}}$). Note
that $y$ is either periodic, $y(z+i)=y(z)$, or anti-periodic
$y(z+i)=-y(z)$ (in this last case $y$ does not descend to $\C
/\langle i\rangle $ and this is what happens with the Riemann
minimal examples, where the Gauss map $g$ restricts to each compact
horizontal section with degree one). A direct computation using
(\ref{u}) gives that $y''+uy=0$. Now consider the PDE system with
unknown $y(z,t)$:
\begin{equation}
\label{PDEy}
\left. \begin{array}{l}
{\mbox{\rm (S)\hspace{1cm}}\displaystyle y''+uy=0,}
 \\
{\mbox{\rm (J)\hspace{1cm}}\rule{0cm}{.5cm}\frac{\partial y}{\partial t}
= {\cal P}_1(u)' y - 2 {\cal P}_1(u)y'}
\end{array}\right\}
\end{equation}
where the function $u$ in (\ref{PDEy}) is the above solution
$u(z,t)$ of the Cauchy problem (\ref{eq:A-G,KdV})-(\ref{jet}). Note
that (S) is a Schr\"{o}dinger type ODE in the variable $z$; on the
contrary, (J) is a PDE in $z,t$. The reason why we consider the
system (\ref{PDEy}) is that its integrability condition is precisely
that $u(z,t)$ solves the KdV equation\footnote{This fact can be
generalized to the $k$-th equation of the KdV hierarchy only by
changing $t$ by $t_k$ and ${\cal P}_1$ by ${\cal P}_k$ in equation
(J).}, as proved by Joshi \cite{joshi1}. Therefore, the Frobenius
theorem implies that (\ref{PDEy}) admits a unique solution
$y=y(z,t)$ with initial condition $y(z,0)=y(z)$. Since $z\mapsto
u(z,t)$ is algebro-geometric for every $t$, part~{\it 2} of
Lemma~\ref{lema7.2} together with equation (S) imply that $y(z,t)$
is defined on $\C \times \D(\varepsilon)$ (for some $\ve >0$) and is
meromorphic in $z$. The uniqueness of solution of an initial value
problem together with the fact that $y(z+i)=\pm y(z)$, give that
$y(z+i,t)=\pm y(z,t)$, with the same choice of signs as for $y(z)$.
Finally, defining
\begin{equation}
\label{eq:gt}
g_t(z)=g(z,t)=y^{-2}(z,t)
\end{equation}
then
\[
\frac{\partial g_t}{\partial t} =
\frac{\partial }{\partial t}\left( \frac{1}{y_t^2}\right) =
-\frac{2}{y_t^3}\frac{\partial y_t}{\partial t}
\stackrel{\rm (J)}{=}
-2\frac{{\cal P}_1(u_t)' y_t - 2 {\cal P}_1(u_t)y_t'}{y_t^3}=
-2\,\partial_z\left(\frac{{\cal P}_1(u_t)}{y_t^2}\right)
=-2\,\partial_z\left( g_t{\cal P}_1(u_t)\right)
\]
\[
\stackrel{(\ref{eq:KdVhie})}{=}-2\,\partial_z\left( g_tu_t\right)
\stackrel{(\star )}{=}-2\,\partial_z\left[ g_t\left( -\frac{3(g_t')^2}{4g_t^2}+\frac{g_t''}{2g_t}
\right) \right]=
-g_t'''+3\frac{g_t'g_t''}{g_t}-
\frac{3}{2}\frac{(g_t')^3}{g_t^2},
\]
where in $(\star )$ we have used (S) and (\ref{eq:gt}) (we cannot
substitute directly~(\ref{u}) since {\it a priori} is only valid for
$u(z,0)=u(z)$). The equality in the last two lines tells us that, up
to a multiplicative constant, $t\mapsto g_t$ satisfies the evolution
equation in item {\it iii)} of Theorem~\ref{thmintegrShiffman}.
Therefore, in order to finish the proof of
Theorem~\ref{thmintegrShiffman} we need to demonstrate items {\it
i), ii)} of that theorem; these are technical issues that we next
sketch.

Equation (\ref{u}) implies that the poles of $u(z)$ coincide with
the zeros and poles of $g(z)$, which are double. A direct
computation gives that the Laurent expansion of $u(z)$ around each
pole $z_0$ is of the form
\begin{equation}
\label{eq:thm6.4*}
u(z)=\frac{-2}{(z-z_0)^2}+\mbox{\rm holomorphic}(z).
\end{equation}
The next step consists of proving that every pole $z_0$ of $u(z)$
propagates holomorphically in $t$ to a curve of poles $z_0(t)$ of
$u_t(z)$ with a similar Laurent expansion as in~(\ref{eq:thm6.4*}).
The argument is purely local: On one hand, as $u_t(z)$ is
algebro-geometric for all $t$, item~{\it 1} of Lemma~\ref{lema7.2}
implies that $u_t$ admits a Laurent expansion of the type
\[
u_t(z)=\frac{-k_j(k_j+1)}{(z-a_j)^2}+\mbox{\rm holomorphic}(z,t)
\]
in a neighborhood of each of the poles $a_1,\ldots ,a_m$ of $u_t$ in
a fixed closed disk $D$ centered at $z_0$, such that $u(z)$ does not
vanish in $D-\{ z_0\} $ and $u_t$ has no zeros in $\partial D$ (both
the number of poles of $u_t$ and the poles themselves may depend on
$t$). Now a continuity argument with respect to $t$ together with
the fact that $k_j$ is integer-valued, give that for $|t|$ small,
$u_t(z)$ has a unique pole in $D$ and $k_1=1$. Once we know that
$u_t$ has a {\it unique} pole in $D$, the holomorphic dependence of
this pole with respect to $t$ is a standard argument.

In order to obtain the desired holomorphicity of the curves of zeros
and poles of $g_t(z)$ with respect to $t$ (which is item {\it i)} of
Theorem~\ref{thmintegrShiffman}), we need to know that these zeros
and poles of $g_t$ (or of $y_t$) coincide with the poles of $u_t$:
we already know that the series expansion of $u_t$ around each of
its poles $z_0(t)$ is of the form (\ref{eq:thm6.4*}) with $z_0(t)$
instead of $z_0$ (of course, the holomorphic term in the
right-hand-side also depends on $t$). Using that $y_t''+u_ty_t=0$
(equation (S)), it is straightforward to expand $y_t$ around
$z_0(t)$ thereby proving that outside of the poles of $u_t$, the
function $y_t$ is holomorphic and its zeros are simple, while at
each pole of $u_t$, the function $y_t$ has either a simple pole or a
double zero. Another continuity argument shows that the possibility
of $y_t$ having a double zero at a pole of $u_t$ cannot occur (for
$t=0$ it does not occur, since $y(z)=1/\sqrt{g(z)}$ has only single
zeros and poles). Therefore the zeros and poles of $g_t$ are double,
and coincide with the poles of $u_t$. In particular, the holomorphic
dependence with respect to $t$ in item~{\it i)} of
Theorem~\ref{thmintegrShiffman} holds.

The quasiperiodicity of $u_t(z)$ is guaranteed by that of $u(z)$
together with the uniqueness of the solution of (\ref{eq:A-G,KdV})
(with $k=1$). In turn, this quasiperiodicity of $u_t(z)$ implies the
same property for $y_t(z)$ and thus for $g_t(z)$. To finish the
proof of items {\it i), ii), iii)} of
Theorem~\ref{thmintegrShiffman}, it only remains to show that
$g_t\in {\cal M}_{\mbox{\rm \footnotesize imm}}$. This follows from
the fact that the map $t\mapsto \mbox{Per}(g_t)$ is constant, see
equation~(\ref{eq:periodconstant}) and Remark~\ref{remark6.12}.

Finally, the last sentence in the statement of Theorem~\ref{thmintegrShiffman}
follows from the maximum principle for minimal surfaces. This finishes the sketch
of proof of Theorem~\ref{thmintegrShiffman}.

\subsection{Why $u$ is algebro-geometric if $g\in {\cal M}_{\mbox{\rm \footnotesize imm}}$.}
\label{subsec7.4} In the last section we proved
Theorem~\ref{thmintegrShiffman} under the additional assumption that
the function $u=u(z)$ given by~(\ref{u}) is an algebro-geometric
potential of the KdV equation. We will devote this section to
explaining why this assumption for $u$ holds for all $g\in {\cal
M}_{\mbox{\rm \footnotesize imm}}$. Briefly, the desired property
for $u$ follows from two facts: firstly, that each infinitesimal
flow $\frac{\partial u}{\partial t_n}$ of $u$ given in
Definition~\ref{def7.1} produces a complex valued, bounded Jacobi
function $v_n$ on $M$ (for instance, $\frac{\partial u}{\partial
t_1}$ produces $S_M+iS_M^*$) which extends smoothly across the zeros
and poles of $g$ to a function in the bounded kernel of an operator
of the type $\Delta +V$ on $\esf^1\times \R $, where $\Delta $ is
the laplacian in the product metric and $V\colon \esf^1\times \R \to
\R $ is a bounded potential. And secondly, that the bounded kernel
of such a Schr\"{o}dinger operator is finite-dimensional. We will
now develop the details of this sketch.

From now on, we consider a function $g\in {\cal M}_{\mbox{\rm
\footnotesize imm}}$ and let $u$ the meromorphic function on $\C
/\langle i\rangle $ given by (\ref{u}). Similarly as in
Definition~\ref{def7.1}, we define the {\it infinitesimal flows for
$g$} as the sequence of meromorphic functions on $\C /\langle
i\rangle $ given by
\begin{equation}
\label{sh}
\left\{ \frac{\partial g}{\partial t_n} =
-2\,\partial_z (g{\cal P}_n(u))\right\} _{n\geq 0}.
\end{equation}
Hence, each $\frac{\partial g}{\partial t_n}$ is a rational
expression in $g$ and its derivatives up to some order. By
substituting (\ref{law}) and (\ref{u}) in (\ref{sh}), one can
compute explicitly the infinitesimal flows for $g$; for instance,
$\frac{\partial g}{\partial t_0}=-g'$ is the infinitesimal
deformation of $g$ in ${\cal W}$ given by translations in the
parameter domain (see Remark~\ref{remarkABC}), and $\frac{\partial
g}{\partial t_1}=-g'''+3\frac{g'g''}{g}-
\frac{3}{2}\frac{(g')^3}{g^2}$ is, up to a multiplicative constant,
the infinitesimal deformation $\dot{g}_S$ given in
equation~(\ref{gpuntodeShiffman}), which corresponds to the (complex
valued) Shiffman function. In particular, both $\frac{\partial
g}{\partial t_0}, \frac{\partial g}{\partial t_1}$ belong to
$T_g{\cal W}$.

For each $n\geq 0$, the infinitesimal flow $\frac{\partial
g}{\partial t_n}$ for $g$ satisfies the following two key
properties:
\begin{enumerate}
\item {\it There exists a meromorphic function $h_n$ on $\C /\langle i\rangle $
which is a rational expression of $g$ and its derivatives up to some
order (depending on $n$), such that $\frac{\partial g}{\partial
t_n}=\partial_z\left( \frac{g^3h_n'}{2g'}\right) $.}
\newline
The proof of this property reduces to integrating the equality
$-2\partial _z(g{\cal P}_n(u))=\partial _z\left(
\frac{g^3h_n'}{2g'}\right) $, which is an ODE in $z$ with unknown
$h_n$ (the recurrence law~(\ref{law}) for the operators ${\cal P}_n$
is useful here).
\item {\it The principal divisor $D$ of the meromorphic
function $\frac{\partial g}{\partial t_n}$ satisfies $D\geq \prod
_jp_jq_j^{-3}$, where the principal divisor of $g$ is $(g)=\prod
_jp_j^2q_j^{-2}$. Therefore, $\frac{\partial g}{\partial t_n} \in
T_g{\cal W}$.}
\newline
This can be deduced from the local expansions of $g$ (which has
double zeros and double poles without residue since $g\in {\cal
M}_{\mbox{\rm \footnotesize imm}}$), of $u$ (given by
equation~(\ref{eq:thm6.4*})) and of ${\cal P}_n(u)$ (proven by
induction on $n$).
\end{enumerate}

The properties 1,2 above together with Proposition~\ref{propos4.8}
produce, for each $n\geq 0$, a complex valued Jacobi function
$f(h_n)=\frac{g^2h'_n}{g'}+\frac{2gh_n}{1+|g|^2}\in {\cal J}_{\C
}(g)$ which is bounded on $\C /\langle i\rangle $. To continue our
argument, we need the following finiteness result, whose proof we
postpone till the end of this section.

\begin{theorem}
\label{bounded} Let $M\subset \R^3$ be a quasiperiodic, immersed
minimal surface of Riemann type. Then, the linear space of bounded
Jacobi functions on $M$ is finite dimensional.
\end{theorem}
Assuming that Theorem~\ref{bounded} holds, we finish the proof that
$u$ is algebro-geometric. We had constructed a sequence of bounded
Jacobi functions $f(h_n)\in {\cal J}_{\C }(g)$. By
Theorem~\ref{bounded}, only finitely many of them can be linearly
independent, hence there exists $n\in \N $ such that $f(h_n)\in
\mbox{Span}\{ f(h_0),\ldots ,f(h_{n-1})\} $. Since the linear map
$h\mapsto f(h)$ given by equation~(\ref{eq:Jacobideh}) is injective,
we conclude that $h_n\in \mbox{Span}\{ h_0,\ldots ,h_{n-1}\} $. As
the map $h\mapsto \dot{g}(h)$ given by (\ref{eq:gpuntodeh}) is
linear, then $\frac{\partial g}{\partial t_n}\in \mbox{Span}\{
\frac{\partial g}{\partial t_0},\ldots ,\frac{\partial g}{\partial
t_{n-1}}\} $. Finally, the equation
\[
\frac{\partial u}{\partial t_n} =
\frac{\partial }{\partial t_n}\left(-\frac{3(g')^2}{4g^2}+\frac{g''}{2g} \right) ,
\]
implies that $\frac{\partial u}{\partial t_n}\in
\mbox{Span}\{ \frac{\partial u}{\partial t_0},\ldots ,\frac{\partial u}{\partial t_{n-1}}\} $,
which proves that $u$ is algebro-geometric.

We finish this section with a sketch of the proof of
Theorem~\ref{bounded} stated above. Suppose that $M\subset \R^3$ is
a quasiperiodic, immersed minimal surface of Riemann type. Then, $M$
is conformally equivalent to $(\C /\langle i\rangle )-g^{-1}(\{
0,\infty \} )$ where $g\in {\cal M}_{\mbox{\footnotesize imm}}$ is
the Gauss map of $M$. Take global coordinates $(\theta ,t)$ on
$\esf^1\times \R $ and consider the product metric $d\theta ^2\times
dt^2$, which is conformal to the metric $ds^2$ on $M$ induced by the
usual inner product of $\R^3$: $ds^2=\l ^2(d\theta ^2+dt^2)$. This
conformality allows us to relate the Jacobi operator $L=\Delta -2K$
of $M$ (here $K$ is the Gaussian curvature of $M$)
to a Schr\"{o}dinger operator $L_M=(\Delta_{\esf^1}+\partial^2_t)+V$ on
$\esf^1\times \R$ by means of the formula $L=\l ^{-2}L_M$, where the
potential $V$ is equal to the square of the norm of the differential
of the Gauss map of $M$ (with respect to $d\theta ^2\times dt^2$).
The quasiperiodicity of $M$ implies that $V_M$ is globally bounded
on $\esf^1\times \R $.

By elliptic regularity, any bounded Jacobi function $v$ on $M$
extends smoothly through the zeros and poles of $g$ to a function
$\widehat{v}$ in the kernel of $L_M$, such that $\widehat{v}$ is
bounded at both ends of $\esf^1\times \R $. Therefore, the space of
bounded Jacobi functions on $M$ (i.e. bounded functions in the
kernel of $L$) identifies naturally with the bounded kernel of
$L_M$, and thus Theorem~\ref{bounded} follows from the following
standard technical result, whose proof we omit here (see
Assertion~5.3 in~\cite{mpr6} for a proof due to Frank Pacard, based
on the paper by Lockhart and McOwen~\cite{loMcOw1}).

\begin{assertion}
\label{assloMcOw2} Let $(\Sigma ,h)$ be a compact Riemannian
manifold and $V\in L^{\infty } (\Sigma \times \R )$. Assume that
there exists $j_0\in \N$ such that\footnote{In the case where
$\Sigma $ is the standard $\esf^1$, then $\l _j=j^2$ and $\l
_{j+1}-\l _j=2j+1$, so the hypothesis (\ref{eq:Pacard3}) is
fulfilled. }
\begin{equation}
\label{eq:Pacard3}
4\| V\| _{L^{\infty }(\Sigma \times \R )}\leq \l _{j_0+1}-\l _{j_0},
\end{equation}
where $\l _0=0<\l _1<\l _2<\ldots $ is the spectrum of $-\Delta _h$
on $\Sigma $. Then, the bounded kernel of $\Delta _h+\partial ^2_t+V$
on $\Sigma \times \R $ is finite dimensional.
\end{assertion}

\section{The asymptotics of the ends of finite genus surfaces.} \label{secas}

In this section we briefly describe the asymptotic behavior of the
ends of a properly embedded minimal surface $M$ in $\rth$ with
finite genus, possibly with compact boundary.
Collin's solution of the Nitsche Conjecture
\cite{col1} implies that if $M$ has more than one end, then each
annular end is asymptotic to the end of a plane or catenoid.

Next assume that $M$ has just one end. On the last page of their
paper~\cite{mr8}, Meeks and Rosenberg claimed that their proof
of the uniqueness of the helicoid could be modified  to prove
the following statement:
\begin{theorem}
\label{thm1endnew}
Any non-planar, properly embedded minimal surface $M$
in $\rth$ with one end, finite topology and infinite total
curvature (without boundary) satisfies the following properties:
\ben
\item $M$ is conformally a compact Riemann surface
$\overline{M}$ punctured in a single point.
\item After a rotation in $\R^3$, the Weierstrass pair $(g,dh)$ of $M$ satisfies
that both differentials $\frac{dg}{g}$, $dh$ extend meromorphically
to $\overline{M}$.
\item $M$ is asymptotic to a helicoid.
\een
\end{theorem}
The technical nature of the proof by Meeks and Rosenberg of the
uniqueness of the helicoid and the absence of a detailed proof of their
claimed generalization above to the case of (non-zero) finite genus,
motivated subsequent investigations. Bernstein and Breiner~\cite{bb2}
have recently given a proof of Theorem~\ref{thm1endnew}, based on
arguments from the original paper by Meeks and Rosenberg together
with a more careful analysis of the multigraph structure of the
end of a surface $M$ under the hypotheses of Theorem~\ref{thm1endnew},
using Colding-Minicozzi theory. With these two ingredients,
Bernstein and Breiner were able to prove that $M$ has finite
type (see Definition~\ref{deffinitetype} below for the notion of
a minimal surface of finite type). Once $M$ is proven to have finite type,
one can use former results by Hauswirth, P\'erez and Romon~\cite{hkp1}
on the geometry of complete embedded, minimal ends of finite type
to complete the proof of Theorem~\ref{thm1endnew}.

\begin{definition}[Finite Type]
\label{deffinitetype}
{\rm A minimal immersion $X\colon M \to \rth$ is said to have
{\em finite type} if it satisfies the following two properties.
\ben
\item The underlying Riemann surface to $M$ is conformally diffeomorphic
to a compact Riemann surface $\overline{M}$ with (possibly empty)
compact boundary, punctured in a finite non-empty set ${\cal E}\subset \Int(M)$.
\item Given an end $e\in {\cal E}$ of $M$, there exists a rotation of the surface in space
such that if $(g,dh)$ is the Weierstrass data of $M$ after this rotation,
then the meromorphic one-forms
$\frac{dg}{g}$ and $dh$ extend across the puncture $e$ to meromorphic one-forms on
a neighborhood of $e$ in $\overline{M}$.
\een
}
\end{definition}

In their survey~\cite{mpe2}, the authors of these notes
outlined the proof by Meeks and Rosenberg
of the uniqueness of the helicoid and at the end of this outline they
mentioned how some difficult parts of the proof could be simplified,
as for instance the facts that the correct conformal structure is
$\C $ or that the height differential can be assumed to be $dz$
with $z$ being the natural coordinate on $\C $, see also
Footnote~\ref{foootn} above. In
the recent paper~\cite{mpe3}, Meeks and P\'erez
not only give the aforementioned simplification
of the proof of the uniqueness of the helicoid, but also tackle the more general problem of
describing the asymptotic behavior, conformal structure and analytic
representation of an annular end of any {\it complete,} injectively immersed minimal surface $M$ in $\rth$ with {\it compact boundary} and finite
topology. Before proceeding to explain the
results in~\cite{mpe3}, it is necessary to make
a remark about properness versus completeness.
Although not explicitly stated in the
paper~\cite{cm35} by Colding and Minicozzi, the results contained there imply that such an $M$ is {\it properly embedded} in $\rth$ (we also remark that in Meeks, P\'erez and Ros~\cite{mpr9} the following
more general result is proven: If $M$ is a complete, connected,
 injectively immersed minimal surface of finite genus, compact boundary and
a countable number of ends in $\R^3$, then $M$ is proper).

We next state the first main result in~\cite{mpe3}.
\begin{theorem}
\label{th1.1}
Let $E\subset \rth$ be a complete, embedded minimal annulus with
infinite total curvature and compact boundary.
Then, the following properties hold:
\begin{enumerate}
\item $E$ is properly embedded in $\R^3$.
\item $E$ is conformally diffeomorphic to
$D(\infty ,R)=\{z\in \C\mid  R\leq |z| \}$.
\item After a suitable homothety and
rigid motion and possibly replacing $M$ by a subend, then:
\ben
\item The height differential $dh=dx_3+idx_3^*$
extends meromorphically across
 infinity with a double pole.
\item The stereographically projected Gauss map $g\colon D(\infty ,R)\to \C\cup\{\infty\}$ of $M$ can be expressed as $g(z)=e^{iz+f(z)}$ for some
holomorphic function $f$ in $D(\infty ,R)$ with $f(\infty )=0$.
\item $E$ is asymptotic to the end of a helicoid if and only if it has zero flux.
\een
\end{enumerate}
\end{theorem}
Note that Theorem~\ref{thm1endnew} above is
a direct consequence of Theorem~\ref{th1.1}
(the zero flux condition in item~{\it 3(c)}
of Theorem~\ref{th1.1} follows from Stokes' theorem).

The main ideas in the proof of
Theorem~\ref{th1.1} are the following ones. First
one proves that the sequence of surfaces
$\{ \l _nE\} _n$ has locally positive injectivity radius in $\R^3-\{ \vec{0}\} $, for every
sequence $\{ \l_n\} _n \subset \R^+$ with
$\l _n\searrow 0$ as $n\to \infty $ (this means that
for every $q \in \R^3-\{ \vec{0}\} $, there exists
$\ve _q>0$ and
$n_q\in \N $ such that for $n>n_q$, the injectivity radius function
of $\l _nE$ restricted to $\{ x\in \R^3\ | \ |x-q|<\ve _q\} \cap (\l _nE)$
is a sequence of functions which is uniformly bounded away from zero; this property is obtained after a
blow-up argument on the
scale of topology similar to the one explained
during the proof of Proposition~\ref{propos2.11}).
The second step in the proof is to apply a result about singular
minimal laminations (namely item 7 of Theorem~1.5 of Meeks, P\'erez and Ros~\cite{mpr11}) to conclude that after extracting a subsequence, the surfaces
$\l _nE$ converge as $n\to \infty $ to a foliation ${\cal F}$ of $\rth$ by parallel planes and the convergence is $C^1$ away from one
straight line orthogonal to the planes in ${\cal F}$
(the presence of boundary prevents us to use Colding-Minicozzi theory in our setting to get this).
Under shrinkings, the boundary of $\l _nE$ collapses into the origin $\vec{0}$, which in turn implies that
the singular set of convergence of the $\l _nE$ to ${\cal F}$ must be a line passing through
$\vec{0}$. Furthermore, this limit foliation
${\cal F}$ is independent of the sequence of positive numbers $\l_n$, hence it can be assumed from now on
that the planes in this foliation are horizontal.

It follows from the previous paragraph that
there is a solid vertical hyperboloid ${\cal H}$ with
axis being the $x_3$-axis, such that  $E-{\cal H}$ consists of two multigraphs over
their projections to the $(x_1,x_2)$-plane $P$. By work of Colding-Minicozzi, any embedded minimal multigraph with a large number of sheets
contains a submultigraph which can be approximated by the multigraph of a helicoid with an additional logarithmic term (see Corollary~14.3 in~\cite{cm34}).
A deeper study of this multigraph structure using
Colding-Minicozzi theory leads to the property that
that each of the two multigraphs $G_1$, $G_2$ in
$E-{\cal H}$ contains  infinite submultigraphs $G_1'$, $G_2'$, respectively
given by functions $u^1(\rho,\t)$, $u^2(\rho,\t)$, such that
$\frac{\partial u^i}{\partial \t}(\rho, \t)>0$ (resp. $<0$)
for $i=1,2$ (this observation was also made by Bernstein and Breiner using the same arguments,
see Proposition~3.3 in \cite{bb2}). The positive
slope property of the curves $\t \mapsto u^i(\rho ,\t )$ is then used to prove that after passing to a
subend (denoted in the same way) of $E$, each horizontal plane $\{ x_3=t\} $ intersects $E$ transversely in either a proper curve at height $t$ or in two proper arcs, each with one extremum on the boundary of $E$. This intersection property, together
with Corollary~1.2 in Meeks and P\'erez~\cite{mpe5}, imply that the conformal structure of $E$ is a punctured disk
$D(\infty ,R)=\{ z\in \C \mid |z|\geq R\} $, $R>0$,
that the height differential $dh$ if $E$ extends meromorphically across $z=\infty $ with a double pole
and (again  after passing to a subend) that the Gauss map
$g$ of $E$ can be written as $g(z)=z^ke^{H(z)}$ for
some $k\in \Z $ and some holomorphic function $H$
in $D(\infty ,R)$.

\begin{remark}
\label{remnew}
At this point in the sketch of proof of Theorem~\ref{th1.1}, one can give a short proof of
the uniqueness of the helicoid among simply-connected, embedded, complete, non-flat minimal surfaces in
$\R^3$: It is clear from the last paragraph that in
this special setting for our minimal surface $M$ in
question,  $M$ is conformally $\C $ and $dh=\l \, dz$ for some
$\l \in \C-\{ 0\} $, or after a change of coordinates,
$dh=dz$. As $dh$ has no zeros,
then the Gauss map $g$ misses $0,\infty $ on the whole surface $M$. Since $M$ is simply-connected, then $g$ lifts through the natural exponential
map $e^w\colon \C ^*\to \C $ and thus, $g(z)=e^{H(z)}$ for some entire function $H$. From this point
one finishes the uniqueness of the helicoid as
in the original proof by Meeks and Rosenberg sketched
just after the statement of Theorem~\ref{ttmr},
 proving that $H$ is a linear function of $z$.
\end{remark}
We will finish our sketch of proof of Theorem~\ref{th1.1} by indicating how to arrive to the
expression $g(z)=e^{iz+f(z)}$ as in item {\it (b)}
of that statement (item {\it (c)} is a consequence
of Theorem~\ref{thm1.3} below). To do this we first prove that $g(z)=e^{H(z)}$, or in other words, $k=0$
in our previous expression $g(z)=z^k e^{H(z)}$. The integer $k$ is the
winding number of $g|_{\partial E}$, so it suffices to prove that this
winding number vanishes. Since the winding number of $g|_{\partial E}$ is
an invariant of the homotopy class of $\partial E$ in $E$, we replace
$E$ by a subend $E'$ where the winding number of $g|_{\partial E'}$ is easier
to calculate and the calculation has a geometric nature. More precisely, we can
choose $E'$ so that $\partial E'$ is a simple closed curve close to the boundary
of a cylinder $C(r,h)=\{(x_1,x_2,x_3)\in \R^3\mid x_1^2+x_2^2\leq r^2, |x_3|\leq h\}$.
Furthermore, $\partial E'$ consists of an arc $\a_T$, which is contained in a plane
parallel to the $(x_1,x_2)$-plane and it is close in the $C^1$-sense to a line
segment parallel to the $x_2$-axis contained in the top disk of $C(r,h)$, a
similar arc $\a_B$ near the bottom disk of $C(r,h)$ and two spiraling arcs
$S_1, S_2$ on the boundary cylinder of $C(r,h)$, where each of these spirals
winds exactly $n$ times around the $x_3$-axis for some large $n\in \N$. Note
that the end points of $\a_t$ and $\a_B$ lie on the cylindrical sides of
$C(r,h)$ and lie vertically over the points $(0,\pm r,0)$.

By construction, the argument of the Gauss map of $E'\mod 2\pi$  restricted to
$\a_T\cup \a_B$ lies in $(-\frac{\pi }{2},\frac{\pi }{2})$; to obtain this property  for this portion
of $\partial E'$,  one uses the fact that under rescalings and translations
the parts of $E'$ in the halfspaces $\{x_3\geq 0\}$, $\{x_3\leq 0\}$ produce
vertical right handed helicoids and also we apply the fact
$\frac{\partial u}{\partial \theta}>0$ holds on certain minimal multigraphs
$u(\rho,\theta)$  with many sheets and defined in polar coordinates on annular
domains in the plane. Finally, using the property
$\frac{\partial u_i}{\partial \theta}>0$
for the two multigraphs $S_1, \, S_2$ corresponding
to multigraphing functions $u_1,u_2$ over
the universal cover of
$\{(x_1,x_2,0) \mid x_1^2 + x_2^2 \geq r^2\}$, we can relate the change of
the arguments of $g$ along $S_1$ and $S_2$ to the number $n$ of their
windings around the $x_3$-axis and these changes of arguments essentially
cancel each other out to make the total change in the argument as $g$
transverses $\partial E'$ to be less than $2\pi$ in absolute value; to
see this holds it is helpful to
notice that if under a parametrization of $\partial E'$,  the winding
number of $S_1$ is $n$ around the $x_3$-axis, then the winding number of $S_2$ is $-n$.
Hence, the winding number of $g|_{\partial E'}$ is $k=0$, which proves that $g(z)=e^{H(z)}$.

Once we have shown $g(z)=e^{H(z)}$, we next consider two cases, depending on whether or
not $H(z)$ has an essential singularity. Following the arguments
of Meeks and Rosenberg in~\cite{mr8}, we find that the case where $H(z)$
has an essential singularity is impossible, while the case that $H(z)$
extends meromorphically across infinity can be reduced, after a change of variables,
to $H(z)=iz+f(z)$ where $f(\infty)=0$, which completes our indication as to why
$g(z)=e^{iz+f(z)}$ with $f(\infty)=0$.

\begin{figure}
\begin{center}
\includegraphics[height=8.3cm]{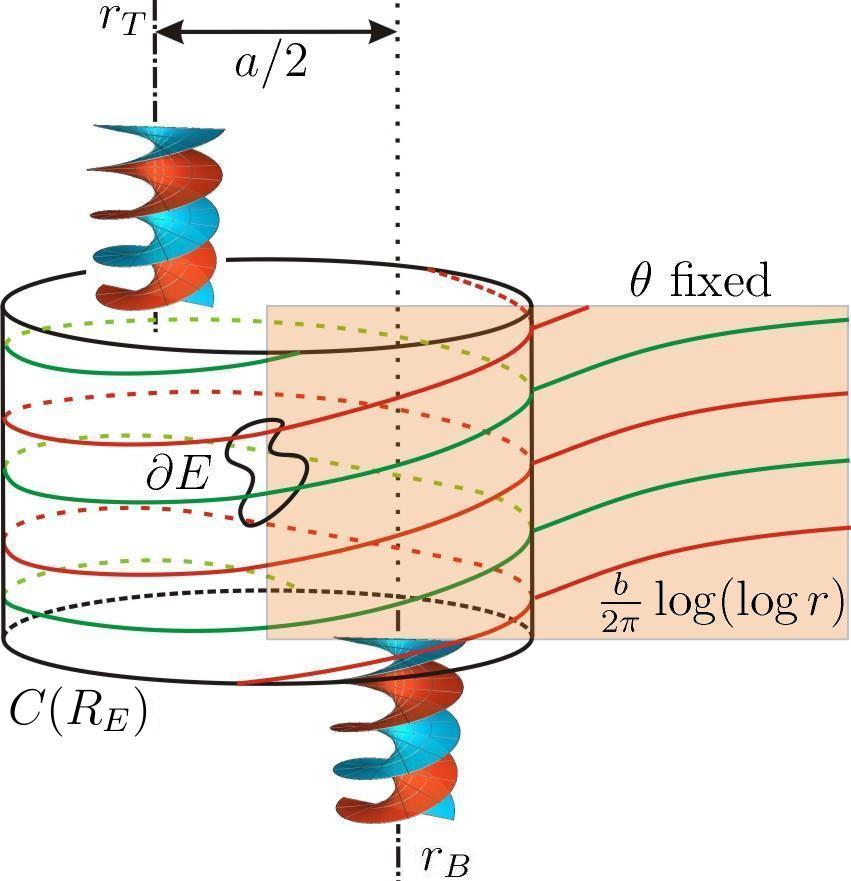}
\caption{The embedded annulus $E$ with flux vector $(a,0,b)$
has the following description.
Outside the cylinder $C(R_E)$, $E$ consists of two horizontal
multigraphs with asymptotic spacing $\pi $ between them.
The translated surfaces $E+(0,0,-2\pi n-\frac{b}{2\pi }\log n)$
(resp. $E+(0,0,2\pi n-\frac{b}{2\pi }\log n)$) converge as $n\to
\infty $ to a vertical helicoid $H_T$ (resp. $H_B$) such that
$H_B=H_T+(0,a/2,0)$ (in the picture, $r_T,r_B$ refer to the axes of
these helicoids).
The intersection of a vertical halfplane containing the $x_3$-axis
with $E-C(R_E)$ consists of an infinite number of curves, each of
which is a graph of a function $u(r)$ that satisfies the property
$\frac{u(r)}{\log (\log r)}$ converges to $\frac{b}{2\pi }$ as the
radial distance $r$ to the $x_3$-axis tends to $\infty $.
} \label{Eab}
\end{center}
\end{figure}

In contrast to the unique asymptotic behavior
of a complete, embedded minimal surface with infinite
total curvature, one end and {\it no boundary,} the second main result of~\cite{mpe3} shows that if we
let the annular end $E$  have non-zero flux (in particular, it is not part of a minimal surface
without boundary), then we have many other asymptotic
models: they are essentially given by the flux vector along the boundary of $E$,
which, after a rotation around the $x_3$-axis, is
$(a,0,b)\in \R^3$.  Therefore the really different asymptotic structures form a
2-parameter family of {\it canonical ends} $\{
E_{a,b}\mid a,b\geq 0\} $. Here, the word {\it canonical}
only refers to item~{\it 3} in the next statement.
The image in Figure~\ref{Eab} describes how the flux vector $(a,0,b)$ of
$E$ influences its geometry. 
%\vspace{-.5cm}

\begin{theorem}[Asymptotics of embedded minimal annular ends]
\label{thm1.3} Given $a,b\geq 0$, there exist a positive number
$R=R(a,b)$ and a properly embedded minimal annulus $E_{a,b}\subset
\R^3$ with compact boundary and flux vector $(a,0,b)$ along its
boundary, such that the following statements hold.
 \ben
\item $E_{a,b}-C(R)$ consists of two disjoint 
multigraphs\footnote{See Footnote~\ref{Ngraph} for
the notion of multigraph.} $\Sigma _1,
\Sigma _2$ over $D(\infty ,R)$ of smooth functions $u_1,u_2\colon
\widetilde{D}(\infty ,R)\to \R $ such that their gradients satisfy
$\nabla u_i(r,\theta )\to 0$ as $r\to \infty $ and the separation
function $w(r,\theta )=u_1(r,\theta )-u_2(r,\theta )$ between both
multigraphs converges to $\pi $ as $r+|\theta |\to \infty $.
Furthermore for $\t $ fixed and $i=1,2$,
\begin{equation}
\label{eq:graphr}
 \lim_{r\to \infty} \frac{u_i(r,\theta)}{\log
(\log(r))} =\lim_{r\to \infty} \frac{u_i(r,\theta)}{\log
(\log(r))}=\frac{b}{2\pi }.
\end{equation}
\item The translated surfaces $E_{a,b}+(0,0,-2\pi n-\frac{b}{2\pi
}\log n)$ (resp. $E_{a,b}+(0,0,2\pi n-\frac{b}{2\pi }\log n)$)
converge as $n\to \infty $ to a vertical helicoid $H_T$ (resp.
$H_B$) such that $H_B=H_T+(0,a/2,0)$. Note that this property
together with item 1 imply that for different values of $a,b$, the
related surfaces $E_{a,b}$ are not asymptotic after a rigid motion
and homothety.
\item Every complete, embedded minimal annulus in $\R^3$ with compact boundary
and infinite total curvature is asymptotic
(up to a rigid motion and homothety) to exactly one of the surfaces $E_{a,b}$.
\een
\end{theorem}

%\begin{theorem}[Asymptotics of embedded minimal annular ends]
%\label{thm1.3}
%Given $a,b\geq 0$, there exist a positive number $R=R(a,b)$
%and a properly embedded minimal annulus $E_{a,b}\subset \R^3$ with compact boundary
%such that the following statements hold.
%\ben
%\item $E_{a,b}-\{ (x_1,x_2,x_3)\in \R^3\mid
%x_1^2+x_2^2<R^2\} $ consists of two disjoint multigraphs\footnote{See Footnote~\ref{Ngraph} for
%the notion of multigraph.} $\Sigma _1,
%\Sigma _2$ over $D(\infty ,R)=\{ z\in \C \mid R\geq |z|\} $ of smooth functions $u_1,u_2\colon
%\{ (r, \theta )\mid r\geq R, \t \in \R \} \to \R $ such that their gradients
%satisfy $\nabla u_i(r,\theta )\to 0$ as $r\to \infty $
%and the separation function $w(r,\theta )=u_1(r,\theta )-u_2(r,\theta )$
%between both multigraphs converges to $\pi $ as $r+|\theta |\to \infty $.
%\item The flux vector of $E_{a,b}$ along its boundary equals $(a,0,b)$,
%and for different values of $a,b$ the related surfaces $E_{a,b}$ are not
%asymptotic after a rigid motion and homothety.
%\item Every complete, embedded minimal annulus in $\R^3$ with compact boundary
%and infinite total curvature is asymptotic
%(up to a rigid motion and homothety) to the surface $E_{a,b}$ as before
%with the related horizontal and vertical fluxes.
%\een
%\end{theorem}

We will finish this section by reporting on
the asymptotic behavior of any properly
embedded minimal surface $M$ in $\R^3$  with
finite genus and an {\it infinite} number of ends.
We have already seen that such an $M$ has exactly two limit ends; it easily follows
that each of its middle ends is asymptotic to a horizontal plane,
after a fixed rotation of $M$. Also, after another rotation of $M$
around a vertical axis followed by a homothety, we can assume that
the flux vector associated to its limit ends has the form $F=(
h,0,1)$, where $h>0$. In this case $M$ is seen to be
conformally a compact Riemann surface  $\overline{M}$ punctured in a
closed countable set with exactly two limit points corresponding to
the two limit ends of $M$, and we can also describe the asymptotic behavior of the ends of $M$:

\begin{theorem}[Asymptotic Limit End Property, Meeks, P\'erez and Ros~\cite{mpr6}]
\label{asympthm}
Let $M$ be a properly embedded minimal surface in $\rth$ with finite genus $g$ and
an infinite number of ends. Then, after a possible rotation and a homothety, the
following statements hold.
\begin{enumerate}
\item  $M$ has two limit ends. In fact, $M$ is conformally
diffeomorphic to  $\overline{M} -{\cal E}_M$, where  $\overline{M}$ is a compact
Riemann surface of genus $g$  and ${\cal E}_M=\{e_n \mid n\in \Z\}
\cup\{e_B,\,e_T\}$ is a countable closed subset of $\overline{M}$
with exactly two limit points $e_T$ and $e_B$. Furthermore,
$\lim_{n\to -\infty}e_n =e_B$, $ \lim_{n\to\infty}e_n= e_T$, and
$e_T$ (resp. $e_B$) corresponds to the top (resp. bottom) end of $M$, while
every $e_n$ with $n\in \Z$ corresponds to a middle end.
\item For each $n \in \Z$, there exists a punctured disk neighborhood $E_n
\subset M\subset \overline{M}$ of $e_n$ which is asymptotic in $\rth$ to a
horizontal plane $P_n$ and which is a graph over its projection to $P_n$.
Furthermore, the usual linear ordering on the index set $\Z$ respects the linear
ordering of the heights of the related planes.
The ordered set of heights $H=\{ h_n
=x_3(P_n)\mid n\in \Z\}$ of these planes naturally corresponds to the set of heights
of the middle ends of $M$.
\item There exists a positive constant $C_M$ such that if $|t|>C_M$, then the
horizontal plane $\{x_3 = t\}$ intersects $M$ in a proper arc when $t\in H$,
or otherwise, $\{x_3 = t\}$ intersects $M$ in a simple closed curve.
\item Let ${\eta}$ denote the unitary outward conormal along the boundary of
$M_t=M\cap \{x_3\leq t\}$. Then the flux vector of $M$, which is defined to be
\[
F_M =\int_{\partial M_t} \eta \, ds
\]
{\rm (}here $ds$ stands for the length element\,{\rm ),} is independent of the
choice of $t$ and has the form $F_M= (h,0,1)$, for some $h>0$.
\item Let $R_h\subset \R^3$ be the Riemann minimal example with horizontal
tangent plane at infinity and flux vector $F=(h,0,1)$ along a compact horizontal
section. Then, there exists a translation vector $v_T\in \rth$ such that as
$t\to \infty $, the function $d_{+}(t)=\sup \left\{
\mbox{\rm dist}(p,{\cal R}_h+v_T) \ | \ p\in M\cap \{ x_3\geq t\} \right\} $
is finite and decays exponentially to zero. In a similar manner, there exists
$v_B\in \rth$ such that as $t\to -\infty$, the function $d_{-}(t)=\sup \left\{
\mbox{\rm dist}(p,{\cal R}_h+v_B) \ | \ p\in M\cap \{ x_3\leq t\} \right\} $ is
finite and decays exponentially to zero. Furthermore, $x_2(v_T) = x_2(v_B).$
\end{enumerate}
\end{theorem}
We now briefly outline the interesting elements of the proof
of Theorem~\ref{asympthm}. The first important ingredient is that  $M$ has
bounded curvature and related uniform local area bounds in balls of
fixed radius. Using these estimates and fixing an end
representative $E_T$ of the top end of $M$, we can find a sequence
of points $p_{n}\in M$ diverging on $E_T$ with Gaussian curvature
bounded by above and away from zero. Our previous arguments show that a
subsequence of the translated surfaces $M(n)=M-p_n$ converges as $n\to \infty $
to a translation of the Riemann minimal example $R_h$ whose flux vector
is equal to $F_M=(h, 0,1)$. Assume after a fixed translation of $M$
that the $M(n)$ converge to $R_h$. Similarly, for a bottom
end representative $E_B$, we can find a divergent sequence of points
$q_n$ such that the surfaces $M'(n)=M-q_n$ converge to $R_h$ as $n\to \infty $.

The second key ingredient in the proof of the theorem is to use the non-zero
Shiffman function $S_{E_{T}}$ of the end representative $E_T$,
to prove that $M(n)$ converges exponentially quickly to $R_h$ as a function of the
$x_3$-coordinate of the points $p_n$ (we are assuming here that the
genus of $M$ is not zero, hence $S_{E_T}$ is not identically zero
because $M$ is not a Riemann minimal example). This is done by first showing that the
norm $|S_{E_{T}}|$ decays exponentially in terms of the
$x_3$-coordinate function of $E_T$. In~\cite{mpr6} we prove that the
bounded Jacobi functions on $R_h$ are all linear and we then
use this property to prove the related exponential decay of
$|S_{E_{T}}|$. Once one has this exponential decay estimate for
$|S_{E_{T}}|$, then elliptic theory can be used to show that $E_{T}$
converges exponentially quickly in terms of $x_3$-coordinates to
$R_h+v_T$ for some vector $v_T\in \rth$. Similarly, $E_{B}$
converges exponentially to $R_h+v_{B}$ for some $v_{B}\in \R^3$.

Finally, a forces argument using the Divergence Theorem applied to a
certain Killing field on $\rth$ shows that the $x_2$-coordinate of
$v_T$ equals the $x_2$-coordinate of $v_B$, which completes our
discussion on the proof of Theorem~\ref{asympthm}.
\par
\vspace{.2cm}
In conclusion, we remark on an important question for a possible
generalization of the classification results discussed in this
paper. This question asks whether or not  a complete embedded
minimal surface $M\subset\R^3$ of finite genus and compact boundary is always properly
embedded. If this is the case, then we can replace in the
statements of our theorems the phrase ``finite genus, properly
embedded minimal surfaces'' by the phrase ``finite genus, complete,
embedded minimal surfaces''. As mentioned above, Colding and Minicozzi~\cite{cm35}
have recently proved that this generalization holds when the surface has finite topology.
Shortly afterwards, Meeks, P\'{e}rez and Ros \cite{mpr9} proved this
result when the surface has finite genus and a countable number of ends; their result
depends on the previously mentioned work in~\cite{cm35} and a related
generalization by Meeks and Rosenberg~\cite{mr13} of the results in \cite{cm35}.
In relation to this problem, recall that Theorem
\ref{thmno1limitend} implies any properly embedded minimal surface
has a countable number of ends. It remains open the outstanding question
of whether or not a complete, embedded minimal surface of finite genus
can have an uncountable number of ends.

\center{William H. Meeks, III at bill@math.umass.edu\\
Mathematics Department, University of Massachusetts, Amherst, MA 01003}
\center{Joaqu\'\i n P\'{e}rez at jperez@ugr.es\\
Department of Geometry and Topology, University of Granada, Granada, Spain}

\bibliographystyle{plain}

\bibliography{bill}

\end{document}